\documentclass[a4paper,12pt]{article}
\usepackage{amssymb}
\usepackage{amsthm}
\usepackage[dvips]{graphicx}
\usepackage{amsmath}
\usepackage{fancyhdr}
\usepackage[all,cmtip]{xy}
\usepackage{stmaryrd}

\usepackage{ascmac}

\theoremstyle{plain}
\newtheorem{theorem}{Theorem}[section]
\newtheorem{proposition}[theorem]{Proposition}
\newtheorem{corollary}[theorem]{Corollary}
\theoremstyle{remark}

\newtheorem{remark}[theorem]{\bf Remark}
\newtheorem{example}[theorem]{\bf Example}
\newtheorem*{acknowledgements}{\bf Acknowledgements}

\renewcommand{\baselinestretch}{1.08}

\newenvironment{reference}[1]{%
\renewcommand{\labelenumi}{[\arabic{enumi}]}
\begin{flushleft}\normalsize{\textbf{References}}\end{flushleft}%
\begin{enumerate}\setlength{\itemsep}{-5pt}\small
}{\end{enumerate}}

\makeatletter
\@addtoreset{equation}{section}
\makeatother

\makeatletter
\renewcommand{\section}{%
\@startsection{section}{1}{\z@}%
{3.5ex \@plus -1ex \@minus -.2ex}%
{2.3ex \@plus.2ex}%
{\reset@font\normalsize\bfseries}}
\renewcommand{\thesection}{\arabic{section}}
\makeatother

\makeatletter
\renewcommand{\subsection}{%
\@startsection{subsection}{2}{\z@}%
{-3.5ex \@plus -1ex \@minus -.2ex}%
{-2.3ex \@plus.2ex}%
{\reset@font\normalsize\bfseries}}
\makeatother

\usepackage[OT2,T1]{fontenc}
\DeclareSymbolFont{cyrletters}{OT2}{wncyr}{m}{n}
\DeclareMathSymbol{\Sha}{\mathalpha}{cyrletters}{"58}

\title{\vspace*{-22mm}
\large{\textbf{
Trivial Iwasawa modules and arithmetic linking numbers
\\
}}
}

\author{
\textsc{\normalsize Yasushi Mizusawa}
\and
\textsc{\normalsize Keiji Okano}
}

\date{}

\begin{document}
{\renewcommand{\baselinestretch}{1.05} \maketitle}

\vspace*{-11mm}
\renewcommand{\abstractname}{}
{\renewcommand{\baselinestretch}{1.05}
\begin{abstract}{\small 
\noindent\textbf{Abstract.} 
We give a computable criterion of the triviality of the unramified Iwasawa module over the maximal multiple $\mathbb Z_p$-extension of a number filed in which an odd prime number $p$ splits completely. This criterion is obtained from the description of the presentation matrix of the Galois coinvariant of the Iwasawa module by arithmetic linking numbers. The description is based on the Koch type presentation of a pro-$p$ Galois group with restricted ramification, which is analogous to the Milnor presentation of a link group from a viewpoint of arithmetic topology. 
}\end{abstract}}

\section{Introduction}

Let $p$ be a prime number, and $k$ a finite extension of the rational number field $\mathbb Q$. 
Let $k_{\mathrm{cyc}}$ be the cyclotomic $\mathbb Z_p$-extension of $k$. 
The Galois group of a $\mathbb Z_p$-extension is isomorphic to the additive group of the ring $\mathbb Z_p$ of $p$-adic integers. 
The maximal multiple $\mathbb Z_p$-extension $\widetilde{k}$ of $k$ is 
the composite of all $\mathbb Z_p$-extensions of $k$, i.e., the maximal free abelian pro-$p$-extension of $k$. 
Then $\mathrm{Gal}(\widetilde{k}/k) \simeq \mathbb Z_p^d$ with $d \ge r_2+1$, 
where $r_2$ is the number of imaginary places of $k$, 
and the equality $d=r_2+1$ is the statement of Leopoldt's conjecture. 
The unramified Iwasawa module $X(\widetilde{k})=\mathrm{Gal}(L(\widetilde{k})/\widetilde{k})$ over $\widetilde{k}$ is the Galois group 
of the maximal unramified abelian pro-$p$-extension $L(\widetilde{k})$ of $\widetilde{k}$ 
with the action of $\widetilde{\varGamma}=\mathrm{Gal}(\widetilde{k}/k)$. 
Greenberg's genelarized conjecture (\cite{Gre01}) asserts that $X(\widetilde{k})$ is pseudonull as a module over the complete group ring $\widetilde{\varLambda} = \mathbb Z_p\llbracket \widetilde{\varGamma} \rrbracket$. 
This conjecture has been verified in various situations 
(\cite{Fuj17,Ito11,Kat17,Min86,Sha08} etc.), 
and the structure (e.g., triviality, $\widetilde{\varLambda}$-cyclicity, and pseudonull $\widetilde{\varLambda}$-submodules) of $X(\widetilde{k})$ has also been studied (\cite{Fuj10,Kle22,MMOO21,MMOO22,Mur23} etc.). 
The case of trivial $X(\widetilde{k})$ is a special case of pseudonull $X(\widetilde{k})$. 
If $p \neq 2$ and $k$ is an abelian $p$-extension of $\mathbb Q$, 
Leopoldt's conjecture holds, and $k$ is totally real, i.e., $d=1$. 
Then $\widetilde{k}=k_{\mathrm{cyc}}$, 
and the cases with trivial $X(k_{\mathrm{cyc}})$ were classified in \cite{Gen00}. 
Additionally in \cite{MY20}, the classification condition was described by arithmetic linking numbers from a viewpoint of arithmetic topology (\cite{Miz19,Mor02,Mor}). 
If $p \neq 2$ and $k$ is a CM-field in which $p$ splits completely, 
the cases with trivial $X(\widetilde{k})$ have been classified in \cite{Oka**}. 
The aim of this paper is to describe the triviality of $X(\widetilde{k})$ by arithmetic linking numbers. 

Let $\widetilde{L}(k_{\mathrm{cyc}})$ be the maximal unramified pro-$p$-extension of $k_{\mathrm{cyc}}$. 
Then $\widetilde{L}(k_{\mathrm{cyc}})$ is a pro-$p$-extension of $k$. 
Put the pro-$p$ Galois group 
\[
\widetilde{G}= \mathrm{Gal}(\widetilde{L}(k_{\mathrm{cyc}})/k).
\]
The Iwasawa module 
$X(k_{\mathrm{cyc}}) = \mathrm{Gal}(L(k_{\mathrm{cyc}})/k_{\mathrm{cyc}})$ 
over $k_{\mathrm{cyc}}$ 
is the abelianization of the subgroup $\widetilde{H}= \mathrm{Gal}(\widetilde{L}(k_{\mathrm{cyc}})/k_{\mathrm{cyc}})$ of $\widetilde{G}$. 
The commutativity of $\widetilde{H}$ has been studied in nonabelian Iwasawa theory of $\mathbb Z_p$-extensions (\cite{Oka06,Oka12I,Oza07,Sha08} etc.). 
Now we assume that $p$ splits completely in $k$. 
It is known that $\widetilde{k}/k_{\mathrm{cyc}}$ is unramified, i.e., $\widetilde{k} \subset \widetilde{L}(k_{\mathrm{cyc}})$. 
Then $\widetilde{L}(k_{\mathrm{cyc}})$ is also the maximal unramified pro-$p$-extension of $\widetilde{k}$, 
and there is a surjective homomorphism $\widetilde{G}/[\widetilde{G},\widetilde{G}] \twoheadrightarrow \mathbb Z_p^d$ with finite kernel, 
where $[\widetilde{G},\widetilde{G}]$ denotes the closed commutator subgroup of $\widetilde{G}$.
The Iwasawa module $X(\widetilde{k})$ over $\widetilde{k}$ can be treated as the abelianization of the subgroup $\mathrm{Gal}(\widetilde{L}(k_{\mathrm{cyc}})/\widetilde{k})$ of $\widetilde{G}$, 
and the triviality of $X(\widetilde{k})$ is a sufficient condition of the commutativity of $\widetilde{H}$. 
Assuming $\widetilde{G}/[\widetilde{G},\widetilde{G}] \simeq \mathbb Z_p^{r_2+1}$, which is a necessary condition of the triviality of $X(\widetilde{k})$ under Leopoldt's conjecture, 
we will recall in Theorem \ref{thm:Ghyoji} that $\widetilde{G}$ has a Koch type presentation of the form 
\[
\widetilde{G} \simeq \langle x_{1-t},x_{2-t},\ldots,x_{r_2+1-t} \,|\, [x_1,y_1], \ldots, [x_{r_1+2r_2},y_{r_1+2r_2}] \rangle 
\]
as a pro-$p$ group (\cite{EM22,Sal08}), 
which is analogous to the Milnor presentation of a link group. 
The $\widetilde{\varGamma}$-coinvariant $X(\widetilde{k})_{\widetilde{\varGamma}}$ of $X(\widetilde{k})$ is the maximal quotient of $X(\widetilde{k})$ on which $\widetilde{\varGamma}$ acts trivially, 
and it can be treated as a subquotient of $\widetilde{G}$. 
Based on the Koch type presentation of $\widetilde{G}$, 
we obtain the following main theorem which shows that a presentation matrix of $X(\widetilde{k})_{\widetilde{\varGamma}}$ is described by arithmetic linking numbers. 


\begin{theorem}\label{thm:main}
Let $r_2$ (resp.\ $r_1$) be the number of imaginary (resp.\ real) places of a 
number field $k$. 
Assume that $p \neq 2$ and that $p$ splits completely in $k$.  
Moreover, we assume that $\widetilde{G}/[\widetilde{G},\widetilde{G}] \simeq \mathbb Z_p^{r_2+1}$. 
(This assumption includes the validity of Leopoldt's conjecture.)
Then there is a $(r_1+2r_2) \times \frac{r_2(r_2+1)}{2}$ presentation matrix $A=(a_{j,(i,i')})_{j,(i,i')}$ of the $\mathbb Z_p$-module $X(\widetilde{k})_{\widetilde{\varGamma}}$ (see \eqref{eq:ajii} and \eqref{eq:A} for definition) such that: 
\begin{enumerate}
\setlength{\itemsep}{0pt}
\renewcommand{\labelenumi}{\rm (\roman{enumi})}
\renewcommand{\theenumi}{\rm (\roman{enumi})}
\item\label{main:a_is_lk}
the entries of $A$ are rational expressions of `arithmetic linking numbers' $l_{j,j'}$, $z_{i,j}$ (see \eqref{eq:zlk} for definition), 
which are values of $p$-adic logarithms of some specific computable elements $\varpi_j$, $\varepsilon_i$ of $k$, 
\item\label{main:lmdrk}
the minimal number of generators of $X(\widetilde{k})$ as a $\widetilde{\varLambda}$-module is $\frac{r_2(r_2+1)}{2}-\mathrm{rank}\,\overline{A}$, 
where $\overline{A}=(\overline{a_{j,(i,i')}})_{j,(i,i')}$ is a matrix with $\mathbb F_p$-entries $\overline{a_{j,(i,i')}}=a_{j,(i,i')}+p\mathbb Z_p$. 
\end{enumerate}
In particular, $\mathrm{rank}\,\overline{A} = \frac{r_2(r_2+1)}{2}$ if and only if $X(\widetilde{k})$ is trivial. 
\end{theorem}

For example, if $k$ is an imaginary bicyclic biquadratic field, 
we obtain the following statement. 
Let $\mathcal{O}_k$ be the ring of algebraic integers of $k$. 

\begin{corollary}[a part of Proposition \ref{prop:4t0}]\label{sttmt:V4t0}
Assume that $p \neq 2$ and $k$ is an imaginary bicyclic biquadratic field in which $p$ splits completely, and that the class number $h$ of $k$ is not divisible by $p$. 
Let $\mathfrak{p}_1,\mathfrak{p}_2,\mathfrak{p}_3,\mathfrak{p}_4$ be the prime ideals of $k$ lying over $p$, 
and $\varpi_1 \in \mathcal{O}_k$ a generator of the principal ideal $\mathfrak{p}_1^h=\varpi_1\mathcal{O}_k$. 
Suppose that $\mathfrak{p}_3,\mathfrak{p}_4$ are the complex conjugations of $\mathfrak{p}_1,\mathfrak{p}_2$ respectively. 
Assume that the fundamental unit $\varepsilon_1$ of $k$ satisfies 
$\varepsilon_1^{p-1} \not\equiv 1 \pmod{\mathfrak{p}_1^2}$. 
Then the following statements hold true: 
\begin{enumerate}
\setlength{\itemsep}{0pt}
\renewcommand{\labelenumi}{\rm (\roman{enumi})}
\renewcommand{\theenumi}{\rm (\roman{enumi})}
\item
$\widetilde{G}/[\widetilde{G},\widetilde{G}] \simeq \mathbb Z_p^{3}$. 
\item
There is an exact sequence 
\[
\entrymodifiers={+!!<0pt,\fontdimen22\textfont2>}
\xymatrix@R=0.8\baselineskip{
\mathbb Z_p^4 \ar[r]^-{A} & \mathbb Z_p^3 \ar[r] & X(\widetilde{k})_{\widetilde{\varGamma}} \ar[r] & 0  
} 
\]
with a presentation matrix $A$ of the form
\begin{equation}\label{eq:V4t0A}
A=h^{-1}\begin{pmatrix}
 l_{1,2}-l_{1,4} & l_{1,3}+l_{1,4} & 0 \\
-l_{1,2}-l_{1,3} & 0 & l_{1,3}+l_{1,4} \\
0 & -l_{1,2}-l_{1,3} & l_{1,2}-l_{1,4} \\
l_{1,3}+l_{1,4} &  l_{1,2}-l_{1,4} & -l_{1,2}-l_{1,3}
\end{pmatrix}, 
\end{equation}
where $\mathfrak{p}_j$-adic logarithm $\log_{\mathfrak{p}_j}$ defines $l_{1,j}=\log_{\mathfrak{p}_j}(\varpi_1)/\log_{\mathfrak{p}_j}(1+p)$ for each $j \in \{2,3,4\}$. 
\item
$X(\widetilde{k})$ is trivial if and only if $(l_{1,2}-l_{1,4})(l_{1,2}+l_{1,3})(l_{1,3}+l_{1,4}) \not\equiv 0 \pmod{p}$. 
\item\label{V4t0:minlam}
$X(k_{\mathrm{cyc}}) \simeq \mathbb Z_p^2$ if and only if $(l_{1,2}+l_{1,3})(l_{1,3}+l_{1,4}) \not\equiv 0 \pmod{p}$. 
\end{enumerate}
\end{corollary}

Theorem \ref{thm:main} is proved in Section \ref{sec:lk}. 
In \S\ref{subsec:settei}, we define arithmetic linking numbers \eqref{eq:zlk}, and give an equivariant condition of the assumption that $\widetilde{G}/[\widetilde{G},\widetilde{G}] \simeq \mathbb Z_p^{r_2+1}$ (Proposition \ref{prop:Vis0}). 
In \S\ref{subsec:Ghyoji}, under the assumption, 
a Koch type presentation of $\widetilde{G}$ is described by arithmetic linking numbers  (Theorem \ref{thm:Ghyoji}). 
In \S\ref{subsec:mainpf}, we prove Theorem \ref{thm:main} with the definition \eqref{eq:A} of the presentation matrix $A$ of $X(\widetilde{k})_{\widetilde{\varGamma}}$ , 
and recall that $[k:\mathbb Q] \le 6$ if $k$ is a CM-field and $X(\widetilde{k})$ is trivial (Corollary \ref{cor:r2<3}). 
Moreover in \S\ref{subsec:minlam}, 
assuming that $p$ does not divide the class number of $k$, 
we also give a criterion of the minimality of the Iwasawa $\lambda$-invariant of $k_{\mathrm{cyc}}/k$, which is the $\mathbb Z_p$-rank of $X(k_{\mathrm{cyc}})$ (Theorem \ref{thm:lambdar2}). 
This criterion (e.g., Corollary \ref{sttmt:V4t0} \ref{V4t0:minlam}) is a generalization of Gold's theorem \cite{Gol74} (Proposition \ref{prop:Gold}) on imaginary quadratic fields. 
In Section \ref{sec:X0}, based on Theorem \ref{thm:main}, 
we give detail criteria of the triviality of $X(\widetilde{k})$ for CM-fields $k$ according to $[k:\mathbb Q] \in \{2,4,6\}$. 
Corollary \ref{sttmt:V4t0} will be proved together with Proposition \ref{prop:4t0} for quartic CM-fields $k$. 
Moreover, we describe the invariants $a,b,c$ of \cite{Oka**} by arithmetic linking numbers (Proposition \ref{prop:abc}). 
In Section \ref{sec:examples}, we give various examples of CM-fields $k$ with trivial (or nontrivial, $\widetilde{\varLambda}$-cyclic) $X(\widetilde{k})$.

\section{Arithmetic linking numbers}\label{sec:lk}

\subsection{Settings and assumptions.}\label{subsec:settei}
Suppose $p \neq 2$. 
Let $k$ be a number field in which $p$ splits completely. 
Let $r_1$ (resp.\ $r_2$) be the number of real (resp.\ imaginary) places of $k$. 
Then the degree of $k$ is $[k:\mathbb Q]=r_1+2r_2$. 
Let $P=\{\mathfrak{p}_1,\mathfrak{p}_2,\ldots,\mathfrak{p}_{r_1+2r_2}\}$ be the set of all prime ideals $\mathfrak{p}_1,\mathfrak{p}_2,\ldots,\mathfrak{p}_{r_1+2r_2}$ of $k$ lying over $p$. 
Since $p$ splits completely in $k/\mathbb Q$, we have $\zeta_p \not\in k$. 
If $k$ is a CM-field, $E_k^{p-1} \subset E_{k^+}$, where $E_k$ and $E_{k^+}$ are the unit groups of $k$ and the maximal real subfield $k^+$ respectively. 
The $p$-rank of $E_k$ is $\dim_{\mathbb F_p}E_k/p=r_1+r_2-1$. 
Let $\{\varepsilon_1,\varepsilon_2,\ldots,\varepsilon_{r_1+r_2-1}\}$ be a system of units of $k$ 
such that 
\[
E_k/p=\langle \varepsilon_1E_k^p,\varepsilon_2E_k^p,\ldots,\varepsilon_{r_1+r_2-1}E_k^p \rangle.
\] 
If $k$ is a CM-field, we may choose the fundamental units of $k^+$ as $\varepsilon_1,\varepsilon_2,\ldots,\varepsilon_{r_2-1}$. 
Let $Cl_k$ be the ideal class group of $k$, 
and put $t=\dim_{\mathbb F_p} Cl_k/p$, the $p$-rank of $Cl_k$. 
Let $h$ be the non-$p$-part of the class number $|Cl_k|$ of $k$. 
If $t=0$, then $h=|Cl_k|$. 
If $t \ge 1$, we choose $t$ prime ideals 
$\mathfrak{a}_{1-t},\ldots,\mathfrak{a}_{-1},\mathfrak{a}_0$ of $k$ such that 
\[
Cl_k/p=\langle [\mathfrak{a}_{1-t}]Cl_k^p,\ldots,[\mathfrak{a}_{-1}]Cl_k^p,[\mathfrak{a}_0]Cl_k^p \rangle
\]
and $\{\mathfrak{a}_{1-t},\ldots,\mathfrak{a}_{-1},\mathfrak{a}_0\} \cap P =\emptyset$. 
Then we also choose $t$ integers $\varepsilon_{1-t},\ldots,\varepsilon_{-1},\varepsilon_0$ of $k$ such that 
\[
\mathfrak{a}_i^{q_i} = \varepsilon_i \mathcal{O}_k 
\]
for each $i \in \{1-t,\ldots,-1,0\}$, 
where $q_i$ is the order of $[\mathfrak{a}_i] \in Cl_k$. 
Then the Sylow $p$-subgroup of $Cl_k$ is a subgroup of $\langle [\mathfrak{a}_{1-t}],\ldots,[\mathfrak{a}_{-1}],[\mathfrak{a}_0] \rangle$, 
and hence $[\mathfrak{p}_i]^h = \prod_{j=1-t}^{0} [\mathfrak{a}_j]^{l_{i,j}}$ with some $0 \le l_{i,j} \in \mathbb Z$ for each $\mathfrak{p}_i \in P$. 
For each $\mathfrak{p}_i \in P$, 
we choose $\varpi_i \in k^{\times}$ such that 
\[
\mathfrak{p}_i^h=
\bigg\{\begin{array}{ll}
\varpi_i \mathcal{O}_k & \text{ if } t=0, \\
(\prod_{j=1-t}^{0} \mathfrak{a}_j^{l_{i,j}}) \varpi_i \mathcal{O}_k & \text{ if } t \ge 1. 
\end{array}
\]
For each prime $v$ of $k$, 
we denote by $k_v$ the localization of $k$ at $v$, 
and by $U_v$ the local unit group of $k_v$. 
Let $J$ be the id\`ele group of $k$. 
Since $k_{\mathfrak{p}} \simeq \mathbb Q_p$ for $\mathfrak{p} \in P$, 
we have $U_{\mathfrak{p}}=W_{\mathfrak{p}} \times (1+p)^{\mathbb Z_p}$, 
where $W_{\mathfrak{p}}$ is the group of $(p-1)$th roots of unity in $k_{\mathfrak{p}}$. 
Put $\widetilde{U} = \prod_{v \not\in P} U_v \times \prod_{\mathfrak{p} \in P} W_{\mathfrak{p}} \subset J$, 
$\widetilde{V} = k^{\times} \cap \widetilde{U}J^p = k^{\times} \cap (\prod_{v \not\in P} U_v)J^p$, 
and $V = k^{\times} \cap (\prod_{v} U_v)J^p$. 
Then there is an exact sequence 
\[
\entrymodifiers={+!!<0pt,\fontdimen22\textfont2>}
\xymatrix@R=0.8\baselineskip{
0 \ar[r] & E_k/p \ar[r] & V/(k^{\times})^p \ar[r] & Cl_k \ar[r]^-{x \mapsto x^p} & Cl_k , 
} 
\]
where the image of $a(k^{\times})^p \in V/(k^{\times})^p$ is well-defined as 
the class $[\mathfrak{a}] \in Cl_k$ of $\mathfrak{a}$ satisfying $\mathfrak{a}^p=a\mathcal{O}_k$ (cf.\ \cite[\S11.2]{Koch}). 
If $t \ge 1$, the kernel of $Cl_k \stackrel{x \mapsto x^p}{\longrightarrow} Cl_k$ is generated by $[\mathfrak{a}_{1-t}^{q_{1-t}/p}],\ldots,[\mathfrak{a}_{-1}^{q_{-1}/p}],[\mathfrak{a}_0^{q_0/p}]$. 
Hence 
$\dim_{\mathbb F_p} V/(k^{\times})^p = r_1+r_2-1+t$ and 
\[
V/(k^{\times})^p = \langle \varepsilon_{1-t}(k^{\times})^p,\varepsilon_{2-t}(k^{\times})^p,\ldots,\varepsilon_{r_1+r_2-1}(k^{\times})^p \rangle.
\] 
For each 
$i \in \{1-t,2-t,\ldots,r_1+r_2-1\}$ and $\mathfrak{p}_j,\mathfrak{p}_{j'} \in P$ ($j \neq j'$), 
put 
\begin{align}\label{eq:zlk}
z_{i,j} 
= \frac{\log_{\mathfrak{p}_j}(\varepsilon_i)}{\log_{\mathfrak{p}_j}(1+p)} \in \mathbb Z_p
\quad \text{and} \quad 
l_{j',j} 
= \frac{\log_{\mathfrak{p}_j}(\varpi_{j'})}{\log_{\mathfrak{p}_j}(1+p)} \in \mathbb Z_p, 
\end{align}
where $\log_{\mathfrak{p}_j}$ denotes the $\mathfrak{p}_j$-adic logarithm. 
Then 
\[
\varepsilon_i (1+p)^{-z_{i,j}} \in W_{\mathfrak{p}_j} 
\quad \text{and} \quad 
\varpi_{j'} (1+p)^{-l_{j',j}} \in W_{\mathfrak{p}_j}. 
\]
In particular, 
\[
\varepsilon_i^{p-1} \equiv  1-z_{i,j}p \pmod{\mathfrak{p}_j^2} 
\quad \text{and} \quad 
\varpi_{j'}^{p-1} \equiv  1-l_{j',j}p \pmod{\mathfrak{p}_j^2}. 
\]
Put $l_{j',j}=0$ if $j=j'$. 
Put $(r_1+r_2-1+t) \times (r_1+2r_2)$ matrices 
\[
Z=(z_{i,j})_{1-t \le i \le r_1+r_2-1 \atop 1 \le j \le r_1+2r_2}
=(\begin{matrix}
\boldsymbol{z}_1 & \boldsymbol{z}_2 & \cdots & \boldsymbol{z}_{r_1+2r_2} 
\end{matrix})
\]
and 
\[
\overline{Z}=(\overline{z_{i,j}})_{1-t \le i \le r_1+r_2-1 \atop 1 \le j \le r_1+2r_2}, 
\]
where $\overline{z_{i,j}}=z_{i,j}+p\mathbb Z_p \in \mathbb F_p$. 
By \cite[Th\'eor\`eme 3.3]{Sal08} (cf.\ \cite[Theorem 1.5 (i)]{BLM13}, or \cite[Theorem 2.2]{EM22}), 
the generator rank 
$d(\widetilde{G})=\dim_{\mathbb F_p} \widetilde{G}/\widetilde{G}^p[\widetilde{G},\widetilde{G}]$ 
of $\widetilde{G}$ is 
\begin{align}\label{eq:genrank}
d(\widetilde{G})= \dim_{\mathbb F_p} \widetilde{V}/(k^{\times})^p +r_2+1.
\end{align}

\begin{proposition}\label{prop:Vis0}
The following three conditions are equivalent:
\begin{enumerate}
\setlength{\itemsep}{0pt}
\renewcommand{\labelenumi}{\rm (\arabic{enumi})}
\renewcommand{\theenumi}{\rm (\arabic{enumi})}
\item\label{genrk}
$d(\widetilde{G})=r_2+1$, i.e., $\widetilde{G}/[\widetilde{G},\widetilde{G}] \simeq \mathbb Z_p^{r_2+1}$. 
\item\label{trivV}
$\widetilde{V}/(k^{\times})^p \simeq 0$. 
\item\label{rkB}
$\mathrm{rank}\,\overline{Z}=r_1+r_2-1+t$, or $r_1+r_2=1$ and $t=0$. 
\end{enumerate}
\end{proposition}

\begin{proof}
By \eqref{eq:genrank}, \ref{genrk} is equivalent to \ref{trivV}. 
If $r_1+r_2=1$ and $t=0$, 
then $V/(k^{\times})^p \simeq 0$ and hence $\ref{trivV}$ holds true. 
Suppose that $r_1+r_2>1$ or $t \ge 1$. 
Put $\varepsilon_{\overline{\mathbf{e}}} = \prod_{i=1-t}^{r_1+r_2-1} \varepsilon_i^{e_i} \in V$ 
for $\overline{\mathbf{e}}=(\overline{e_{1-t}},\overline{e_{2-t}},\ldots,\overline{e_{r_1+r_2-1}}) \in \mathbb F_p^{r_1+r_2-1+t}$ 
with $e_i \in \{0,1,\ldots,p-1\}$. 
Note that $\varepsilon_{\overline{\mathbf{e}}} \in (k^{\times})^p$ if and only if $\overline{\mathbf{e}}=(\overline{0},\overline{0},\ldots,\overline{0})$. 
Since 
\[
\varepsilon_{\overline{\mathbf{e}}}^{-1} (1+p)^{\sum_{i=1-t}^{r_1+r_2-1} e_i z_{i,j}} = \prod_{i=1-t}^{r_1+r_2-1}  (\varepsilon_i (1+p)^{-z_{i,j}})^{-e_i} \in W_{\mathfrak{p}_j} \subset (k_{\mathfrak{p}_j}^{\times})^p
\]
for each $\mathfrak{p}_j \in P$, 
$\varepsilon_{\overline{\mathbf{e}}} \in \widetilde{V}$ if and only if $\sum_{i=1-t}^{r_1+r_2-1} e_i z_{i,j} \equiv 0 \pmod{p}$ for all $\mathfrak{p}_j \in P$, i.e., 
$\overline{\mathbf{e}} \overline{Z} = \boldsymbol{0} \in \mathbb F_p^{r_1+2r_2}$. 
Therefore $\varepsilon_{\overline{\mathbf{e}}} \in \widetilde{V} \setminus (k^{\times})^p$ 
if and only if $\overline{\mathbf{e}} \neq (\overline{0},\overline{0},\ldots,\overline{0})$ and 
$\overline{\mathbf{e}} \overline{Z} = \boldsymbol{0}$. 
This implies that $\widetilde{V}/(k^{\times})^p \not\simeq 0$ if and only if $\mathrm{rank}\,\overline{Z}<r_1+r_2-1+t$. 
Hence \ref{trivV} is equivalent to \ref{rkB}. 
\end{proof}

\begin{remark}
If $r_2=0$, the condition \ref{genrk} is equivalent to that 
$\widetilde{k}=k_{\mathrm{cyc}}$ and $X(k_{\mathrm{cyc}})$ is trivial. 
If $k$ is a real quadratic field with $t=0$, 
then $Z=(z_{1,1},-z_{1,1})$, and the condition \ref{rkB} is equivalent to that $\varepsilon_1^{p-1} \not\equiv 1 \pmod{\mathfrak{p}_1^2}$. 
When $r_2=0$, under Leopoldt's conjecture, 
Proposition \ref{prop:Vis0} can be regarded as a generalization of a part of Fukuda-Komatsu's criterion \cite{FK86}. 
\end{remark}

Now we assume that the conditions \ref{genrk} \ref{trivV} \ref{rkB} of Proposition \ref{prop:Vis0} are satisfied. 
Then, since $k_{\mathrm{cyc}}/k$ is totally ramified, we have $t < r_2+1$. 
In particular, $r_1+r_2-1+t < r_1+2r_2$, i.e., $Z$ is a horizontal matrix. 
By renumbering the elements of $P$, we may assume that 
the $(r_1+r_2-1+t) \times (r_1+r_2-1+t)$ submatrix
\begin{align*}
B
&=(z_{i,j})_{1-t \le i \le r_1+r_2-1 \atop r_2+2-t \le j \le r_1+2r_2}
= (\begin{matrix}
\boldsymbol{z}_{r_2+2-t} & \boldsymbol{z}_{r_2+3-t} & \cdots & \boldsymbol{z}_{r_1+2r_2} 
\end{matrix})
\\
&= 
\begin{pmatrix}
z_{1-t,r_2+2-t} & z_{1-t,r_2+3-t} & \cdots & z_{1-t,r_1+2r_2} \\
z_{2-t,r_2+2-t} & z_{2-t,r_2+3-t} & \cdots & z_{2-t,r_1+2r_2} \\
\vdots & \vdots & \ddots & \vdots \\
z_{r_1+r_2-1,r_2+2-t} & z_{r_1+r_2-1,r_2+3-t} & \cdots & z_{r_1+r_2-1,r_1+2r_2} 
\end{pmatrix}
\end{align*}
of $Z$ satisfies $\det B \not\equiv 0 \pmod{p}$. 
Let $b_{j,i}$ be the entries of 
\[
B^{-1}=(b_{j,i})_{r_2+2-t \le j \le r_1+2r_2 \atop 1-t \le i \le r_1+r_2-1}
\]
such that 
\[
\sum_{i=1-t}^{r_1+r_2-1} b_{j,i} z_{i,j'} = 
\bigg\{
\begin{array}{ll}
1 & \text{if $j=j'$}, \\
0 & \text{if $j \neq j'$ and $r_2+2-t \le j' \le r_1+2r_2$}. 
\end{array}
\]
For $j \in \{r_2+2-t,\ldots,r_1+2r_2\}$ and $i \in \{1-t,\ldots,r_2+1-t\}$, 
put 
\begin{align*}
c_{j,i}
=-\sum_{\kappa=1-t}^{r_1+r_2-1} b_{j,\kappa} z_{\kappa,i} 
= -
\frac{\left|\begin{matrix}
\boldsymbol{z}_{r_2+2-t} & \cdots & \boldsymbol{z}_{j-1} & \boldsymbol{z}_i & \boldsymbol{z}_{j+1} & \cdots & \boldsymbol{z}_{r_1+2r_2} 
\end{matrix}\right|}
{\left|\begin{matrix}
\boldsymbol{z}_{r_2+2-t} & \cdots & \boldsymbol{z}_{j-1} & \boldsymbol{z}_j & \boldsymbol{z}_{j+1} & \cdots & \boldsymbol{z}_{r_1+2r_2} 
\end{matrix}\right|} 
\end{align*}
if $i \ge 1$, and put $c_{j,i} = -q_i b_{j,i}$ if $i \le 0$. 
Then we define 
\[
\widetilde{l}_{j,i}
=(l_{j,i}+\sum_{\ell=r_2+2-t}^{r_1+2r_2} c_{\ell,i} l_{j,\ell} ) h^{-1} 
\]
for $j \in \{1,2,\ldots,r_1+2r_2\}$ and $i \in \{1-t,\ldots,r_2+1-t\}$.

\subsection{Presentation of Galois groups.}\label{subsec:Ghyoji}
We recall the settings and the assumptions:
\begin{enumerate}
\setlength{\itemsep}{0pt}
\renewcommand{\labelenumi}{\rm (a\arabic{enumi})}
\renewcommand{\theenumi}{\rm (a\arabic{enumi})}
\item\label{asm:podd}
$p \neq 2$. 
\item\label{asm:spcomp}
$p$ splits completelly in $k$ as $p\mathcal{O}_k=\prod_{i=1}^{r_1+2r_2} \mathfrak{p}_i$. 
\item\label{asm:tprk}
$t < r_2+1$, where $t$ is the $p$-rank of $Cl_k$. 
\item\label{asm:B}
$\det B \not\equiv 0 \pmod{p}$ for the rightmost $(r_1+r_2-1+t) \times (r_1+r_2-1+t)$ submatrix $B$ of $Z$. 
Then $\widetilde{V}/(k^{\times})^p \simeq 0$ by Proposition \ref{prop:Vis0}. 
\end{enumerate}

By the id\`elic class field theory, 
there exists a subgroup $Y$ of $J$ such that $\widetilde{G}/[\widetilde{G},\widetilde{G}] \simeq J/Y$. 
Let $f : \widetilde{G} \rightarrow J/Y$ be the natural surjective homomorphism.
Let $\sigma_i \in \widetilde{G}$ be a Frobenius element of a prime $\widetilde{\mathfrak{p}}_i$ lying over $\mathfrak{p}_i \in P$ 
such that 
$f(\sigma_i^h)=\varpi_i Y$.
Let $\alpha_i=(\alpha_{i,v})_v$ be an element of $J$ such that $\alpha_{i,\mathfrak{p}_i}=1+p$ and $\alpha_{i,v}=1$ for all $v \neq \mathfrak{p}_i$. 
There is a generator $\tau_i \in \widetilde{G}$ of an inertia group of $\widetilde{\mathfrak{p}}_i$ 
such that $f(\tau_i)=\alpha_iY$. 
If $t \ge 1$, for each $i \in \{1-t,\ldots,-1,0\}$, 
we put $\alpha_i=(\alpha_{i,v})_v \in J$ such that $\alpha_{i,v}$ is a uniformizer of $k_v$ if $v=\mathfrak{a}_i$, and $\alpha_{i,v}=1$ for all $v \neq \mathfrak{a}_i$. 
Then there are elements $\tau_{1-t},\ldots,\tau_{-1},\tau_0 \in \widetilde{G}$ such that $f(\tau_i)=\alpha_iY$ for each $i \le 0$ and $\mathrm{Gal}(L(k)/k)=\langle \tau_{1-t}|_{L(k)},\ldots,\tau_{-1}|_{L(k)},\tau_0|_{L(k)} \rangle$, 
where $L(k)$ is the maximal abelian $p$-extension of $k$. 

For a pro-$p$ group $G$ and its subgroup $H$, 
let $[H,G]$ be the minimal closed subgroup of $G$ containing $[h,g]=h^{-1}g^{-1}hg$ for all $h \in H$ and all $g \in G$. 

\begin{theorem}\label{thm:Ghyoji}
Under the settings in \S\ref{subsec:settei} and the assumptions 
{\ref{asm:podd}--\ref{asm:B}}
above, 
$\widetilde{G}$ has the following presentation as a finitely generated pro-$p$ group:
\[
\entrymodifiers={+!!<0pt,\fontdimen22\textfont2>}
\xymatrix@R=0.8\baselineskip{
1 \ar[r] & R \ar[r] & F \ar[r]^-{\pi} & \widetilde{G} \ar[r] & 1,  
} 
\]
where $F$ is a free pro-$p$ group on $r_2+1$ letters $x_{1-t},x_{2-t},\ldots,x_{r_2+1-t}$, 
$R$ is a closed normal subgroup of $F$ which is normally generated by 
$r_1+2r_2$ elements $[x_i,y_i]$ $(i \in \{1,2,\ldots,r_1+2r_2\})$, 
where $x_i$ and $y_i$ are pro-$p$ words in $F$ such that $\pi(x_i)=\tau_i$, $\pi(y_i)=\sigma_i$, 
\begin{align}\label{eq:xj}
x_j \equiv 
\prod_{i=1-t}^{r_2+1-t} x_i^{c_{j,i}} \pmod{[F,F]}
\end{align}
if $r_2+2-t \le j \le r_1+2r_2$, and 
\begin{align}\label{eq:yj}
y_j 
\equiv \prod_{i=1-t}^{r_1+2r_2} x_i^{l_{j,i} h^{-1}} 
\equiv \prod_{i=1-t}^{r_2+1-t} x_i^{\widetilde{l}_{j,i}} \pmod{[F,F]} 
\end{align}
for $j \in \{1,2,\ldots,r_1+2r_2\}$. 
\end{theorem}

\begin{proof}
This theorem is a special case of \cite[Theorem 1.1]{EM22}. 
One can easily see that the settings and the assumptions of Theorem \ref{thm:Ghyoji} satisfy them of \cite[\S3 Proof of Theorem 1.1]{EM22}. 
In particular, the congruences \eqref{eq:xj}, \eqref{eq:yj} are the same as \cite[(3.6), (3.7)]{EM22}. 
\end{proof}

\begin{remark}
Recall that $l_{j',j}$ (resp.\ $z_{i,j}$) modulo $p$ describes the $p$-th power residue of $\varpi_{j'}^{p-1}$ (resp.\ $\varepsilon_i^{p-1}$) modulo $\mathfrak{p}_j^2$. 
According to the sense of \cite{Mor} (and \cite{CKKPPY19} etc.), 
we call $z_{i,j}$, $l_{j',j}$ arithmetic linking numbers, 
which form analogues $\widetilde{l}_{j,i}$ of linking numbers in analogy between the presentation $\widetilde{G} \simeq F/R$ and the Milnor presentation of a link group. 
\end{remark}

\subsection{Proof of Theorem \ref{thm:main}}\label{subsec:mainpf}
The assumptions {\ref{asm:podd}--\ref{asm:B}} are satisfied by Proposition \ref{prop:Vis0}. 
By Theorem \ref{thm:Ghyoji}, $\widetilde{G} \simeq F/R$. 
Since $\widetilde{\varGamma} \simeq \widetilde{G}/[\widetilde{G},\widetilde{G}] \simeq F/[F,F]$, 
we have $R \subset [F,F]$ and $[\widetilde{G},\widetilde{G}]=\mathrm{Gal}(\widetilde{L}(k_{\mathrm{cyc}})/\widetilde{k})$. 
Then 
\[
X(\widetilde{k}) \simeq [\widetilde{G},\widetilde{G}]/[[\widetilde{G},\widetilde{G}],[\widetilde{G},\widetilde{G}]] \simeq [F,F]/[[F,F],[F,F]]R , 
\]
and hence 
$X(\widetilde{k})_{\widetilde{\varGamma}} \simeq [F,F]/[[F,F],F]R$. 
Put an lexicographically ordered set 
\[
\mathcal{I}=\{(i,i') \in \mathbb Z^2 \mid 1-t \le i < i' \le r_2+1-t\} .
\] 
Then $|\mathcal{I}|=\frac{r_2(r_2+1)}{2}$, and 
$\{[x_i,x_{i'}] \mid (i,i') \in \mathcal{I} \}$ forms a basis of the free $\mathbb Z_p$-module $[F,F]/[[F,F],F] \simeq \mathbb Z_p^{\frac{r_2(r_2+1)}{2}}$. 
There is an exact sequence 
\[
\entrymodifiers={+!!<0pt,\fontdimen22\textfont2>}
\xymatrix@R=0.8\baselineskip{
0 \ar[r] & R[[F,F],F]/[[F,F],F] \ar[r]^-{f_A} & [F,F]/[[F,F],F] \ar[r] & X(\widetilde{k})_{\widetilde{\varGamma}} \ar[r] & 0 \\
& \mathbb Z_p^{r_1+2r_2} \ar@{->>}[u] \ar[r]^-{A} & \mathbb Z_p^{\frac{r_2(r_2+1)}{2}} \ar@{}[u]|-{\rotatebox[origin=c]{90}{$\simeq$}} & & 
} 
\]
with a presentation matrix $A$ of $X(\widetilde{k})_{\widetilde{\varGamma}}$. 
If $j \in \{1,2,\ldots,r_2+1-t\}$, 
\begin{align*}
[x_j,y_j] \equiv [x_j,\prod_{i=1-t}^{r_2+1-t} x_i^{\widetilde{l}_{j,i}}] 
\equiv 
\prod_{i=1-t}^{j-1} [x_i, x_j]^{-\widetilde{l}_{j,i}} 
\prod_{i'=j+1}^{r_2+1-t} [x_j, x_{i'}]^{\widetilde{l}_{j,i'}} 
\pmod{[[F,F],F]}. 
\end{align*}
If $j \in \{r_2+2-t,\ldots,r_1+2r_2\}$, 
\begin{align*}
[x_j,y_j] 
&\equiv [\prod_{i=1-t}^{r_2+1-t} x_i^{c_{j,i}},\prod_{i'=1-t}^{r_2+1-t} x_{i'}^{\widetilde{l}_{j,i'}}] 
\equiv \prod_{(i,i') \in \mathcal{I}} [x_i, x_{i'}]^{c_{j,i}\widetilde{l}_{j,i'}-c_{j,i'}\widetilde{l}_{j,i}}
\pmod{[[F,F],F]}. 
\end{align*}
Therefore, for each $j \in \{1,2,\ldots,r_1+2r_2\}$, 
\begin{align*}
[x_j,y_j] \equiv \prod_{(i,i') \in \mathcal{I}} [x_i, x_{i'}]^{a_{j,(i,i')}} \pmod{[[F,F],F]} 
\end{align*}
where 
\begin{align}\label{eq:ajii}
a_{j,(i,i')} =
\left\{
\begin{array}{cl}
\widetilde{l}_{j,i'} & \text{if $j=i$}, \\
-\widetilde{l}_{j,i} & \text{if $j=i'$}, \\
0 & \text{if $j \in \{1,2,\ldots,r_2+1-t\} \setminus \{i,i'\}$}, \\
c_{j,i}\widetilde{l}_{j,i'}-c_{j,i'}\widetilde{l}_{j,i} & \text{if $j \in \{r_2+2-t,\ldots,r_1+2r_2\}$}. 
\end{array}\right.
\end{align}
Then 
\begin{align}\label{eq:A}
A=(a_{j,(i,i')})_{1 \le j \le r_1+2r_2 \atop (i,i') \in \mathcal{I}}
\end{align}
is a $(r_1+2r_2) \times \frac{r_2(r_2+1)}{2}$ presentation matrix of $X(\widetilde{k})_{\widetilde{\varGamma}}$. 
Hence $A$ satisfies the condition \ref{main:a_is_lk}. 

Let $\mathfrak{m}$ be the maximal ideal of $\widetilde{\varLambda}$. 
Since 
\[
[F,F]/[F,F]^p[[F,F],F]R \simeq X(\widetilde{k})_{\widetilde{\varGamma}}/p \simeq X(\widetilde{k})/\mathfrak{m},
\]
the following exact sequence is induced:
\[
\entrymodifiers={+!!<0pt,\fontdimen22\textfont2>}
\xymatrix@R=0.8\baselineskip{
0 \ar[r] & R/R \cap [F,F]^p[[F,F],F] \ar[r]^-{f_{\overline{A}}} & [F,F]/[F,F]^p[[F,F],F] \ar[r] & X(\widetilde{k})/\mathfrak{m} \ar[r] & 0 \\
& \mathbb F_p^{r_1+2r_2} \ar@{->>}[u] \ar[r]^-{\overline{A}} & \mathbb F_p^{\frac{r_2(r_2+1)}{2}} \ar@{}[u]|-{\rotatebox[origin=c]{90}{$\simeq$}} & & 
} 
\]
Then 
$\dim_{\mathbb F_p} X(\widetilde{k})/\mathfrak{m} = \mathrm{dim}_{\mathbb F_p} \mathrm{Coker}\,f_{\overline{A}}=\frac{r_2(r_2+1)}{2}-\mathrm{rank}\,\overline{A}$. 
By Nakayama's lemma, 
$\dim_{\mathbb F_p} X(\widetilde{k})/\mathfrak{m}$ 
is the minimal number of generators of the $\widetilde{\varLambda}$-module $X(\widetilde{k})$. 
Therefore $A$ also satisfies the conditions \ref{main:lmdrk}. 
Thus Theorem \ref{thm:main} has been proved. 

\begin{corollary}[cf.\ {\cite[Lemma 3.9]{Fuj22}} or {\cite[Th\'eor\`eme 10]{Jau24}}]\label{cor:r2<3}
Under the assumptions of Theorem \ref{thm:main}, 
$r_2 \le \frac{3+\sqrt{9+8(r_1+m)}}{2}$ if $X(\widetilde{k})$ is generated by at most $m$ elements as a $\widetilde{\varLambda}$-module. 
In particular, $r_2 \le 3$ if $r_1+m \le 1$, for example, if $k$ is a CM-field and $X(\widetilde{k})$ is trivial. 
\end{corollary}

\begin{proof}
If there is a surjective homomorphism $\widetilde{\varLambda}^m \twoheadrightarrow X(\widetilde{k})$, 
we have $\dim_{\mathbb F_p} X(\widetilde{k})/\mathfrak{m} \le m$. 
Then $\frac{r_2(r_2+1)}{2}-m \le \mathrm{rank}\,\overline{A} \le r_1+2r_2$, 
and hence $r_2 \le \frac{3+\sqrt{9+8(r_1+m)}}{2}$. 
Since $3 < \frac{3+\sqrt{17}}{2} < 4$, we obtain the latter statement. 
\end{proof}

\subsection{Minimum $\lambda$-invariants.}\label{subsec:minlam} 
Suppose that $r_2 \ge 2$ and $t=0$ in the situation of Theorem \ref{thm:Ghyoji}. 
Put 
\[
\mathcal{I}_2=\{(j,j') \in \mathbb Z^2 \mid 2 \le j < j' \le r_2+1\}, 
\] 
which is an lexicographically ordered subset of $\mathcal{I}$. 
Put a $\frac{r_2(r_2-1)}{2} \times \frac{r_2(r_2+1)}{2}$ matrix 
\[
A'=(a'_{(j,j'),(i,i')})_{(j,j') \in \mathcal{I}_2 \atop (i,i') \in \mathcal{I}} 
\]
with entries 
\[
a'_{(j,j'),(i,i')} = 
\left\{
\begin{array}{cl}
1 & \text{if $(i,i') \in \{(1,j), (j,j')\}$}, \\
-1 & \text{if $(i,i')=(1,j')$}, \\
0 & \text{otherwise}. 
\end{array}\right.
\]

\begin{theorem}\label{thm:lambdar2}
Assume that $r_2 \ge 2$, $t=0$ and {\ref{asm:podd}, \ref{asm:spcomp}, \ref{asm:B}} 
under the settings in \S\ref{subsec:settei}. 
Put a $(r_1+\frac{r_2(r_2+3)}{2}) \times \frac{r_2(r_2+1)}{2}$ matrix 
$A^+=\left(\begin{smallmatrix} A \\ A' \end{smallmatrix}\right)$,  
which is a concatenation of $A$ defined as $\eqref{eq:A}$ and $A'$ defined as above. 
Then $X(k_{\mathrm{cyc}}) \simeq \mathbb Z_p^{r_2}$ if and only if $\mathrm{rank}\,\overline{A^+} = \frac{r_2(r_2+1)}{2}$ for $\mathbb F_p$-matrix $\overline{A^+}$. 
\end{theorem}

\begin{proof}
Let $H$ be the normal subgroup of $F$ normally generated by $r_2$ elements $x_1x_i^{-1}$ $(i \in \{2,\ldots,r_2+1\})$. 
Since $F/H \simeq \mathbb Z_p$, we have $[F,F] \subset H$, and hence $R \subset H$. 
By the choices of $\alpha_i$ and $\tau_i$, 
$\tau_i|_{k_{\mathrm{cyc}}}=\widetilde{\tau}_i|_{k_{\mathrm{cyc}}}$ 
for $\widetilde{\tau}_i \in \mathrm{Gal}(k_{\mathfrak{p}_i}(\zeta_{p^{\infty}})/k_{\mathfrak{p}_i})$ 
satisfying $\widetilde{\tau}_i(\zeta_{p^n}^{1+p})=\zeta_{p^n}$ for all $n$. 
Since $\tau_i|_{k_{\mathrm{cyc}}}=\tau_1|_{k_{\mathrm{cyc}}}$ for all $i \in \{2,\ldots,r_2+1\}$, 
we have $H/R \simeq \widetilde{H}=\mathrm{Gal}(\widetilde{L}(k_{\mathrm{cyc}})/k_{\mathrm{cyc}})$ under the presentation $\pi$. 
Then $[H,H]R/R \simeq \mathrm{Gal}(\widetilde{L}(k_{\mathrm{cyc}})/L(k_{\mathrm{cyc}}))$, 
and $[H,H][[F,F],F]/[[F,F],F]$ is generated by $\frac{(r_2-1)r_2}{2}$ elements 
$[x_1x_j^{-1},x_1x_{j'}^{-1}][[F,F],F]$ $((j,j') \in \mathcal{I}_2)$. 
Since 
\[
[x_1x_j^{-1},x_1x_{j'}^{-1}] \equiv [x_1,x_j][x_1,x_{j'}]^{-1}[x_j,x_{j'}] \pmod{[[F,F],F]} , 
\]
$A'$ is a presentation matrix of $[F,F]/[H,H][[F,F],F]$. 
Put $\mathcal{G}=\mathrm{Gal}(L(k_{\mathrm{cyc}})/k) \simeq F/[H,H]R$. 
There is an exact sequence 
\[
\entrymodifiers={+!!<0pt,\fontdimen22\textfont2>}
\xymatrix@R=0.8\baselineskip{
0 \ar[r] & R[H,H][[F,F],F]/[[F,F],F] \ar[r]^-{f_{A^+}} & [F,F]/[[F,F],F] \ar[r] & [\mathcal{G},\mathcal{G}]/[[\mathcal{G},\mathcal{G}],\mathcal{G}] \ar[r] & 0 \\
& \mathbb Z_p^{r_1+2r_2+\frac{(r_2-1)r_2}{2}} \ar@{->>}[u] \ar[r]^-{A^+} & \mathbb Z_p^{\frac{r_2(r_2+1)}{2}} \ar@{}[u]|-{\rotatebox[origin=c]{90}{$\simeq$}} & & 
} 
\]
with a presentation matrix $A^+$ of $[\mathcal{G},\mathcal{G}]/[[\mathcal{G},\mathcal{G}],\mathcal{G}]$. 
Then  
\[
X(k_{\mathrm{cyc}}) \simeq \mathbb Z_p^{r_2} 
\ \Leftrightarrow \ L(k_{\mathrm{cyc}})=\widetilde{k} 
\ \Leftrightarrow \ [\mathcal{G},\mathcal{G}] \simeq 1 
\ \Leftrightarrow \ [\mathcal{G},\mathcal{G}]/[[\mathcal{G},\mathcal{G}],\mathcal{G}] \simeq 1 , 
\]
and hence we obtain the statement. 
\end{proof}

\section{Triviality of Iwasawa modules}\label{sec:X0}

We focus on the triviality of $X(\widetilde{k})$ for a CM-field $k$. 
In addition to {\ref{asm:podd}--\ref{asm:B}}, we assume the following conditions: 
\begin{enumerate}
\setlength{\itemsep}{0pt}
\setcounter{enumi}{4}
\renewcommand{\labelenumi}{\rm (a\arabic{enumi})}
\renewcommand{\theenumi}{\rm (a\arabic{enumi})}
\item\label{asm:cm}
$k$ is a CM-field of degree $2r_2 \le 6$. Then $r_1=0$. 
\item\label{asm:t1}
$t \le 1$. 
\item\label{asm:ift1}
If $t=1$, $Cl_k/p=\langle [\mathfrak{p}_i]Cl_k^p \rangle$ for some $\mathfrak{p}_i|p$, i.e., $l_{i,0} \not\equiv 0 \pmod{p}$ for some $i$. 
\end{enumerate}
These are necessary conditions of the triviality of $X(\widetilde{k})$ for a CM-field $k$ (cf.\ \cite{Oka**}). 
Let $\iota$ denote the complex conjugation. 
Symbolic computations in this section can be done by using \texttt{PARI/GP} \cite{PARI}.

\subsection{Quadratic fields.} 

In addition to the settings in \S\ref{subsec:settei} and the assumptions 
\ref{asm:podd}--\ref{asm:ift1}, 
we assume that $r_2=1$. 
Then $P=\{\mathfrak{p}_1,\mathfrak{p}_2\}$, 
and $X(\widetilde{k}) \simeq 0$ if and only if $X(k_{\mathrm{cyc}}) \simeq \mathbb Z_p$. 

If $t=0$, then $V/(k^{\times})^p \simeq 0$ and $Z$ is an empty matrix. 
Hence $\widetilde{G}/[\widetilde{G},\widetilde{G}] \simeq \mathbb Z_p^2$, 
$\widetilde{G} \simeq \langle x_1,x_2 \,|\, [x_1,y_1], [x_2,y_2] \rangle$, 
and $\widetilde{l}_{j,i}=l_{j,i}h^{-1}$ for $j,i \in \{1,2\}$. 
Since $\mathfrak{p}_2=\iota(\mathfrak{p}_1)$ and $\iota(\varpi_2)=\pm \varpi_1$, 
we have $l_{1,2}=l_{2,1}$. 
Then 
\[
A
= \begin{pmatrix} a_{1,(1,2)} \\ a_{2,(1,2)} \end{pmatrix}
= \begin{pmatrix} \widetilde{l}_{1,2} \\ -\widetilde{l}_{2,1} \end{pmatrix}
= h^{-1}\begin{pmatrix} l_{1,2} \\ -l_{1,2} \end{pmatrix}. 
\]
Thus we revisit Gold's theorem: 

\begin{proposition}[{cf.\ \cite{Gol74}, see also \cite[Theorem 4.2]{Miz19}}]\label{prop:Gold}
Assume that $p \neq 2$, $k$ is an imaginary quadratic field in which $p$ splits as $p\mathcal{O}_k=\mathfrak{p}_1\mathfrak{p}_2$, and $t=0$, i.e., $p \nmid |Cl_k|$. 
Then $X(\widetilde{k}) \simeq 0$ if and only if $l_{1,2} \not\equiv 0 \pmod{p}$, i.e., 
$\varpi_1^{p-1} \not\equiv 1 \pmod{\mathfrak{p}_2^2}$. 
\end{proposition}

Suppose $t=1$. 
Then $Z=(z_{0,1}, z_{0,2})$. 
By the assumption \ref{asm:B}, 
$B=(z_{0,2})$, $\det B = z_{0,2} \in \mathbb Z_p^{\times}$, and $b_{2,0}=z_{0,2}^{-1}$. 
Then 
$\widetilde{G}/[\widetilde{G},\widetilde{G}] \simeq \mathbb Z_p^2$, 
$\widetilde{G} \simeq \langle x_0,x_1 \,|\, [x_1,y_1], [x_2,y_2] \rangle$, 
$c_{2,1}=-b_{2,0} z_{0,1}=-z_{0,2}^{-1}z_{0,1}$, 
and 
\begin{align*}
\widetilde{l}_{j,1}&=(l_{j,1} +c_{2,1}  l_{j,2} ) h^{-1} 
=\bigg\{
\begin{array}{ll}
- z_{0,2}^{-1} z_{0,1}  l_{1,2} h^{-1} & \text{if $j=1$}, \\
l_{2,1} h^{-1} & \text{if $j=2$},
\end{array}
\\
\widetilde{l}_{j,0}&=(l_{j,0}-q_0 b_{2,0} l_{j,2}) h^{-1} 
=\bigg\{
\begin{array}{ll}
(l_{1,0}-q_0 z_{0,2}^{-1} l_{1,2}) h^{-1}& \text{if $j=1$}, \\
l_{2,0}h^{-1}& \text{if $j=2$}. 
\end{array}
\end{align*}
Hence 
\[
A=
\begin{pmatrix}
a_{1,(0,1)} \\
a_{2,(0,1)} 
\end{pmatrix}
=
\begin{pmatrix}
-\widetilde{l}_{1,0} \\
-q_0 b_{2,0} \widetilde{l}_{2,1}-c_{2,1}\widetilde{l}_{2,0}
\end{pmatrix}
\equiv
h^{-1}
\begin{pmatrix}
-l_{1,0} \\
z_{0,2}^{-1}z_{0,1}l_{2,0}
\end{pmatrix}
\pmod{p}. 
\]
Since $[\mathfrak{p}_1]=[\mathfrak{p}_2]^{-1}$, 
$l_{1,0} \in \mathbb Z_p^{\times}$ by the assumption \ref{asm:ift1}. 
Therefore $\mathrm{rank}\,\overline{A}=1$, i.e., $X(\widetilde{k}) \simeq 0$. 
Thus we obtain the following. 

\begin{proposition}\label{prop:2t1}
Assume that $p \neq 2$, $k$ is an imaginary quadratic field in which $p$ splits as $p\mathcal{O}_k=\mathfrak{p}_1\mathfrak{p}_2$, and $L(k)/k$ is a nontrivial cyclic extension in which $\mathfrak{p}_1$ is inert. 
Then $X(\widetilde{k}) \simeq 0$ if $\widetilde{G}/[\widetilde{G},\widetilde{G}] \simeq \mathbb Z_p^2$.
\end{proposition}

\subsection{Quartic fields with $t=0$.} 
In addition to the settings in \S\ref{subsec:settei} and the assumptions 
\ref{asm:podd}--\ref{asm:ift1}, 
we assume that $r_2=2$ and $t=0$. 
Then $P=\{\mathfrak{p}_1,\mathfrak{p}_2,\mathfrak{p}_3,\mathfrak{p}_4\}$, 
and $Z=(z_{1,1},z_{1,2},z_{1,3},z_{1,4})$. 
By the assumption \ref{asm:B}, 
$B=(z_{1,4})$, $\det B = z_{1,4} \in \mathbb Z_p^{\times}$, and $b_{4,1}=z_{1,4}^{-1}$. 
Then 
$\widetilde{G}/[\widetilde{G},\widetilde{G}] \simeq \mathbb Z_p^3$, and 
$\widetilde{G} \simeq \langle x_1,x_2,x_3 \,|\, [x_1,y_1], [x_2,y_2], [x_3,y_3], [x_4,y_4] \rangle$. 
By renumbering $\mathfrak{p}_1,\mathfrak{p}_2,\mathfrak{p}_3$, 
we may assume that 
$\iota(\mathfrak{p}_1)=\mathfrak{p}_3$, $\iota(\varpi_1)=\varpi_3$ 
and $\iota(\mathfrak{p}_2)=\mathfrak{p}_4$, $\iota(\varpi_2)=\varpi_4$. 
Since $1+p, \varepsilon_1^{p-1} \in k^+$, we have 
$z_{1,3}=z_{1,1}$, $z_{1,4}=z_{1,2} \in \mathbb Z_p^{\times}$, 
$(l_{3,1},l_{3,2},l_{3,4})=(l_{1,3},l_{1,4},l_{1,2})$ 
and $(l_{4,1},l_{4,2},l_{4,3})=(l_{2,3},l_{2,4},l_{2,1})$. 
For $j \in \{1,2,3,4\}$ and $i \in \{1,2,3\}$, 
\begin{align*}
c_{4,i}
&= - b_{4,1} z_{1,i} = - z_{1,4}^{-1} z_{1,i} , 
\\
\widetilde{l}_{j,i}
&=(l_{j,i} + c_{4,i} l_{j,4}) h^{-1}
=(l_{j,i}-z_{1,4}^{-1} z_{1,i}  l_{j,4}) h^{-1} . 
\end{align*}
Then 
\begin{align*}
A
=(a_{j,(i,i')})_{1 \le j \le 4 \atop (i,i') \in \mathcal{I}}
=\begin{pmatrix}
\widetilde{l}_{1,2} & \widetilde{l}_{1,3} & 0 \\
-\widetilde{l}_{2,1} & 0 & \widetilde{l}_{2,3} \\
0 & -\widetilde{l}_{3,1} & -\widetilde{l}_{3,2} \\
c_{4,1}\widetilde{l}_{4,2}-c_{4,2}\widetilde{l}_{4,1}&
c_{4,1}\widetilde{l}_{4,3}-c_{4,3}\widetilde{l}_{4,1}&
c_{4,2}\widetilde{l}_{4,3}-c_{4,3}\widetilde{l}_{4,2}
\end{pmatrix}, 
\end{align*}
where $\mathcal{I}=\{(1,2),(1,3),(2,3)\}$. 
Moreover, we have $A^+=\left(\begin{smallmatrix} A \\ A' \end{smallmatrix}\right)$ with $A'=(1,-1,1)$. 

\begin{proposition}\label{prop:4t0}
Assume that $p \neq 2$, $k$ is a quartic CM-field in which $p$ splits completely as $p\mathcal{O}_k=\mathfrak{p}_1\mathfrak{p}_2\mathfrak{p}_3\mathfrak{p}_4$, and $t=0$, i.e., $p \nmid |Cl_k|$. 
Moreover, we assume that 
$\iota(\varpi_1)=\varpi_3$, $\iota(\varpi_2)=\varpi_4$, 
and 
$\varepsilon_1^{p-1} \not\equiv 1 \pmod{\mathfrak{p}_4^2}$ 
for a unit $\varepsilon_1$ such that $E_k/p=\langle \varepsilon_1E_k^p \rangle$. 
Put
\[
D=
(l_{1,2}z_{1,1} - l_{1,3}z_{1,2})(l_{2,4}z_{1,1} - l_{2,3}z_{1,2})
+
(l_{1,4}z_{1,1} - l_{1,3}z_{1,2})(l_{2,4}z_{1,1} - l_{2,1}z_{1,2}) . 
\]
Then $\widetilde{G}/[\widetilde{G},\widetilde{G}] \simeq \mathbb Z_p^3$, and 
$X(\widetilde{k}) \simeq 0$ if and only if the following two conditions are satisfied: 
\begin{itemize}
\item 
$l_{1,2}-l_{1,4} \not\equiv 0 \pmod{p}$ or $l_{2,1}-l_{2,3} \not\equiv 0 \pmod{p}$. 
\item
$D \not\equiv 0 \pmod{p}$. 
\end{itemize}
In addition, assume that $k/\mathbb Q$ is abelian, and that $\{\varpi_1,\varpi_2,\varpi_3,\varpi_4\}$ is a $\mathrm{Gal}(k/\mathbb Q)$-set. 
Then $X(k_{\mathrm{cyc}}) \simeq \mathbb Z_p^2$ if and only if the following condition is satisfied: 
\begin{itemize}
\item
$( l_{1,2} + l_{1,3})(l_{1,3} + l_{1,4}) \not\equiv 0 \pmod{p}$ if $\mathrm{Gal}(k/\mathbb Q) \simeq V_4$. 
\item
$( l_{1,2} + l_{1,3})^2+(l_{1,3} + l_{1,4})^2 \not\equiv 0 \pmod{p}$ if $\mathrm{Gal}(k/\mathbb Q) \simeq C_4$. 
\end{itemize}
Moreover, $X(\widetilde{k}) \simeq 0$ if and only if $X(k_{\mathrm{cyc}}) \simeq \mathbb Z_p^2$ and 
$l_{1,2}-l_{1,4} \not\equiv 0 \pmod{p}$.
\end{proposition}

\begin{proof}[Proof of Proposition \ref{prop:4t0} (and Corollary \ref{sttmt:V4t0})]
Recall that $z_{1,2} \in \mathbb Z_p^{\times}$. 
By a symbolic computation, we have 
\[
\left(\begin{array}{cccc|c}
1&0&0&0&0 \\
0&z_{1,2}&0&0&0 \\
h&0&h&0&0 \\
0&h&0&h&0 \\
\hline
0&0&0&0&1 \\
\end{array}\right)
A^+
\begin{pmatrix}
z_{1,1} & 0 & 0 \\
-z_{1,2} & z_{1,2} & 0 \\
0 & -z_{1,1} & 1
\end{pmatrix}
=
\left(\begin{array}{ccc}
-z_{1,2}\widetilde{l}_{3,1} & z_{1,2}\widetilde{l}_{1,3} & 0 \\
-z_{1,1}z_{1,2}\widetilde{l}_{2,1} &
-z_{1,1}z_{1,2}\widetilde{l}_{2,3} &
z_{1,2}\widetilde{l}_{2,3} \\
0&0& l_{1,2} - l_{1,4} \\
0&0& -l_{2,1} + l_{2,3} \\
\hline
z_{1,1}+z_{1,2} & -z_{1,1}-z_{1,2} & 1 
\end{array}\right) , 
\]
and 
\[
\det \begin{pmatrix}
-z_{1,2}\widetilde{l}_{3,1} & z_{1,2}\widetilde{l}_{1,3} \\
-z_{1,1}z_{1,2}\widetilde{l}_{2,1} & -z_{1,1}z_{1,2}\widetilde{l}_{2,3}  
\end{pmatrix} 
=z_{1,1} h^{-2} D. 
\]
If $z_{1,1}+z_{1,2} \equiv 0 \pmod{p}$, 
we have $\mathrm{rank}\,\overline{A^+}=3$ if and only if $D \not\equiv 0 \pmod{p}$. 

Since $X(\widetilde{k}) \simeq 0$ if and only if $\mathrm{rank}\,\overline{A}=3$, 
we obtain the first-half statement if $z_{1,1} \in \mathbb Z_p^{\times}$. 
If $z_{1,1} \equiv 0 \pmod{p}$, 
\[
\begin{pmatrix}
1&0&0&0 \\
0&1&0&0 \\
1&0&1&0 \\
0&1&0&1 \\
\end{pmatrix}
hA
\begin{pmatrix}
1 & 0 & 0 \\
0 & 1 & 0 \\
-1& 0 & 1
\end{pmatrix}
\equiv
\begin{pmatrix}
 l_{1,2} - l_{1,4} & l_{1,3} & 0 \\
-l_{2,1}-l_{2,3} & 0 & l_{2,3} \\
0&0& l_{1,2} - l_{1,4} \\
0&0& -l_{2,1} + l_{2,3} 
\end{pmatrix}
\pmod{p}, 
\]
and $D \equiv z_{1,2}^2 l_{1,3}(l_{2,1} + l_{2,3}) \pmod{p}$. 
Hence we obtain the first-half statement also when $z_{1,1} \equiv 0 \pmod{p}$. 

If $k/\mathbb Q$ is a $V_4$-extension, 
there is an element $g$ of $\mathrm{Gal}(k/\mathbb Q)$ such that
$g(\mathfrak{p}_1)=\mathfrak{p}_2$ and $g(\mathfrak{p}_3)=\mathfrak{p}_4$. 
Then $\mathrm{Gal}(k/\mathbb Q)=\{id,g,\iota|_k,g\iota|_k\}$. 
Since $\varepsilon_1^{(p-1)(1+g)}=1$, we have $z_{1,2}=-z_{1,1}$. 
In particular, $\varepsilon_1^{p-1} \not\equiv 1 \pmod{\mathfrak{p}_1^2}$ if and only if $\varepsilon_1^{p-1} \not\equiv 1 \pmod{\mathfrak{p}_4^2}$. 
Since $\{\varpi_1,\varpi_2,\varpi_3,\varpi_4\}$ is a $\mathrm{Gal}(k/\mathbb Q)$-set, 
$\varpi_2=g(\varpi_1)$ and $\varpi_4=g(\varpi_3)$. 
Then 
$l_{1,2}=l_{2,1}=l_{3,4}=l_{4,3}$, 
$l_{1,3}=l_{3,1}=l_{2,4}=l_{4,2}$, 
$l_{1,4}=l_{4,1}=l_{2,3}=l_{3,2}$, 
and hence 
$l_{2,1}-l_{2,3} = l_{1,2}-l_{1,4}$ and 
$D=2z_{1,1}^2(l_{1,2} + l_{1,3})(l_{1,3} + l_{1,4})$. 
Thus $A$ is written in the form \eqref{eq:V4t0A}, 
and the statements of Corollary \ref{sttmt:V4t0} are also obtained. 

If $k/\mathbb Q$ is a $C_4$-extension, 
there is a generator $g$ of $\mathrm{Gal}(k/\mathbb Q)$ such that 
$g^{i-1}(\mathfrak{p}_1)=\mathfrak{p}_i$ for $i \in \{1,2,3,4\}$. 
Then $g^2=\iota|_k$. 
Since $\varepsilon_1^{(p-1)(1+g)}=1$, we have $z_{1,2}=-z_{1,1}$. 
Since $\{\varpi_1,\varpi_2,\varpi_3,\varpi_4\}$ is a $\mathrm{Gal}(k/\mathbb Q)$-set, 
$\varpi_i=g^{i-1}(\varpi_1)$. 
Then 
$l_{1,2}=l_{2,3}=l_{3,4}=l_{4,1}$, 
$l_{1,3}=l_{2,4}=l_{3,1}=l_{4,2}$, 
$l_{1,4}=l_{2,1}=l_{3,2}=l_{4,3}$, 
and hence 
$l_{2,1}-l_{2,3} = -(l_{1,2}-l_{1,4})$ and 
$D=z_{1,1}^2(( l_{1,2} + l_{1,3})^2+(l_{1,3} + l_{1,4})^2)$. 

By Theorems \ref{thm:main} and \ref{thm:lambdar2}, 
we obtain the latter-half statement for abelian extension $k/\mathbb Q$. 
\end{proof}

\begin{remark}
Also when $k/\mathbb Q$ is a $C_4$-extension, we have 
\[
A=h^{-1}\begin{pmatrix}
 l_{1,2}-l_{1,4} & l_{1,3}+l_{1,4} & 0 \\
-l_{1,3}-l_{1,4} & 0 & l_{1,2}+l_{1,3} \\
0 & -l_{1,2}-l_{1,3} & l_{1,2}-l_{1,4} \\
l_{1,2}+l_{1,3} &  -l_{1,2}+l_{1,4} & -l_{1,3}-l_{1,4}
\end{pmatrix}
\]
with similar statements to Corollary \ref{sttmt:V4t0}. 
\end{remark}

\subsection{Quartic fields with $t=1$.} 
In addition to the settings in \S\ref{subsec:settei} and the assumptions 
\ref{asm:podd}--\ref{asm:ift1}, 
we assume that $r_2=2$ and $t=1$. 
Then $P=\{\mathfrak{p}_1,\mathfrak{p}_2,\mathfrak{p}_3,\mathfrak{p}_4\}$, 
and 
\[
Z=\begin{pmatrix}
z_{0,1} & z_{0,2} & z_{0,3} & z_{0,4} \\
z_{1,1} & z_{1,2} & z_{1,3} & z_{1,4} 
\end{pmatrix} . 
\]
By the assumption \ref{asm:B}, 
$\det B = z_{0,3}z_{1,4}-z_{0,4}z_{1,3} \in \mathbb Z_p^{\times}$ and 
\[
\begin{pmatrix}
b_{3,0} & b_{3,1} \\ 
b_{4,0} & b_{4,1}  
\end{pmatrix}
=B^{-1}
=\frac{1}{\det B}
\begin{pmatrix}
z_{1,4} & -z_{0,4} \\ 
-z_{1,3} & z_{0,3}  
\end{pmatrix} . 
\]
Then 
$\widetilde{G}/[\widetilde{G},\widetilde{G}] \simeq \mathbb Z_p^3$, and 
$\widetilde{G} \simeq \langle x_0,x_1,x_2 \,|\, [x_1,y_1], [x_2,y_2], [x_3,y_3], [x_4,y_4] \rangle$, 
\begin{align*}
c_{j,i}&=- b_{j,0} z_{0,i} - b_{j,1} z_{1,i} &&(j \in \{3,4\},\ i \in \{1,2\}), \\
\widetilde{l}_{j,i}&=(l_{j,i}+ c_{3,i}  l_{j,3}+ c_{4,i}  l_{j,4} ) h^{-1} &&(j \in \{1,2,3,4\},\ i \in \{1,2\}), \\
\widetilde{l}_{j,0}&=(l_{j,0} -q_0b_{3,0} l_{j,3}-q_0b_{4,0} l_{j,4} ) h^{-1} &&(j \in \{1,2,3,4\}), 
\end{align*}
and 
\[
A=(a_{j,(i,i')})_{1 \le j \le 4 \atop (i,i') \in \mathcal{I}}
=\begin{pmatrix}
-\widetilde{l}_{1,0} & 0 & \widetilde{l}_{1,2} \\
0 & -\widetilde{l}_{2,0} & -\widetilde{l}_{2,1} \\
-q_0 b_{3,0} \widetilde{l}_{3,1}-c_{3,1}\widetilde{l}_{3,0} & 
-q_0 b_{3,0} \widetilde{l}_{3,2}-c_{3,2}\widetilde{l}_{3,0} & 
c_{3,1}\widetilde{l}_{3,2}-c_{3,2}\widetilde{l}_{3,1}
\\
-q_0 b_{4,0} \widetilde{l}_{4,1}-c_{4,1}\widetilde{l}_{4,0} & 
-q_0 b_{4,0} \widetilde{l}_{4,2}-c_{4,2}\widetilde{l}_{4,0} & 
c_{4,1}\widetilde{l}_{4,2}-c_{4,2}\widetilde{l}_{3,1}
\end{pmatrix} , 
\]
where $\mathcal{I}=\{(0,1),(0,2),(1,2)\}$. 
By a symbolic computation, 
one can see that 
\[
(\det B)h A \equiv 
\begin{pmatrix}
\left|\begin{matrix}\begin{smallmatrix}
0 & l_{1,0} & 0 & 0 \\
1 & 0 & 0 & 0 
\end{smallmatrix} \\ Z \end{matrix}\right|
& 
0
&
-\left|\begin{matrix}\begin{smallmatrix}
0 & l_{1,2} & l_{1,3} & l_{1,4} \\
1 & 0 & 0 & 0 
\end{smallmatrix} \\ Z \end{matrix}\right|
\\	
0 
& 
-\left|\begin{matrix}\begin{smallmatrix}
l_{2,0} & 0 & 0 & 0 \\
0 & 1 & 0 & 0 
\end{smallmatrix} \\ Z \end{matrix}\right|
&
-\left|\begin{matrix}\begin{smallmatrix}
l_{2,1} & 0 & l_{2,3} & l_{2,4} \\
0 & 1 & 0 & 0 
\end{smallmatrix} \\ Z \end{matrix}\right|
\\
\left|\begin{matrix}\begin{smallmatrix}
0 & l_{3,0} & 0 & 0 \\
0 & 0 & 1 & 0 
\end{smallmatrix} \\ Z \end{matrix}\right|
& 
-\left|\begin{matrix}\begin{smallmatrix}
l_{3,0} & 0 & 0 & 0 \\
0 & 0 & 1 & 0 
\end{smallmatrix} \\ Z \end{matrix}\right|
&
-\left|\begin{matrix}\begin{smallmatrix}
l_{3,1} & l_{3,2} & 0 & l_{3,4} \\
0 & 0 & 1 & 0 
\end{smallmatrix} \\ Z \end{matrix}\right|
\\
\left|\begin{matrix}\begin{smallmatrix}
0 & l_{4,0} & 0 & 0 \\
0 & 0 & 0 & 1 
\end{smallmatrix} \\ Z \end{matrix}\right|
& 
-\left|\begin{matrix}\begin{smallmatrix}
l_{4,0} & 0 & 0 & 0 \\
0 & 0 & 0 & 1 
\end{smallmatrix} \\ Z \end{matrix}\right|
&
-\left|\begin{matrix}\begin{smallmatrix}
l_{4,1} & l_{4,2} & l_{4,3} & 0 \\
0 & 0 & 0 & 1 
\end{smallmatrix} \\ Z \end{matrix}\right|
\end{pmatrix}
\pmod{p}. 
\]

\subsection{Sextic fields with $t=0$.} 

In addition to the settings in \S\ref{subsec:settei} and the assumptions 
\ref{asm:podd}--\ref{asm:ift1}, 
we assume that $r_2=3$ and $t=0$. 
Then $P=\{\mathfrak{p}_1,\mathfrak{p}_2,\mathfrak{p}_3,\mathfrak{p}_4,\mathfrak{p}_5,\mathfrak{p}_6\}$, 
and 
\[
Z=\begin{pmatrix}
z_{1,1} & z_{1,2} & z_{1,3} & z_{1,4}  & z_{1,5}  & z_{1,6} \\
z_{2,1} & z_{2,2} & z_{2,3} & z_{2,4}  & z_{2,5}  & z_{2,6} 
\end{pmatrix}. 
\]
By the assumption \ref{asm:B}, 
$\det B = z_{1,5}z_{2,6}-z_{1,6}z_{2,5} \in \mathbb Z_p^{\times}$ and 
\[
\begin{pmatrix}
b_{5,1} & b_{5,2} \\ 
b_{6,1} & b_{6,2}  
\end{pmatrix}
=B^{-1}
=\frac{1}{\det B}
\begin{pmatrix}
z_{2,6} & -z_{1,6} \\ 
-z_{2,5} & z_{1,5}  
\end{pmatrix} . 
\]
Then 
$\widetilde{G}/[\widetilde{G},\widetilde{G}] \simeq \mathbb Z_p^4$, and 
$\widetilde{G} \simeq \langle x_1,x_2,x_3,x_4 \,|\, [x_1,y_1], \ldots, [x_6,y_6] \rangle$. 
Since $1+p, \varepsilon_1^{p-1}, \varepsilon_2^{p-1} \in k^+$, 
we have $z_{i,j}=z_{i,j'}$ if $\iota(\mathfrak{p}_j)=\mathfrak{p}_{j'}$. 
Therefore $\iota(\mathfrak{p}_5) \neq \mathfrak{p}_6$. 
By renumbering the elements of $P$, we may assume that 
$\iota(\mathfrak{p}_j)=\mathfrak{p}_{j+3}$ and $\iota(\varpi_j)=\varpi_{j+3}$ for $j \in \{1,2,3\}$. 
Then 
\[
\begin{pmatrix}
z_{1,1} & z_{1,2} & z_{1,3} \\
z_{2,1} & z_{2,2} & z_{2,3}  
\end{pmatrix}
=
\begin{pmatrix}
z_{1,4}  & z_{1,5}  & z_{1,6} \\
z_{2,4}  & z_{2,5}  & z_{2,6} 
\end{pmatrix}.
\]
and 
\[
\begin{pmatrix}
l_{4,1} & l_{4,2} & l_{4,3} & l_{4,4} & l_{4,5} & l_{4,6} \\
l_{5,1} & l_{5,2} & l_{5,3} & l_{5,4} & l_{5,5} & l_{5,6} \\
l_{6,1} & l_{6,2} & l_{6,3} & l_{6,4} & l_{6,5} & l_{6,6} 
\end{pmatrix}
=
\begin{pmatrix}
l_{1,4} & l_{1,5} & l_{1,6} &       0 & l_{1,2} & l_{1,3} \\
l_{2,4} & l_{2,5} & l_{2,6} & l_{2,1} &       0 & l_{2,3} \\
l_{3,4} & l_{3,5} & l_{3,6} & l_{3,1} & l_{3,2} &       0 
\end{pmatrix}. 
\]
Recall that 
\begin{align*}
c_{j,i}&=-b_{j,1} z_{1,i}-b_{j,2} z_{2,i} 
&& (j \in \{5,6\},\ i \in \{1,2,3,4\}), \\
\widetilde{l}_{j,i} 
&=( l_{j,i} + c_{5,i} l_{j,5} + c_{6,i} l_{j,6} ) h^{-1} 
&& (j \in \{1,\ldots,6\},\ i \in \{1,2,3,4\}) .
\end{align*}
Then 
\begin{align*}
A&=(a_{j,(i,i')})_{1 \le j \le 6 \atop (i,i') \in \mathcal{I}} \\
&=\left(\begin{smallmatrix}
\widetilde{l}_{1,2} &
\widetilde{l}_{1,3} & 
\widetilde{l}_{1,4} & 
&
&
\\
-\widetilde{l}_{2,1} &
&
&
\widetilde{l}_{2,3} &
\widetilde{l}_{2,4}
\\
&
-\widetilde{l}_{3,1} &
&
-\widetilde{l}_{3,2} &
&
\widetilde{l}_{3,4} \\
&
&
-\widetilde{l}_{4,1} &
&
-\widetilde{l}_{4,2} &
-\widetilde{l}_{4,3} 
\\
{c_{5,1}\widetilde{l}_{5,2}-c_{5,2}\widetilde{l}_{5,1}}&
{c_{5,1}\widetilde{l}_{5,3}-c_{5,3}\widetilde{l}_{5,1}}&
{c_{5,1}\widetilde{l}_{5,4}-c_{5,4}\widetilde{l}_{5,1}}&
{c_{5,2}\widetilde{l}_{5,3}-c_{5,3}\widetilde{l}_{5,2}}&
{c_{5,2}\widetilde{l}_{5,4}-c_{5,4}\widetilde{l}_{5,2}}&
{c_{5,3}\widetilde{l}_{5,4}-c_{5,4}\widetilde{l}_{5,3}}
\\
{c_{6,1}\widetilde{l}_{6,2}-c_{6,2}\widetilde{l}_{6,1}}&
{c_{6,1}\widetilde{l}_{6,3}-c_{6,3}\widetilde{l}_{6,1}}&
{c_{6,1}\widetilde{l}_{6,4}-c_{6,4}\widetilde{l}_{6,1}}&
{c_{6,2}\widetilde{l}_{6,3}-c_{6,3}\widetilde{l}_{6,2}}&
{c_{6,2}\widetilde{l}_{6,4}-c_{6,4}\widetilde{l}_{6,2}}&
{c_{6,3}\widetilde{l}_{6,4}-c_{6,4}\widetilde{l}_{6,3}}
\end{smallmatrix} \right), 
\end{align*}
where $\mathcal{I}=\{(1,2),(1,3),(1,4),(2,3),(2,4),(3,4)\}$. 
Moreover, we have $A^+=\left(\begin{smallmatrix} A \\ A' \end{smallmatrix}\right)$ with 
\[
A'=
\begin{pmatrix}
1 & -1 &  0 &  1 &  0 &  0 \\
1 &  0 & -1 &  0 &  1 &  0 \\
0 &  1 & -1 &  0 &  0 &  1 
\end{pmatrix} . 
\]

\begin{proposition}\label{prop:6t0}
Assume that $p \neq 2$, $k$ is a sextic CM-field in which $p$ splits completely as $p\mathcal{O}_k=\mathfrak{p}_1\mathfrak{p}_2\mathfrak{p}_3\mathfrak{p}_4\mathfrak{p}_5\mathfrak{p}_6$, and $t=0$, i.e., $p \nmid |Cl_k|$. 
Moreover, we assume that 
$\iota(\varpi_1)=\varpi_4$, $\iota(\varpi_2)=\varpi_5$, $\iota(\varpi_3)=\varpi_6$, 
and $z_{1,5}z_{2,6} \not\equiv z_{1,6}z_{2,5} \pmod{p}$. 
Put
\[
M
=\begin{pmatrix}
\left|\begin{matrix}\begin{smallmatrix} l_{1,4} & l_{1,5} & l_{1,3} \end{smallmatrix} \\ Z' \end{matrix}\right|&
-\left|\begin{matrix}\begin{smallmatrix} l_{1,4} & l_{1,5} & l_{1,6} \end{smallmatrix} \\ Z' \end{matrix}\right|&
-\left|\begin{matrix}\begin{smallmatrix} l_{1,4} & l_{1,2} & l_{1,3} \end{smallmatrix} \\ Z' \end{matrix}\right|\\
-\left|\begin{matrix}\begin{smallmatrix} l_{2,4} & l_{2,5} & l_{2,3} \end{smallmatrix} \\ Z' \end{matrix}\right|&
\left|\begin{matrix}\begin{smallmatrix} l_{2,4} & l_{2,5} & l_{2,6} \end{smallmatrix} \\ Z' \end{matrix}\right|&
-\left|\begin{matrix}\begin{smallmatrix} l_{2,1} & l_{2,5} & l_{2,6} \end{smallmatrix} \\ Z' \end{matrix}\right|\\
-\left|\begin{matrix}\begin{smallmatrix} l_{3,1} & l_{3,2} & l_{3,6} \end{smallmatrix} \\ Z' \end{matrix}\right|&
-\left|\begin{matrix}\begin{smallmatrix} l_{3,4} & l_{3,5} & l_{3,6} \end{smallmatrix} \\ Z' \end{matrix}\right|&
\left|\begin{matrix}\begin{smallmatrix} l_{3,1} & l_{3,5} & l_{3,6} \end{smallmatrix} \\ Z' \end{matrix}\right|
\end{pmatrix}, 
\]
where $Z'=\left(\begin{smallmatrix}
z_{1,1} & z_{1,2} & z_{1,3} \\
z_{2,1} & z_{2,2} & z_{2,3}  
\end{smallmatrix}\right)$. 
Then $\widetilde{G}/[\widetilde{G},\widetilde{G}] \simeq \mathbb Z_p^4$, and 
$X(\widetilde{k}) \simeq 0$ if and only if the following two conditions are satisfied: 
\begin{itemize}
\item 
$(l_{1,2}-l_{1,5})(l_{2,3}-l_{2,6})(l_{3,1}-l_{3,4}) \not\equiv (l_{1,3}-l_{1,6})(l_{2,1}-l_{2,4})(l_{3,2}-l_{3,5}) \pmod{p}$. 
\item
$\det M \not\equiv 0 \pmod{p}$. 
\end{itemize}
If $k/\mathbb Q$ is a $C_6$-extension, 
we may assume that $\varpi_{i+1}=g(\varpi_i)$ for $i \in \{1,\ldots,5\}$ with a generator $g$ of $\mathrm{Gal}(k/\mathbb Q)$ 
after a suitable permutation on $P$. 
Then $X(k_{\mathrm{cyc}}) \simeq \mathbb Z_p^3$ if and only if 
$m_1m_2 \not\equiv 0 \pmod{p}$ for 
\begin{align*}
m_1 &= l_{1,2} +l_{1,4} +l_{1,6}, \\ 
m_2 &=(l_{1,2} +2l_{1,3} +2l_{1,4}+l_{1,5})(l_{1,3} +2l_{1,4} +2l_{1,5}+l_{1,6})+(l_{1,2}+l_{1,3}-l_{1,5}-l_{1,6})^2 . 
\end{align*}
Moreover, $X(\widetilde{k}) \simeq 0$ if and only if $X(k_{\mathrm{cyc}}) \simeq \mathbb Z_p^3$ and 
$(l_{1,2}-l_{1,5})^3+(l_{1,3}-l_{1,6})^3 \not\equiv 0 \pmod{p}$. 
\end{proposition}

\begin{proof}
Let $d_1,d_2,d_3$ be the components of the cross product 
\[
(d_1,d_2,d_3)=
\left(
\left|\begin{matrix}
z_{1,2} & z_{1,3} \\
z_{2,2} & z_{2,3}  
\end{matrix}\right|
,
\left|\begin{matrix}
z_{1,3} & z_{1,1} \\
z_{2,3} & z_{2,1}  
\end{matrix}\right|
,
\left|\begin{matrix}
z_{1,1} & z_{1,2} \\
z_{2,1} & z_{2,2}  
\end{matrix}\right|
\right)
= \boldsymbol{z}'_1 \times \boldsymbol{z}'_2
\]
of the row vectors $\boldsymbol{z}'_1=(z_{1,1},z_{1,2},z_{1,3})$, $\boldsymbol{z}'_2=(z_{2,1},z_{2,2},z_{2,3})$ of $Z'$. 
Then $d_1=\det B \in \mathbb Z_p^{\times}$, and 
\[
M
=\begin{pmatrix}
 d_1l_{1,4} +d_2l_{1,5} +d_3l_{1,3}&
-d_1l_{1,4} -d_2l_{1,5} -d_3l_{1,6}&
-d_1l_{1,4} -d_2l_{1,2} -d_3l_{1,3}\\
-d_1l_{2,4} -d_2l_{2,5} -d_3l_{2,3}&
 d_1l_{2,4} +d_2l_{2,5} +d_3l_{2,6}&
-d_1l_{2,1} -d_2l_{2,5} -d_3l_{2,6}\\
-d_1l_{3,1} -d_2l_{3,2} -d_3l_{3,6}&
-d_1l_{3,4} -d_2l_{3,5} -d_3l_{3,6}&
 d_1l_{3,1} +d_2l_{3,5} +d_3l_{3,6}
\end{pmatrix}. 
\]
By a symbolic computation, one can see that 
\[
{\scriptsize \begin{pmatrix}
1& 0& 0&0&0&0\\
0& 1& 0&0&0&0\\
0& 0& 1&0&0&0\\
-1& 0& 0&1&0&0\\
0&-1& 0&0&1&0\\
0& 0&-1&0&0&1
\end{pmatrix}}
d_1 hA
{\scriptsize \begin{pmatrix}
1& 0& 0&0&0&0\\
0& 1& 0&0&0&0\\
0& 0& 0&1&0&0\\
0& 0& 1&0&0&0\\
1& 0& 0&0&1&0\\
0& 1& 0&0&0&1
\end{pmatrix}}
=
\begin{pmatrix}
A_{11} & A_{12} \\
O & A_{22} 
\end{pmatrix}
\]
with 
\begin{align*}
A_{11}
&=\begin{pmatrix}
(l_{1,2}-l_{1,5})d_1 & 
(l_{1,3}-l_{1,6})d_1 &
0\\
-(l_{2,1}-l_{2,4})d_1 &
0&
(l_{2,3}-l_{2,6})d_1 \\
0&
-(l_{3,1}-l_{3,4})d_1 &
-(l_{3,2}-l_{3,5})d_1 
\end{pmatrix},
\\
A_{12}
&=\begin{pmatrix}
l_{1,4}d_1+l_{1,5}d_2+l_{1,6}d_3 & 
0&
0\\
0&
l_{2,4}d_1+l_{2,5}d_2+l_{2,6}d_3 &
0\\
0&
0&
l_{3,4}d_1+l_{3,5}d_2+l_{3,6}d_3 
\end{pmatrix},
\\
A_{22}
&=\begin{pmatrix}
-(l_{1,4}d_1+l_{1,2}d_2+l_{1,3}d_3 \atop
+l_{1,4}d_1+l_{1,5}d_2+l_{1,6}d_3)&
(l_{1,2}-l_{1,5})d_1 &
(l_{1,3}-l_{1,6})d_1 \\
(l_{2,1}-l_{2,4})d_2 & 
-(l_{2,1}d_1+l_{2,5}d_2+l_{2,3}d_3 \atop
+l_{2,4}d_1+l_{2,5}d_2+l_{2,6}d_3)&
(l_{2,3}-l_{2,6})d_2 \\
(l_{3,1}-l_{3,4})d_3 & 
(l_{3,2}-l_{3,5})d_3 & 
-(l_{3,1}d_1+l_{3,2}d_2+l_{3,6}d_3 \atop
+l_{3,4}d_1+l_{3,5}d_2+l_{3,6}d_3)
\end{pmatrix},
\end{align*}
and that 
\begin{align*}
\det A_{11}&=d_1^3\big(
(l_{1,2}-l_{1,5})(l_{2,3}-l_{2,6})(l_{3,1}-l_{3,4})-(l_{1,3}-l_{1,6})(l_{2,1}-l_{2,4})(l_{3,2}-l_{3,5}) 
\big), 
\\
\det A_{22} &=2 \det M. 
\end{align*}
Hence we obtain the first-half statement. 

Suppose that $k/\mathbb Q$ is a $C_6$-extension. 
Let $g_0$ be a generator of $\mathrm{Gal}(k/\mathbb Q)$. 
Then $g_0^{\pm 3} = \iota|_k$ and $\mathfrak{p}_6 \in \{g_0^{\pm 1}(\mathfrak{p}_5), g_0^{\pm 4}(\mathfrak{p}_5)\}$. 
Hence $\mathfrak{p}_6 \in \{ g(\mathfrak{p}_5),g\iota (\mathfrak{p}_5) \}$ 
with a generator $g$ of $\mathrm{Gal}(k/\mathbb Q)$. 
We renumber the elements of $P$ as follows: 
If $\mathfrak{p}_6=g\iota(\mathfrak{p}_5)$, we replace $\mathfrak{p}_5$ with $\iota(\mathfrak{p}_5)=\mathfrak{p}_2$. 
(Note that $z_{i,5}=z_{i,2}$.) 
Hence we may assume that $\mathfrak{p}_6=g(\mathfrak{p}_5)$. 
Put $\mathfrak{p}_4=g^2(\mathfrak{p}_5)$ and $\mathfrak{p}_j=\iota(\mathfrak{p}_{j+3})$ for $j \in \{1,2,3\}$. 
Then $\mathfrak{p}_{i+1}=g(\mathfrak{p}_i)$ for $i \in \{1,\ldots,5\}$. 
Since $\varepsilon_i^{(p-1)(1+g+g^2)}=1$, 
we have $z_{i,1}+z_{i,2}+z_{i,3}=0$, and hence $d_1=d_2=d_3=\det B \in \mathbb Z_p^{\times}$. 
Since $l_{j,i}=l_{1,i'}$ for $i' \equiv 1+i-j \pmod{6}$, 
we have 
\[
\det A_{11}
=d_1^3\big(
-(l_{1,2}-l_{1,5})^3-(l_{1,3}-l_{1,6})^3 
\big)
\]
and 
\[
M=d_1 \begin{pmatrix}
 l_{1,4} +l_{1,5} +l_{1,3}&
-l_{1,4} -l_{1,5} -l_{1,6}&
-l_{1,4} -l_{1,2} -l_{1,3}\\
-l_{1,3} -l_{1,4} -l_{1,2}&
 l_{1,3} +l_{1,4} +l_{1,5}&
-l_{1,6} -l_{1,4} -l_{1,5}\\
-l_{1,5} -l_{1,6} -l_{1,4}&
-l_{1,2} -l_{1,3} -l_{1,4}&
 l_{1,5} +l_{1,3} +l_{1,4}
\end{pmatrix}. 
\]
One can see that $\det M=-d_1^3m_1m_2$. 
Moreover, by a symbolic computation, 
\[
\begin{pmatrix}
hE_3 & O & O \\
hE_3 & hE_3 & O \\
O & O & E_3 \\
\end{pmatrix}
A^+
{\scriptsize \begin{pmatrix}
0& 0&  1& 0& 0&  0\\
1& 0&  0& 0& 0&  0\\
1& 1&  1& 0& 0& -1\\
1& 0& -1& 1& 0&  0\\
1& 1&  0& 0& 1& -1\\
0& 1&  1& 0& 0&  0
\end{pmatrix}}
=
\begin{pmatrix}
M' & \ast \\
O & \ast \\
O & E_3 \\
\end{pmatrix}
\]
with 
\[
M'=
\begin{pmatrix}
l_{1,3} + l_{1,4} + l_{1,5} & l_{1,4} + l_{1,5} + l_{1,6} & l_{1,2} + l_{1,4} + l_{1,6} \\
l_{1,2} + l_{1,3} + l_{1,4} & l_{1,3} + l_{1,4} + l_{1,5} & -(l_{1,2} + l_{1,4} + l_{1,6}) \\
-(l_{1,4} + l_{1,5} + l_{1,6}) & l_{1,2} + l_{1,3} + l_{1,4} & l_{1,2} + l_{1,4} + l_{1,6} 
\end{pmatrix}
\]
satisfying $\det M' = m_1m_2$, 
where $E_3=\left(\begin{smallmatrix} 1&0&0 \\ 0&1&0 \\ 0&0&1 \end{smallmatrix}\right)$. 
Hence $\mathrm{rank}\,\overline{A^+}=6$ if and only if $m_1m_2 \not\equiv 0 \pmod{p}$. 
Thus we obtain the statement for $C_6$-extension $k/\mathbb Q$. 
\end{proof}

\begin{remark}
There is no CM-field which is a $S_3$-extension of $\mathbb Q$. 
If $k$ is a CM-field and $\mathrm{Gal}(k/\mathbb Q) \simeq S_3$, 
$k^+$ is totally real, and the conjugate ${k^+}'$ is also real. 
Then $k=k^+{k^+}' \subset \mathbb R$. This is a contradiction. 
\end{remark}

\subsection{Sextic fields with $t=1$.} 

In addition to the settings in \S\ref{subsec:settei} and the assumptions 
\ref{asm:podd}--\ref{asm:ift1}, 
we assume that $r_2=3$ and $t=1$. 
Then $P=\{\mathfrak{p}_1,\mathfrak{p}_2,\mathfrak{p}_3,\mathfrak{p}_4,\mathfrak{p}_5,\mathfrak{p}_6\}$, 
and $Z=(z_{i,j})_{0 \le i \le 2 \atop 1 \le j \le 6}$ is a $3 \times 6$ matrix. 
By the assumption \ref{asm:B}, $\det B \in \mathbb Z_p^{\times}$ for 
\[
B=\begin{pmatrix}
z_{0,4} & z_{0,5} & z_{0,6} \\
z_{1,4} & z_{1,5} & z_{1,6} \\
z_{2,4} & z_{2,5} & z_{2,6} 
\end{pmatrix}, 
\]
and 
\[
\begin{pmatrix}
b_{4,0} & b_{4,1} & b_{4,2} \\ 
b_{5,0} & b_{5,1} & b_{5,2} \\ 
b_{6,0} & b_{6,1} & b_{6,2}  
\end{pmatrix}
=B^{-1}
=
\frac{1}{\det B}
\begin{pmatrix}
\left|\begin{smallmatrix}
z_{1,5} & z_{1,6} \\
z_{2,5} & z_{2,6} 
\end{smallmatrix}\right|
&
-\left|\begin{smallmatrix}
z_{0,5} & z_{0,6} \\
z_{2,5} & z_{2,6} 
\end{smallmatrix}\right|
&
\left|\begin{smallmatrix}
z_{0,5} & z_{0,6} \\
z_{1,5} & z_{1,6} \\
\end{smallmatrix}\right|
\\
-\left|\begin{smallmatrix}
z_{1,4} & z_{1,6} \\
z_{2,4} & z_{2,6} 
\end{smallmatrix}\right|
&
\left|\begin{smallmatrix}
z_{0,4} & z_{0,6} \\
z_{2,4} &z_{2,6} 
\end{smallmatrix}\right|
&
-\left|\begin{smallmatrix}
z_{0,4} & z_{0,6} \\
z_{1,4} & z_{1,6} \\
\end{smallmatrix}\right|
\\
\left|\begin{smallmatrix}
z_{1,4} & z_{1,5} \\
z_{2,4} & z_{2,5} 
\end{smallmatrix}\right|
&
-\left|\begin{smallmatrix}
z_{0,4} & z_{0,5} \\
z_{2,4} & z_{2,5} 
\end{smallmatrix}\right|
&
\left|\begin{smallmatrix}
z_{0,4} & z_{0,5} \\
z_{1,4} & z_{1,5} 
\end{smallmatrix}\right|
\end{pmatrix}
\]
Then 
$\widetilde{G}/[\widetilde{G},\widetilde{G}] \simeq \mathbb Z_p^4$, and 
$\widetilde{G} \simeq \langle x_0,x_1,x_2,x_3 \,|\, [x_1,y_1], \ldots, [x_6,y_6] \rangle$, 
\begin{align*}
c_{j,i}&= -b_{j,0} z_{0,i} -b_{j,1} z_{1,i} -b_{j,2} z_{2,i} &&(j \in \{4,5,6\},\ i \in \{1,2,3\}), \\
\widetilde{l}_{j,i}&=
(l_{j,i} +c_{4,i} l_{j,4} +c_{5,i} l_{j,5} +c_{6,i} l_{j,6}) h^{-1} &&(j \in \{1,\ldots,6\},\ i \in \{1,2,3\}), \\
\widetilde{l}_{j,0}&
=(l_{j,0}-q_0b_{4,0} l_{j,4}-q_0b_{5,0} l_{j,5}-q_0b_{6,0} l_{j,6}) h^{-1} &&(j \in \{1,\ldots,6\}), 
\end{align*}
and 
$A=(a_{j,(i,i')})_{1 \le j \le 6 \atop (i,i') \in \mathcal{I}}$ 
is written in the form 
\[
\left(
\begin{smallmatrix}
-\widetilde{l}_{1,0}&
&
&
\widetilde{l}_{1,2} &
\widetilde{l}_{1,3} & 
\\
&
-\widetilde{l}_{2,0}&
&
-\widetilde{l}_{2,1} &
&
\widetilde{l}_{2,3} \\
&
&
-\widetilde{l}_{3,0}&
&
-\widetilde{l}_{3,1}&
-\widetilde{l}_{3,2}\\
- q_0 b_{4,0}\widetilde{l}_{4,1} - c_{4,1}\widetilde{l}_{4,0}&
- q_0 b_{4,0}\widetilde{l}_{4,2} - c_{4,2}\widetilde{l}_{4,0}&
- q_0 b_{4,0}\widetilde{l}_{4,3} - c_{4,3}\widetilde{l}_{4,0}&
c_{4,1} \widetilde{l}_{4,2}-c_{4,2} \widetilde{l}_{4,1}&
c_{4,1} \widetilde{l}_{4,3}-c_{4,3} \widetilde{l}_{4,1}&
c_{4,2} \widetilde{l}_{4,3}-c_{4,3} \widetilde{l}_{4,2}\\
- q_0 b_{5,0}\widetilde{l}_{5,1} - c_{5,1}\widetilde{l}_{5,0}&
- q_0 b_{5,0}\widetilde{l}_{5,2} - c_{5,2}\widetilde{l}_{5,0}&
- q_0 b_{5,0}\widetilde{l}_{5,3} - c_{5,3}\widetilde{l}_{5,0}&
c_{5,1} \widetilde{l}_{5,2}-c_{5,2} \widetilde{l}_{5,1}&
c_{5,1} \widetilde{l}_{5,3}-c_{5,3} \widetilde{l}_{5,1}&
c_{5,2} \widetilde{l}_{5,3}-c_{5,3} \widetilde{l}_{5,2}\\
- q_0 b_{6,0}\widetilde{l}_{6,1} - c_{6,1}\widetilde{l}_{6,0}&
- q_0 b_{6,0}\widetilde{l}_{6,2} - c_{6,2}\widetilde{l}_{6,0}&
- q_0 b_{6,0}\widetilde{l}_{6,3} - c_{6,3}\widetilde{l}_{6,0}&
c_{6,1} \widetilde{l}_{6,2}-c_{6,2} \widetilde{l}_{6,1}&
c_{6,1} \widetilde{l}_{6,3}-c_{6,3} \widetilde{l}_{6,1}&
c_{6,2} \widetilde{l}_{6,3}-c_{6,3} \widetilde{l}_{6,2}
\end{smallmatrix}
\right), 
\]
where $\mathcal{I}=\{(0,1),(0,2),(0,3),(1,2),(1,3),(2,3)\}$. 
By a symbolic computation, 
one can see that 
\[
(\det B)hA \equiv 
\begin{pmatrix}
-l_{1,0}d_{456} & 
0 & 
0 &
\widetilde{d}_{1,2456} & 
\widetilde{d}_{1,3456} & 
0 \\
0 & 
-l_{2,0}d_{456} & 
0 &
-\widetilde{d}_{2,1456} & 
0 &
\widetilde{d}_{2,3456} \\
0 & 
0 & 
-l_{3,0}d_{456} &
0 & 
-\widetilde{d}_{3,1456} &
-\widetilde{d}_{3,2456} \\ 
l_{4,0}d_{156} & 
l_{4,0}d_{256} & 
l_{4,0}d_{356} &
\widetilde{d}_{4,1256} & 
\widetilde{d}_{4,1356} &
\widetilde{d}_{4,2356} \\
-l_{5,0}d_{146} &  
-l_{5,0}d_{246} & 
-l_{5,0}d_{346} &
-\widetilde{d}_{5,1246} & 
-\widetilde{d}_{5,1346} &
-\widetilde{d}_{5,2346} \\ 
l_{6,0}d_{145} & 
l_{6,0}d_{245} & 
l_{6,0}d_{345} &
\widetilde{d}_{6,1245} & 
\widetilde{d}_{6,1345} &
\widetilde{d}_{6,2345} 
\end{pmatrix}
\pmod{p}, 
\]
where 
\[
d_{i_1i_2i_3}=\left|
\begin{matrix}
z_{0,{i_1}} & z_{0,{i_2}} & z_{0,{i_3}} \\
z_{1,{i_1}} & z_{1,{i_2}} & z_{1,{i_3}} \\
z_{2,{i_1}} & z_{2,{i_2}} & z_{2,{i_3}} \\
\end{matrix}
\right| , 
\]
and 
\[
\widetilde{d}_{j,i_1i_2i_3i_4}
=\left|\begin{matrix}
l_{j,i_1} & l_{j,i_2} & l_{j,i_3} & l_{j,i_4} \\
z_{0,{i_1}} & z_{0,{i_2}} & z_{0,{i_3}} & z_{0,{i_4}} \\
z_{1,{i_1}} & z_{1,{i_2}} & z_{1,{i_3}} & z_{1,{i_4}} \\
z_{2,{i_1}} & z_{2,{i_2}} & z_{2,{i_3}} & z_{2,{i_4}} 
\end{matrix}\right|
= l_{j,i_1}d_{i_2i_3i_4}-l_{j,i_2}d_{i_1i_3i_4}+l_{j,i_3}d_{i_1i_2i_4}-l_{j,i_4}d_{i_1i_2i_3} .
\]
Note that $d_{456}=\det B$.

\subsection{Comparison with previous criteria.} 

In the situations of Propositions \ref{prop:4t0} and \ref{prop:6t0}, 
assuming that $k/\mathbb Q$ is abelian, 
we compare the propositions with the conditions of \cite{Oka**}. 
By Theorem \ref{thm:Ghyoji}, $\widetilde{G}/[\widetilde{G},\widetilde{G}] \simeq \mathrm{Gal}(\widetilde{k}/k)=\langle \overline{x_1},\ldots,\overline{x_{r_2+1}} \rangle$, $x_i[\widetilde{G},\widetilde{G}] \leftrightarrow \overline{x_i}$ ($r_2 \in \{2,3\}$). 
The decomposition group of $\mathfrak{p}_i$ in $\mathrm{Gal}(\widetilde{k}/k)$ is $D_i=\langle \overline{x_i},\overline{y_i} \rangle$. 
Under the isomorphism $\mathrm{Gal}(\widetilde{k}/k) \simeq J/Y$, 
$\overline{x_i}$ and $\overline{y_i}^h$ correspond to $\alpha_iY$ and $\varpi_iY$ respectively. 
If $g(\varpi_i)=\varpi_j$ for $g \in \mathrm{Gal}(k/\mathbb Q)$, 
then $g(\alpha_iY)=\alpha_jY$ and $g(\varpi_iY)=\varpi_jY$, 
and hence $g(\overline{x_i})=\overline{x_j}$ and $g(\overline{y_i})=\overline{y_j}$. 
Suppose that $g(\varpi_1)=\varpi_2$, and that $g(\varpi_2)=\varpi_3$ if $r_2=3$. 
Put $\gamma = \overline{x_1}$ 
and $x=\overline{y_1}(\overline{x_1})^{-\sum_{i=1}^{r_2+1} \widetilde{l}_{1,i}}$. 
Then $D_1=\langle \gamma,x \rangle$. 
Since $\overline{x_i}|_{k_{\mathrm{cyc}}} = \gamma|_{k_{\mathrm{cyc}}}$ by the choices of $\alpha_i$, 
and since $\overline{y_1} = \prod_{i=1}^{r_2+1} \overline{x_i}^{\widetilde{l}_{1,i}}$ by \eqref{eq:yj}, 
we have $\overline{y_1}|_{k_{\mathrm{cyc}}}=\gamma|_{k_{\mathrm{cyc}}}^{\sum_{i=1}^{r_2+1} \widetilde{l}_{1,i}}$, 
i.e., $x \in \mathrm{Gal}(\widetilde{k}/k_{\mathrm{cyc}})$. 
Put $y=g(x)=\overline{y_2}(\overline{x_2})^{-\sum_{i=1}^{r_2+1} \widetilde{l}_{1,i}}$, 
and put $z=g(y)=\overline{y_3}(\overline{x_3})^{-\sum_{i=1}^{4} \widetilde{l}_{1,i}}$ if $r_2=3$. 
Note that $\mathrm{Gal}(\widetilde{k}/k_{\mathrm{cyc}})=\langle \overline{x_2}\gamma^{-1},\ldots,\overline{x_{r_2+1}}\gamma^{-1} \rangle$. 

Now we assume that $X(k_{\mathrm{cyc}}) \simeq \mathbb Z_p^{r_2}$ (i.e., $\mathrm{rank}\,\overline{A^+}=\frac{r_2(r_2+1)}{2}$) as in \cite{Oka**}. 
Then 
$\mathrm{Gal}(\widetilde{k}/k_{\mathrm{cyc}})=\langle x,y \rangle$ if $r_2=2$, 
and 
$\mathrm{Gal}(\widetilde{k}/k_{\mathrm{cyc}})=\langle x,y,z \rangle$ if $r_2=3$
(see \cite[Lemmas 3.1.3, 3.2.1 and 4.1.4]{Oka**}). 
Since $\overline{y_j} = \prod_{i=1}^{r_2+1} \overline{x_i}^{\widetilde{l}_{j,i}}$ by \eqref{eq:yj}, 
we have 
\begin{align*}
x &=(\overline{x_2}\gamma^{-1})^{\widetilde{l}_{1,2}} \cdots (\overline{x_{r_2+1}}\gamma^{-1})^{\widetilde{l}_{1,{r_2+1}}}, 
\\
y &=
(\overline{x_2}\gamma^{-1})^{\widetilde{l}_{2,2}-\sum_{i=1}^{r_2+1} \widetilde{l}_{1,i}}
(\overline{x_3}\gamma^{-1})^{\widetilde{l}_{2,3}} \cdots (\overline{x_{r_2+1}}\gamma^{-1})^{\widetilde{l}_{2,{r_2+1}}}
\gamma^{\sum_{i=1}^{r_2+1} (\widetilde{l}_{2,i}-\widetilde{l}_{1,i})} , 
\\
z &=
(\overline{x_2}\gamma^{-1})^{\widetilde{l}_{3,2}}
(\overline{x_3}\gamma^{-1})^{\widetilde{l}_{3,3}-\sum_{i=1}^{4} \widetilde{l}_{1,i}}
(\overline{x_4}\gamma^{-1})^{\widetilde{l}_{3,4}}
\gamma^{\sum_{i=1}^{4} (\widetilde{l}_{3,i}-\widetilde{l}_{1,i})} . 
\end{align*}
Since $y,z \in \mathrm{Gal}(\widetilde{k}/k_{\mathrm{cyc}})$, 
we have 
\begin{align}\label{eq:l2l1}
\sum_{i=1}^{r_2+1} (\widetilde{l}_{2,i}-\widetilde{l}_{1,i})=0 , 
\end{align}
and
\begin{align}\label{eq:l3l1}
\sum_{i=1}^{4} (\widetilde{l}_{3,i}-\widetilde{l}_{1,i})=0
\end{align}
if $r_2=3$. 
If $r_2=2$, 
$(x,y)=(\overline{x_2}\gamma^{-1}, \overline{x_3}\gamma^{-1})M_2$ with 
\[
M_2=\begin{pmatrix} 
\widetilde{l}_{1,2} & \widetilde{l}_{2,2}-\sum_{i=1}^{3} \widetilde{l}_{1,i} \\
\widetilde{l}_{1,3} & \widetilde{l}_{2,3}
\end{pmatrix}
\]
satisfying $\det M_2 \not\equiv 0 \pmod{p}$. 
By a symbolic computation, 
\[
h^2\det M_2=D z_{1,1}^{-2}=
\bigg\{\begin{array}{ll}
2(l_{1,2}+l_{1,3})(l_{1,3}+l_{1,4}) & \text{ if } \mathrm{Gal}(k/\mathbb Q) \simeq V_4, \\
(l_{1,2}+l_{1,3})^2+(l_{1,3}+l_{1,4})^2 & \text{ if } \mathrm{Gal}(k/\mathbb Q) \simeq C_4. 
\end{array}
\]
If $r_2=3$, 
$(x,y,z)=(\overline{x_2}\gamma^{-1}, \overline{x_3}\gamma^{-1}, \overline{x_4}\gamma^{-1})M_3$ with 
\[
M_3=\begin{pmatrix} 
\widetilde{l}_{1,2} & \widetilde{l}_{2,2}-\sum_{i=1}^{4} \widetilde{l}_{1,i} & \widetilde{l}_{3,2} \\
\widetilde{l}_{1,3} & \widetilde{l}_{2,3} & \widetilde{l}_{3,3}-\sum_{i=1}^{4} \widetilde{l}_{1,i} \\
\widetilde{l}_{1,4} & \widetilde{l}_{2,4} & \widetilde{l}_{3,4}
\end{pmatrix}
\]
satisfying $\det M_3 \not\equiv 0 \pmod{p}$. 
By a symbolic computation, $h^3\det M_3 =m_1m_2$. 
Hence the conditions of $X(k_{\mathrm{cyc}}) \simeq \mathbb Z_p^{r_2}$ in Propositions \ref{prop:4t0}, \ref{prop:6t0} are certainly satisfied. 

The conditions of \cite{Oka**} are described by $a,b,c \in \mathbb Z_p$ satisfying:
\begin{itemize}
\item
$g(\gamma)=\gamma x^a y^b$ if $\mathrm{Gal}(k/\mathbb Q) \simeq C_4$, 
\item
$g(\gamma)=\gamma x^a y^{-a}$ and $\iota(\gamma)=\gamma x^{a+b} y^{b-a}$ if $\mathrm{Gal}(k/\mathbb Q) \simeq V_4$, 
\item
$g(\gamma)=\gamma x^a y^b z^c$ if $\mathrm{Gal}(k/\mathbb Q) \simeq C_6$. 
\end{itemize}

\begin{proposition}\label{prop:abc}
Under the settings above, $a,b$ and $c$ are described by arithmetic linking numbers as follows: 
\begin{itemize}
\item
If $\mathrm{Gal}(k/\mathbb Q) \simeq C_4$, 
\[
a=\frac{h(l_{1,2} + l_{1,3})}{(l_{1,2}+l_{1,3})^2+(l_{1,3}+l_{1,4})^2}, \quad 
b=\frac{-h(l_{1,3} + l_{1,4})}{(l_{1,2}+l_{1,3})^2+(l_{1,3}+l_{1,4})^2} . 
\]
Then $a+b \not\equiv 0 \pmod{p}$ if and only if $l_{1,2}-l_{1,4} \not\equiv 0 \pmod{p}$. 
\item
If $\mathrm{Gal}(k/\mathbb Q) \simeq V_4$, 
\[
a = \frac{h}{2(l_{1,2} + l_{1,3})}, \quad 
b = \frac{h}{2(l_{1,3} + l_{1,4})} . 
\]
Then $a-b \not\equiv 0 \pmod{p}$ if and only if $l_{1,2}-l_{1,4} \not\equiv 0 \pmod{p}$. 
\item
If $\mathrm{Gal}(k/\mathbb Q) \simeq C_6$, 
we have $m_1m_2={a'}^3+{b'}^3+{c'}^3-3a'b'c'$ and 
\[
a = \frac{h( {a'}^2 - b'c' )}{m_1m_2}, \quad 
b = \frac{-h( {c'}^2 - a'b' )}{m_1m_2}, \quad 
c = \frac{h( {b'}^2 - a'c' )}{m_1m_2}, 
\]
where 
$a' = l_{1,2} + l_{1,3} + l_{1,4}$, 
$b' = -(l_{1,3} + l_{1,4} + l_{1,5})$, 
$c' = l_{1,4} + l_{1,5} + l_{1,6}$. 
Then $a+b \not\equiv 0 \pmod{p}$ if and only if $l_{1,2}+l_{1,3}-l_{1,5}-l_{1,6} \not\equiv 0 \pmod{p}$. 
Moreover, 
$(a+b)((a+b)^2+3c^2) \not\equiv 0 \pmod{p}$ if and only if $(l_{1,2} - l_{1,5})^3 +(l_{1,3} - l_{1,6})^3 \not\equiv 0 \pmod{p}$. 
\end{itemize}
\end{proposition}

\begin{proof}
If $\mathrm{Gal}(k/\mathbb Q) \simeq C_4$, using \eqref{eq:l2l1}, we obtain 
\begin{align*}
\overline{x_1}^{-1}\overline{x_2}
&=\overline{x_1}^{-1}g(\overline{x_1})=\gamma^{-1}g(\gamma)=x^a y^b 
\\
&=  
\overline{x_1}^{a(
-\widetilde{l}_{1,2}-\widetilde{l}_{1,3}
)+b\widetilde{l}_{2,1}} 
\overline{x_2}^{a\widetilde{l}_{1,2}+b(\widetilde{l}_{2,2}-\sum_{i=1}^{3} \widetilde{l}_{1,i})} 
\overline{x_3}^{a\widetilde{l}_{1,3}+b\widetilde{l}_{2,3}} . 
\end{align*}
If $\mathrm{Gal}(k/\mathbb Q) \simeq V_4$, using \eqref{eq:l2l1}, we obtain 
\begin{align*}
\overline{x_1}^{-1}\overline{x_2}
&=\overline{x_1}^{-1}g(\overline{x_1})=\gamma^{-1}g(\gamma)=x^a y^{-a} 
\\
&=  
\overline{x_1}^{a(-\widetilde{l}_{1,2}-\widetilde{l}_{1,3}-\widetilde{l}_{2,1})} 
\overline{x_2}^{a(\widetilde{l}_{1,2}-\widetilde{l}_{2,2}+\sum_{i=1}^{3} \widetilde{l}_{1,i})} 
\overline{x_3}^{a(\widetilde{l}_{1,3}-\widetilde{l}_{2,3})} , 
\\
\overline{x_2}^{-1}\overline{x_3}
&=g(\overline{x_1})^{-1}\iota(\overline{x_1})=g(\gamma)^{-1}\iota(\gamma)=x^b y^{b} 
\\
&=  
\overline{x_1}^{b(-\widetilde{l}_{1,2}-\widetilde{l}_{1,3}+\widetilde{l}_{2,1})} 
\overline{x_2}^{b(\widetilde{l}_{1,2}+\widetilde{l}_{2,2}-\sum_{i=1}^{3} \widetilde{l}_{1,i})} 
\overline{x_3}^{b(\widetilde{l}_{1,3}+\widetilde{l}_{2,3})} . 
\end{align*}
If $\mathrm{Gal}(k/\mathbb Q) \simeq C_6$, using \eqref{eq:l2l1} and \eqref{eq:l3l1}, we obtain 
\begin{align*}
\overline{x_1}^{-1}\overline{x_2}
=\ &\overline{x_1}^{-1}g(\overline{x_1})=\gamma^{-1}g(\gamma)=x^a y^b z^c
\\
=\ & 
\overline{x_1}^{a(\widetilde{l}_{1,1}-\sum_{i=1}^{4}\widetilde{l}_{1,i})+b\widetilde{l}_{2,1}+c\widetilde{l}_{3,1}}
\overline{x_2}^{a\widetilde{l}_{1,2}+b(\widetilde{l}_{2,2}-\sum_{i=1}^{4}\widetilde{l}_{1,i})+c\widetilde{l}_{3,2}}
\\
& \times 
\overline{x_3}^{a\widetilde{l}_{1,3}+b\widetilde{l}_{2,3}+c(\widetilde{l}_{3,3}-\sum_{i=1}^{4}\widetilde{l}_{1,i})}
\overline{x_4}^{a\widetilde{l}_{1,4}+b\widetilde{l}_{2,4}+c\widetilde{l}_{3,4}}. 
\end{align*}
By symbolic computations, we obtain $a,b,c$ as the solutions of simultaneous linear equations. 
If $\mathrm{Gal}(k/\mathbb Q) \simeq C_6$, we have 
$m_1m_2={a'}^3+{b'}^3+{c'}^3-3a'b'c'$, $(a+b)m_2 = h(a'-c')=h(l_{12}+l_{13}-l_{15}-l_{16})$ and 
\[
(a+b)((a+b)^2+3c^2)m_1^2m_2^2 = h^3((l_{1,2} - l_{1,5})^3 +(l_{1,3} - l_{1,6})^3). 
\]
Thus we obtain the statements. 
\end{proof}

\section{Examples}\label{sec:examples}

By using \texttt{PARI/GP} \cite{PARI}, we obtain the following examples.

\subsection{Quartic fields with $t=0$.} 

\begin{example}\label{ex:3V4t0}
Put $p=3$ and $k=\mathbb Q(\sqrt{7},\sqrt{-11})=\mathbb Q(\alpha)$ with $\alpha=\sqrt{7}+\sqrt{-11}$ 
satisfying $\alpha^4 + 8\alpha^2 + 324=0$. 
Then $\mathcal{O}_k=\mathbb Z[\sqrt{7},\frac{1+\sqrt{-11}}{2}]$, and $p\mathcal{O}_k$ splits into four prime ideals 
\[
\mathfrak{p}_1=3\mathcal{O}_k+\big(\sqrt{7}+\textstyle{\frac{1+\sqrt{-11}}{2}}\big)\mathcal{O}_k, 
\]
$\mathfrak{p}_2=g(\mathfrak{p}_1)$, $\mathfrak{p}_3=\iota(\mathfrak{p}_1)$ and $\mathfrak{p}_4=g\iota(\mathfrak{p}_1)$, 
where $g$ is a generator of $\mathrm{Gal}(k/\mathbb Q(\sqrt{-11}))$.  
The class number of $k^+=\mathbb Q(\sqrt{7})$ is $1$, 
and the findamental unit of $k$ is $\varepsilon_1=8+3\sqrt{7}$. 
Since $\mathfrak{p}_1\mathfrak{p}_3=(2-\sqrt{7})\mathcal{O}_k$ 
and $(2-\sqrt{7})^2=(4+\sqrt{7})\varepsilon_1^{-1}$, 
we have $\mathfrak{p}_1^2\mathfrak{p}_3^2=(4+\sqrt{7})\mathcal{O}_k$, 
and hence 
\[
\sqrt{7} \equiv 5 \pmod{\mathfrak{p}_1^2\mathfrak{p}_3^2}, \qquad 
\sqrt{7} \equiv 4 \pmod{\mathfrak{p}_2^2\mathfrak{p}_4^2}. 
\]
Since $\varepsilon_1^{p-1} =127+48\sqrt{7}$ satisfies 
\begin{align*}
\varepsilon_1^{p-1} &\equiv 127+48 \cdot 5 \equiv 7 \equiv (1+p)^{(p-1)1} \pmod{\mathfrak{p}_1^2\mathfrak{p}_3^2}, \\
\varepsilon_1^{p-1} &\equiv 127+48 \cdot 4 \equiv 4 \equiv (1+p)^{(p-1)2} \pmod{\mathfrak{p}_2^2\mathfrak{p}_4^2},  
\end{align*}
we have 
\[
Z=(z_{1,1},z_{1,2},z_{1,3},z_{1,4}) \equiv (1,2,1,2) \pmod{p}. 
\]
In particular, $B= (z_{1,4})$, and $b_{4,1}=z_{1,4}^{-1} \equiv 2 \pmod{p}$. 
Since $Cl_k \simeq \mathbb Z/4\mathbb Z$, $t=0$ and $h=4$. 
The ideal class $[\mathfrak{p}_1]$ is a generator of $Cl_k$. 
The principal ideals $\mathfrak{p}_i^h=\varpi_i\mathcal{O}_k$ have generators 
\[
\varpi_1 = 2-\sqrt{7} +\textstyle{\frac{1-\sqrt{-11}}{2}}, 
\]
$\varpi_2=g(\varpi_1)$, $\varpi_3=\iota(\varpi_1)$ and $\varpi_4=g\iota(\varpi_1)$, respectively. 
Since $\mathfrak{p}_1\mathfrak{p}_2=\frac{1-\sqrt{-11}}{2}\mathcal{O}_k$, 
we have $\mathfrak{p}_1^2\mathfrak{p}_2^2=\frac{5+\sqrt{-11}}{2}\mathcal{O}_k$, 
and hence 
\[
\sqrt{-11} \equiv 4 \pmod{\mathfrak{p}_1^2\mathfrak{p}_2^2}, \qquad 
\sqrt{-11} \equiv 5 \pmod{\mathfrak{p}_3^2\mathfrak{p}_4^2}. 
\]
Since $\varpi_1 = 7-\sqrt{7} -5\sqrt{-11} \pmod{9}$, we have 
\begin{align*}
\varpi_1^{p-1} &=(7-4 -5\cdot 4)^2 \equiv 1 \equiv (1+p)^{(p-1)0}\pmod{\mathfrak{p}_2^2}, \\
\varpi_1^{p-1} &=(7-5 -5\cdot 5)^2 \equiv 7 \equiv (1+p)^{(p-1)1} \pmod{\mathfrak{p}_3^2}, \\
\varpi_1^{p-1} &=(7-4 -5\cdot 5)^2 \equiv 7 \equiv (1+p)^{(p-1)1} \pmod{\mathfrak{p}_4^2}. 
\end{align*}
Then $(l_{1,2},l_{1,3},l_{1,4}) \equiv (0,1,1) \pmod{p}$, and hence we have 
\begin{itemize}
\item
$( l_{1,2} + l_{1,3})(l_{1,3} + l_{1,4}) \equiv 2 \not\equiv 0 \pmod{p}$, 
\item 
$l_{1,2}-l_{1,4} \equiv 2 \not\equiv 0 \pmod{p}$, 
\end{itemize}
and $(a,b) \equiv (2, 1) \pmod{p}$ by Proposition \ref{prop:abc}. 
Therefore $X(k_{\mathrm{cyc}}) \simeq \mathbb Z_p^2$ and $X(\widetilde{k}) \simeq 0$ by Proposition \ref{prop:4t0}. 
In fact, $\mathrm{rank}\,\overline{A} = 3$ for a $\mathbb F_3$-matrix 
\[
\overline{A}=(\overline{a_{j,(i,i')}})_{j,(i,i')} = 
\begin{pmatrix}
2 & 2 & 0 \\
2 & 0 & 2 \\
0 & 2 & 2 \\
2 & 2 & 2 
\end{pmatrix}. 
\]
\end{example}

\begin{example}\label{ex:3V4not0}
Suppose $p=3$, and put $k=\mathbb Q(\sqrt{7},\sqrt{-m})=\mathbb Q(\alpha)$ with $m \in \{17, 65\}$ and $\alpha=\sqrt{7}+\sqrt{-m}$. 
Then $\alpha^4+2(m-7)\alpha^2+(m+7)^2=0$, and $p\mathcal{O}_k=\mathfrak{p}_1\mathfrak{p}_2\mathfrak{p}_3\mathfrak{p}_4$ 
with  
$\mathfrak{p}_1=3\mathcal{O}_k+\big(\frac{\alpha}{2} + 1\big)\mathcal{O}_k$, 
$\mathfrak{p}_2=g(\mathfrak{p}_1)$, $\mathfrak{p}_3=\iota(\mathfrak{p}_1)$ and $\mathfrak{p}_4=g\iota(\mathfrak{p}_1)$, 
where $g$ is a generator of $\mathrm{Gal}(k/\mathbb Q(\sqrt{-m}))$.  
Put $\varepsilon_1=8+3\sqrt{7}$. 
Then $Z \equiv (1,2,1,2) \pmod{p}$ and $b_{4,1} \equiv 2 \pmod{p}$. 

Suppose $m=17$. 
Since $Cl_k \simeq \mathbb Z/20\mathbb Z$, $t=0$ and $h=20$. 
Then $\mathfrak{p}_1^h=\varpi_1\mathcal{O}_k$ with 
$\varpi_1 = \frac{1}{8}\alpha^3 + \frac{39}{4}\alpha^2 - 19\alpha + 11$, 
and $(l_{1,2},l_{1,3},l_{1,4}) \equiv (0,1,0) \pmod{p}$. 
Since $(l_{1,2}+l_{1,3})(l_{1,3}+l_{1,4}) \not\equiv 0 \pmod{p}$, 
$X(k_{\mathrm{cyc}}) \simeq \mathbb Z_p^2$ and 
$(a,b) \equiv (1, 1) \pmod{p}$ by Propositions \ref{prop:4t0} and \ref{prop:abc}. 
Moreover, since $l_{1,2}-l_{1,4} \equiv 0 \pmod{p}$, 
$X(\widetilde{k}) \not\simeq 0$ by Proposition \ref{prop:4t0}. 
In fact, $\mathrm{rank}\,\overline{A} = 2$ for 
$\overline{A}= 
\left(\begin{smallmatrix}
0 & 2 & 0 \\
1 & 0 & 2 \\
0 & 1 & 0 \\
2 & 0 & 1 
\end{smallmatrix}\right)$, 
and 
$\mathrm{rank}\,\overline{A^+} = 3$ for 
$\overline{A^+}=\left(\begin{smallmatrix}
\hbox{$\overline{A}$} \\ \begin{smallmatrix} 1 & 2 & 1 \end{smallmatrix}
\end{smallmatrix}\right)$. 

Suppose $m=65$. 
Since $Cl_k \simeq \mathbb Z/40\mathbb Z \times \mathbb Z/2\mathbb Z$, $t=0$ and $h=80$. 
Then $\mathfrak{p}_1^h=\varpi_1\mathcal{O}_k$ with 
$\varpi_1 = \frac{6664171}{2}\alpha^3 - 62398544\alpha^2 + 688638986\alpha - 4014380975$, 
and $(l_{1,2},l_{1,3},l_{1,4}) \equiv (1,1,2) \pmod{p}$. 
Since $(l_{1,2} + l_{1,3})(l_{1,3} + l_{1,4}) \equiv 0 \pmod{p}$, 
we have $X(k_{\mathrm{cyc}}) \not\simeq \mathbb Z_p^2$ and $X(\widetilde{k}) \not\simeq 0$ by Proposition \ref{prop:4t0}. 
In fact, $\mathrm{rank}\,\overline{A} = \mathrm{rank}\,\overline{A^+} = 2$ for 
$\overline{A}= 
\left(\begin{smallmatrix}
1 & 0 & 0 \\
2 & 0 & 0 \\
0 & 2 & 1 \\
0 & 1 & 2 
\end{smallmatrix}\right)$ 
and 
$\overline{A^+}=\left(\begin{smallmatrix}
\hbox{$\overline{A}$} \\ \begin{smallmatrix} 1 & 2 & 1 \end{smallmatrix}
\end{smallmatrix}\right)$. 

In both cases of $m=17$ and $m=65$, since $\mathrm{rank}\,\overline{A} = 2 = \frac{r_2(r_2+1)}{2}-1$, 
$X(\widetilde{k})$ is $\widetilde{\varLambda}$-cyclic by the property \ref{main:lmdrk} of $A$. 
\end{example}

\begin{example}\label{ex:3V4d4}
Put $p=3$ and $k=\mathbb Q(\sqrt{43},\sqrt{-m})$ with a squarefree integer $m \equiv 2 \pmod{3}$. 
Then $p\mathcal{O}_k=\mathfrak{p}_1\mathfrak{p}_2\mathfrak{p}_3\mathfrak{p}_4$ with four prime ideals $\mathfrak{p}_i$. 
Since the fundamental unit $\varepsilon_1=3482+531\sqrt{43}$ of $k^+$ satisfies 
$\varepsilon_1 \equiv -1 \pmod{p^2}$, 
we have $\overline{Z} = (0,0,0,0)$. 
Then $\widetilde{V}/(k^{\times})^p \simeq E_k/p \simeq \mathbb Z/p\mathbb Z$ and $d(\widetilde{G})=4$ by Proposition \ref{prop:Vis0}. 
In particular, $X(\widetilde{k}) \not\simeq 0$. 
\end{example}

\begin{example}\label{ex:3C4t0}
Put $p=3$ and $k=\mathbb Q(\alpha) \subset \mathbb Q(\zeta_{13},\sqrt{10})$ 
with $\alpha=\sqrt{\frac{-13+3\sqrt{13}}{2}10}=\sqrt{-65+15\sqrt{13}}$ 
satisfying $\alpha^4 + 130\alpha^2 + 1300=0$. 
Note that $\mathbb Q(\sqrt{\frac{-13+3\sqrt{13}}{2}}) \subset \mathbb Q(\zeta_{13})$, 
and that $\mathrm{Gal}(k/\mathbb Q) \simeq C_4$. 
The fundamental unit of $k$ is $\varepsilon_1=\frac{3+\sqrt{13}}{2}=\frac{\alpha^2+110}{30}$. 
Then $p\mathcal{O}_k=\mathfrak{p}_1\mathfrak{p}_2\mathfrak{p}_3\mathfrak{p}_4$ with 
$\mathfrak{p}_1
=3\mathcal{O}_k+\big( \frac{\alpha^3}{15} + \frac{\alpha^2}{30} + \frac{16\alpha}{3} + \frac{11}{3} \big)\mathcal{O}_k$, 
$\mathfrak{p}_2=g(\mathfrak{p}_1)$, $\mathfrak{p}_3=g^2(\mathfrak{p}_1)$ and $\mathfrak{p}_4=g^3(\mathfrak{p}_1)$, 
where $g$ is a generator of $\mathrm{Gal}(k/\mathbb Q)$ such that $g(\alpha)=\alpha \varepsilon_1=\sqrt{-65-15\sqrt{13}}$.  
One can see that 
$Z \equiv (1,2,1,2) \pmod{p}$, and $b_{4,1} \equiv 2 \pmod{p}$. 
Since $Cl_k \simeq \mathbb Z/34\mathbb Z \times \mathbb Z/2\mathbb Z$, $t=0$ and $h=68$. 
Then $\mathfrak{p}_1^h=\varpi_1\mathcal{O}_k$ with 
$\varpi_1 = 
\frac{1194376}{3}\alpha^3 - \frac{4274169}{2}\alpha^2 + \frac{30168356}{3}\alpha + 5545534$. 
One can see that $(l_{1,2},l_{1,3},l_{1,4}) \equiv (2, 2, 1) \pmod{p}$, 
$( l_{1,2} + l_{1,3})^2+(l_{1,3} + l_{1,4})^2 \equiv 1 \not\equiv 0 \pmod{p}$, $l_{1,2}-l_{1,4} \equiv 1 \not\equiv 0 \pmod{p}$, 
and $a \equiv 2 \pmod{p}$, $b \equiv 0 \pmod{p}$ by Proposition \ref{prop:abc}. 
By Proposition \ref{prop:4t0}, we have $X(k_{\mathrm{cyc}}) \simeq \mathbb Z_p^2$ and $X(\widetilde{k}) \simeq 0$. 
In fact, $\mathrm{rank}\,\overline{A} = 3$ for 
$\overline{A}= 
\left(\begin{smallmatrix}
2 & 0 & 0 \\
0 & 0 & 2 \\
0 & 1 & 2 \\
2 & 1 & 0 
\end{smallmatrix}\right)$. 
\end{example}

\subsection{Quartic fields with $t=1$.} 

\begin{example}\label{ex:3V4t1}
Put $p=3$ and $k=\mathbb Q(\sqrt{7},\sqrt{-23})=\mathbb Q(\alpha)$ with $\alpha=\sqrt{7}+\sqrt{-23}$ 
satisfying $\alpha^4 + 32\alpha^2 + 900=0$. 
Then $p\mathcal{O}_k=\mathfrak{p}_1\mathfrak{p}_2\mathfrak{p}_3\mathfrak{p}_4$ 
with  
$\mathfrak{p}_1=3\mathcal{O}_k+\big( \sqrt{7} + \frac{1+\sqrt{-23}}{2} \big)\mathcal{O}_k$, 
$\mathfrak{p}_2=g(\mathfrak{p}_1)$, $\mathfrak{p}_3=\iota(\mathfrak{p}_1)$ and $\mathfrak{p}_4=g\iota(\mathfrak{p}_1)$, 
where $g$ is a generator of $\mathrm{Gal}(k/\mathbb Q(\sqrt{-23}))$.  
The ideal class $[\mathfrak{p}_1]$ generates $Cl_k \simeq \mathbb Z/24\mathbb Z$. 
We have $t=1$ and $h=8$. 
Put
$\mathfrak{a}_0=29\mathcal{O}_k+ \alpha\mathcal{O}_k$, 
which is a prime ideal of $k$ lying over $29$. 
The order of $[\mathfrak{a}_0] \in Cl_k$ is $q_0=12$. 
In particular, $Cl_k/p = \langle [\mathfrak{a}_0]Cl_k^p \rangle$. 
Then $\mathfrak{a}_0^{q_0}=\varepsilon_0\mathcal{O}_k$ with 
$\varepsilon_0 = 
\frac{88779+287\sqrt{-23}}{2} + \frac{-34839+891\sqrt{-23}}{2}\sqrt{7}$. 
The fundamental unit of $k$ is $\varepsilon_1=8+3\sqrt{7}$. 
Then 
\[
Z=\begin{pmatrix}
z_{0,1} & z_{0,2} & z_{0,3} & z_{0,4} \\
z_{1,1} & z_{1,2} & z_{1,3} & z_{1,4} 
\end{pmatrix} 
\equiv
\begin{pmatrix}
2 & 2 & 1 & 1 \\
1 & 2 & 1 & 2 
\end{pmatrix}
\pmod{p}, 
\]
and hence 
\[
\begin{pmatrix}
b_{3,0} & b_{3,1} \\ 
b_{4,0} & b_{4,1}  
\end{pmatrix}
=B^{-1}
\equiv
\begin{pmatrix}
2 & 2 \\ 
2 & 1  
\end{pmatrix}
\pmod{p} . 
\]
Since $\mathfrak{p}_i^h=\mathfrak{a}_0^{l_{i,0}} \varpi_i\mathcal{O}_k$ for $i \in \{1,2,3,4\}$ 
with 
$(l_{1,0},l_{2,0},l_{3,0},l_{4,0})=(4, 4, 8, 8)$ and 
\begin{align*}
\varpi_1 &= \big(103404-28948\sqrt{-23} + (82468-10091\sqrt{-23})\sqrt{7} \big) 29^{-4}, \\ 
\varpi_2 &= \big({-32004}+13556\sqrt{-23} + (-77809+5996\sqrt{-23})\sqrt{7} \big) 29^{-4}, \\
\varpi_3 &= \big(\textstyle{\frac{-4780928131+2975830185\sqrt{-23}}{2} + \frac{-3178361441+1135888013\sqrt{-23}}{2}\sqrt{7}} \big)29^{-8}, \\
\varpi_4 &= \big(\textstyle{\frac{99743021-2164820295\sqrt{-23}}{2} + \frac{2096910929-767151197\sqrt{-23}}{2}\sqrt{7}}\big)29^{-8} , 
\end{align*}
we have 
\[
(l_{i,j})_{1 \le i,j \le 4}=
\begin{pmatrix}
0 & l_{1,2} & l_{1,3} & l_{1,4} \\
l_{2,1} & 0 & l_{2,3} & l_{2,4} \\
l_{3,1} & l_{3,2} & 0 & l_{3,4} \\
l_{4,1} & l_{4,2} & l_{4,3} & 0 
\end{pmatrix}
\equiv
\begin{pmatrix}
0 & 2 & 0 & 1 \\
1 & 0 & 2 & 2 \\
2 & 0 & 0 & 2 \\ 
1 & 1 & 1 & 0 
\end{pmatrix}
\pmod{p}, 
\]
and 
\[
A = (a_{j,(i,i')})_{j,(i,i')} \equiv
\begin{pmatrix}
1& 0& 1 \\
0& 1& 0 \\
0& 2& 1 \\
2& 0& 1 
\end{pmatrix}
\pmod{p}. 
\]
Since $\mathrm{rank}\,\overline{A} = 3$, we have $X(\widetilde{k}) \simeq 0$. 
\end{example}

\begin{example}\label{ex:5C4t1}
Put $p=5$ and $k=\mathbb Q(\alpha) \subset \mathbb Q(\zeta_{29},\sqrt{3})$ 
with $\alpha=\sqrt{\frac{-29+5\sqrt{29}}{2} 3}$ 
satisfying $\alpha^4 + 87\alpha^2 + 261=0$. 
Note that $\mathbb Q(\sqrt{\frac{-29+5\sqrt{29}}{2}}) \subset \mathbb Q(\zeta_{29})$, 
and that $\mathrm{Gal}(k/\mathbb Q) \simeq C_4$. 
Then $p\mathcal{O}_k=\mathfrak{p}_1\mathfrak{p}_2\mathfrak{p}_3\mathfrak{p}_4$ with 
$\mathfrak{p}_1
=5\mathcal{O}_k+\Big( 
\frac{2}{15}\alpha^3 + \frac{1}{15}\alpha^2 + \frac{44}{5}\alpha + \frac{17}{5} \Big)\mathcal{O}_k$, 
$\mathfrak{p}_2=g(\mathfrak{p}_1)$, $\mathfrak{p}_3=g^2(\mathfrak{p}_1)$ and $\mathfrak{p}_4=g^3(\mathfrak{p}_1)$, 
where $g$ is a generator of $\mathrm{Gal}(k/\mathbb Q)$ such that $g(\alpha)=-\sqrt{\frac{-29-5\sqrt{29}}{2} 3}$.  
The ideal class $[\mathfrak{p}_1]$ generates $Cl_k \simeq \mathbb Z/10\mathbb Z \times \mathbb Z/2\mathbb Z$. 
We have $t=1$ and $h=4$. 
Put
$\mathfrak{a}_0=23\mathcal{O}_k+(\alpha-9)\mathcal{O}_k$, 
which is a prime ideal of $k$ lying over $23$. 
The order of $[\mathfrak{a}_0] \in Cl_k$ is $q_0=10$. 
In particular, $Cl_k/p = \langle [\mathfrak{a}_0]Cl_k^p \rangle$. 
Then $\mathfrak{a}_0^{q_0}=\varepsilon_0\mathcal{O}_k$ with 
$\varepsilon_0 = 
-\frac{41}{15}\alpha^3 + \frac{322}{15}\alpha^2 - \frac{3627}{5}\alpha + \frac{2954}{5}$. 
The fundamental unit of $k$ is $\varepsilon_1=\frac{5-\sqrt{29}}{2}=-\frac{\alpha^2+6}{15}$. 
Then 
$Z=\left(\begin{smallmatrix}
z_{0,1} & z_{0,2} & z_{0,3} & z_{0,4} \\
z_{1,1} & z_{1,2} & z_{1,3} & z_{1,4} 
\end{smallmatrix}\right)
\equiv
\left(\begin{smallmatrix}
3 & 3 & 4 & 0 \\
3 & 2 & 3 & 2
\end{smallmatrix}\right)
\pmod{p}$, 
and hence 
$\left(\begin{smallmatrix}
b_{3,0} & b_{3,1} \\ 
b_{4,0} & b_{4,1}  
\end{smallmatrix}\right)
=B^{-1}
\equiv
\left(\begin{smallmatrix}
4 & 0 \\ 
4 & 3  
\end{smallmatrix}\right)
\pmod{p}$. 
Since $\mathfrak{p}_i^h=\mathfrak{a}_0^{l_{i,0}} \varpi_i\mathcal{O}_k$ for $i \in \{1,2,3,4\}$ 
with 
$(l_{1,0},l_{2,0},l_{3,0},l_{4,0})=(4, 2, 6, 8)$ and 
\begin{align*}
\varpi_1 &= \big(\textstyle{ \frac{4892}{15}\alpha^3 + \frac{16957}{15}\alpha^2 + \frac{91424}{5}\alpha + \frac{40034}{5} }\big) 23^{-4}, \\ 
\varpi_2 &= \big(\textstyle{ {-\frac{19}{15}\alpha^3} -\frac{194}{15}\alpha^2 - \frac{673}{5}\alpha - \frac{5758}{5} }\big) 23^{-2}, \\ 
\varpi_3 &= \big(\textstyle{ \frac{114234}{5}\alpha^3 - \frac{1468532}{15}\alpha^2 + \frac{7497169}{5}\alpha - \frac{51582334}{5} }\big) 23^{-6}, \\ 
\varpi_4 &= \big(\textstyle{ \frac{141541772}{15}\alpha^3 + \frac{221625902}{15}\alpha^2 + \frac{3865644704}{5}\alpha + \frac{4691851209}{5} }\big) 23^{-8}, 
\end{align*}
we have 
$(l_{i,j})_{1 \le i,j \le 4}
\equiv
\left(\begin{smallmatrix}
0 & 3 & 2 & 1 \\
1 & 0 & 4 & 4 \\
1 & 4 & 0 & 4 \\
1 & 0 & 2 & 0
\end{smallmatrix}\right)
\pmod{p}$, 
and 
$A = (a_{j,(i,i')})_{j,(i,i')} \equiv
\left(\begin{smallmatrix}
4 & 0 & 4 \\
0 & 2 & 4 \\
3 & 3 & 0 \\
2 & 1 & 0
\end{smallmatrix}\right)
\pmod{p}$. 
Since $\mathrm{rank}\,\overline{A} = 3$, we have $X(\widetilde{k}) \simeq 0$. 
\end{example}

\subsection{Sextic fields with $t=0$.} 

\begin{example}\label{ex:3C6t0}
Let $k^+$ be the unique cubic subfield of $\mathbb Q(\zeta_{61})$. 
Put $p=3$ and $k=k^+(\sqrt{-2})$. 
Then $k=\mathbb Q(\alpha)$ with $\alpha$ satisfying 
$\alpha^6 - 2\alpha^5 - 33\alpha^4 + 50\alpha^3 + 398\alpha^2 - 396\alpha + 1089 =0$,  
and $p\mathcal{O}_k=\prod_{i=1}^{6}\mathfrak{p}_i$ with 
$\mathfrak{p}_1=3\mathcal{O}_k+(\alpha -2)\mathcal{O}_k$, 
$\mathfrak{p}_i=g^{i-1}(\mathfrak{p}_1)$ for $i \in \{2,\ldots,6\}$, 
where $g$ is a generator of $\mathrm{Gal}(k/\mathbb Q)$ such that 
$g(\alpha)=-\frac{4}{\delta}\alpha^5 - \frac{25}{\delta}\alpha^4 + \frac{440}{3\delta}\alpha^3 + \frac{1010}{\delta}\alpha^2 - \frac{10265}{3\delta}\alpha - \frac{7425}{\delta}$ 
for $\delta=11 \cdot 241$. 
The fundamental units of $k$ are 
\begin{align*}
\varepsilon_1 &= 
\textstyle{-\frac{232}{9\delta}\alpha^5 + \frac{1201}{9\delta}\alpha^4 + \frac{4972}{9\delta}\alpha^3 - \frac{36856}{9\delta}\alpha^2 + \frac{4852}{\delta}\alpha - \frac{8085}{\delta}} , \\
\varepsilon_2 &= 
\textstyle{-\frac{46}{9\delta}\alpha^5 - \frac{1613}{9\delta}\alpha^4 + \frac{3454}{9\delta}\alpha^3 + \frac{59333}{9\delta}\alpha^2 - \frac{5894}{\delta}\alpha - \frac{92994}{\delta}} .
\end{align*}
One can see that 
$Z=(z_{i,j})_{i,j} \equiv 
\left(\begin{smallmatrix}
2 & 0 & 1 & 2 & 0 & 1 \\
1 & 2 & 0 & 1 & 2 & 0 
\end{smallmatrix}\right)
\pmod{p}$, 
and hence 
$\left(\begin{smallmatrix}
b_{5,1} & b_{5,2} \\ 
b_{6,1} & b_{6,2}  
\end{smallmatrix}\right)
=B^{-1}
\equiv
\left(\begin{smallmatrix}
0 & 2 \\ 
1 & 0  
\end{smallmatrix}\right)
\pmod{p}$. 
Since $Cl_k \simeq \mathbb Z/43\mathbb Z$, $t=0$ and $h=43$. 
Then $\mathfrak{p}_1^h=\varpi_1\mathcal{O}_k$ with 
\[
\varpi_1 = 
\textstyle{\frac{29881}{3\delta}\alpha^5 - \frac{93587}{3\delta}\alpha^4 - \frac{516835}{3\delta}\alpha^3 + \frac{958130}{3\delta}\alpha^2 - \frac{953187}{\delta}\alpha + \frac{1625195}{\delta}} , 
\]
and $\varpi_i=g^{i-1}(\varpi_1)$ for $i \in \{2,\ldots,6\}$. 
Since 
$(l_{1,2},l_{1,3},l_{1,4},l_{1,5},l_{1,6}) \equiv (2, 0, 0, 1, 2) \pmod{p}$, 
we have 
$m_1 \equiv m_2 \equiv 1 \not\equiv 0 \pmod{p}$, 
$(l_{1,2}-l_{1,5})^3+(l_{1,3}-l_{1,6})^3 \equiv 2 \not\equiv 0 \pmod{p}$, 
and $(a,b,c) \equiv (1, 1, 1) \pmod{p}$ by Proposition \ref{prop:abc}. 
By Proposition \ref{prop:6t0}, $X(k_{\mathrm{cyc}}) \simeq \mathbb Z_p^3$ and $X(\widetilde{k}) \simeq 0$. 
In fact, $\mathrm{rank}\,\overline{A} = 6$ for 
\[
A = (a_{j,(i,i')})_{j,(i,i')} \equiv
\begin{pmatrix}
1& 1& 0& 0& 0& 0 \\
0& 0& 0& 1& 1& 0 \\
0& 2& 0& 1& 0& 2 \\
0& 0& 1& 0& 1& 1 \\
2& 2& 2& 1& 2& 1 \\
1& 1& 2& 1& 2& 0 
\end{pmatrix}
\pmod{p}. 
\]
\end{example}

\begin{example}\label{ex:7C6t0}
Let $k^+$ be the unique cubic subfield of $\mathbb Q(\zeta_{19})$. 
Put $p=7$ and $k=k^+(\sqrt{-3})$. 
Then $k=\mathbb Q(\alpha)$ with $\alpha$ satisfying 
$\alpha^6 - 2\alpha^5 - 2\alpha^4 + 14\alpha^3 + 55\alpha^2 - 192\alpha + 343 =0$,  
and $p\mathcal{O}_k=\prod_{i=1}^{6}\mathfrak{p}_i$ with 
$\mathfrak{p}_1=7\mathcal{O}_k+\alpha\mathcal{O}_k$, 
$\mathfrak{p}_i=g^{i-1}(\mathfrak{p}_1)$ for $i \in \{2,\ldots,6\}$, 
where $g$ is a generator of $\mathrm{Gal}(k/\mathbb Q)$ such that 
$g(\alpha)=-\frac{272}{\delta}\alpha^5 - \frac{305}{\delta}\alpha^4 + \frac{860}{\delta}\alpha^3 + \frac{1150}{\delta}\alpha^2 - \frac{24925}{\delta}\alpha + \frac{6300}{\delta}$ for $\delta=7 \cdot 1699$.  
The fundamental units of $k$ are 
\begin{align*}
\varepsilon_1 &= \textstyle{\frac{-272}{\delta}\alpha^5 - \frac{305}{\delta}\alpha^4 + \frac{860}{\delta}\alpha^3 + \frac{1150}{\delta}\alpha^2 - \frac{13032}{\delta}\alpha + \frac{6300}{\delta}} ,\\
\varepsilon_2 &= \textstyle{\frac{-162}{\delta}\alpha^5 - \frac{444}{\delta}\alpha^4 + \frac{862}{\delta}\alpha^3 + \frac{2259}{\delta}\alpha^2 - \frac{14058}{\delta}\alpha - \frac{19684}{\delta}} .
\end{align*}
One can see that 
$Z=(z_{i,j})_{i,j} \equiv 
\left(\begin{smallmatrix}
1 & 1 & 5 & 1 & 1 & 5 \\
1 & 5 & 1 & 1 & 5 & 1
\end{smallmatrix}\right)
\pmod{p}$, 
and hence 
$\left(\begin{smallmatrix}
b_{5,1} & b_{5,2} \\ 
b_{6,1} & b_{6,2}  
\end{smallmatrix}\right)
=B^{-1}
\equiv
\left(\begin{smallmatrix}
2 & 4 \\ 
4 & 2  
\end{smallmatrix}\right)
\pmod{p}$. 
Since $Cl_k \simeq \mathbb Z/3\mathbb Z$, $t=0$ and $h=3$. 
Then $\mathfrak{p}_1^h=\varpi_1\mathcal{O}_k$ with $\varpi_1=-\alpha$, 
and $\varpi_i=g^{i-1}(\varpi_1)$ for $i \in \{2,\ldots,6\}$. 
Since 
$(l_{1,2},l_{1,3},l_{1,4},l_{1,5},l_{1,6}) \equiv (1, 0, 2, 6, 1) \pmod{p}$, 
we have 
$m_1 \equiv 4 \pmod{p}$, 
$m_2 \equiv 6 \pmod{p}$, 
$(l_{1,2}-l_{1,5})^3+(l_{1,3}-l_{1,6})^3 \equiv 0 \pmod{p}$, 
and $(a,b,c) \equiv (4, 0, 2) \pmod{p}$ by Proposition \ref{prop:abc}. 
By Proposition \ref{prop:6t0}, $X(k_{\mathrm{cyc}}) \simeq \mathbb Z_p^3$ but $X(\widetilde{k}) \not\simeq 0$. 
In fact, $\mathrm{rank}\,\overline{A} = 5 = \frac{r_2(r_2+1)}{2}-1$ for 
\[
A = (a_{j,(i,i')})_{j,(i,i')} \equiv
\begin{pmatrix}
3& 2& 3& 0& 0& 0 \\
4& 0& 0& 3& 5& 0 \\
0& 2& 0& 2& 0& 1 \\
0& 0& 6& 0& 3& 2 \\
1& 4& 5& 3& 1& 3 \\
2& 6& 4& 2& 5& 4 
\end{pmatrix}
\pmod{p}, 
\]
and one can easily see that $\mathrm{rank}\,\overline{A^+} = 6$. 
Consequently, by the property \ref{main:lmdrk} of $A$, $X(\widetilde{k})$ is a nontrivial $\widetilde{\varLambda}$-cyclic $\widetilde{\varLambda}$-module. 
\end{example}

\subsection{Sextic fields with $t=1$.} 

\begin{example}\label{ex:3C6t1}
Let $k^+$ be the unique cubic subfield of $\mathbb Q(\zeta_{61})$. 
Put $p=3$ and $k=k^+(\sqrt{-23})$. 
Then $k=\mathbb Q(\alpha)$ with $\alpha$ satisfying 
$\alpha^6 - 2\alpha^5 + 30\alpha^4 - 34\alpha^3 + 2015\alpha^2 - 1740\alpha + 43551 =0$, 
and $p\mathcal{O}_k=\prod_{i=1}^{6}\mathfrak{p}_i$ with 
$\mathfrak{p}_1=3\mathcal{O}_k+(-\frac{469}{2\delta}\alpha^5 + \frac{971}{3\delta}\alpha^4 + \frac{9895}{3\delta}\alpha^3 + \frac{171851}{6\delta}\alpha^2 - \frac{981047}{6\delta}\alpha + \frac{221379}{2\delta})\mathcal{O}_k$,  
$\mathfrak{p}_i=g^{i-1}(\mathfrak{p}_1)$ for $i \in \{2,\ldots,6\}$, 
where $\delta=487 \cdot 499$, 
and $g$ is a generator of $\mathrm{Gal}(k/\mathbb Q)$ such that 
$g(\alpha)=-(\frac{388}{9\delta}\alpha^5 + \frac{5855}{9\delta}\alpha^4 + \frac{8360}{9\delta}\alpha^3 + \frac{24460}{9\delta}\alpha^2 + \frac{800795}{3\delta}\alpha + \frac{623540}{\delta})$. 
Since $Cl_k \simeq \mathbb Z/36\mathbb Z \times \mathbb Z/4\mathbb Z$, we have $t=1$ and $h=16$. 
Put
$\mathfrak{a}_0=41\mathcal{O}_k+(\alpha-5)\mathcal{O}_k$, 
which is a prime ideal of $k$ lying over $41$. 
The orders of $[\mathfrak{p}_1]$ and $[\mathfrak{a}_0]$ in $Cl_k$ are $q_0=36$. 
In particular, $Cl_k/p = \langle [\mathfrak{a}_0]Cl_k^p \rangle$. 
Then $\mathfrak{a}_0^{q_0}=\varepsilon_0\mathcal{O}_k$ with 
\begin{align*}
\varepsilon_0 = \ &
\textstyle{\frac{1888524068369}{18\delta}\alpha^5 + \frac{6943574603762}{9\delta}\alpha^4 - \frac{28025383797112}{9\delta}\alpha^3} \\
&\textstyle{- \frac{8978735955151}{18\delta}\alpha^2 + \frac{929227301965079}{6\delta}\alpha + \frac{4220669428821629}{2\delta}}
\end{align*}
The fundamental units of $k$ are written as 
\begin{align*}
\varepsilon_1 &= \textstyle{-\frac{354}{\delta}\alpha^5 - \frac{1584}{\delta}\alpha^4 - \frac{12638}{\delta}\alpha^3 + \frac{231}{\delta}\alpha^2 - \frac{176362}{\delta}\alpha + \frac{1138112}{\delta}} ,\\
\varepsilon_2 &= \textstyle{-\frac{320}{\delta}\alpha^5 + \frac{2687}{\delta}\alpha^4 - \frac{16916}{\delta}\alpha^3 + \frac{24922}{\delta}\alpha^2 - \frac{137456}{\delta}\alpha - \frac{860384}{\delta}} .
\end{align*}
Then 
\[
Z=\begin{pmatrix}
z_{0,1} & z_{0,2} & z_{0,3} & z_{0,4} & z_{0,5} & z_{0,6} \\
z_{1,1} & z_{1,2} & z_{1,3} & z_{1,4} & z_{1,5} & z_{1,6} \\
z_{2,1} & z_{2,2} & z_{2,3} & z_{2,4} & z_{2,5} & z_{2,6} 
\end{pmatrix} 
\equiv
\begin{pmatrix}
2 & 0 & 2 & 0 & 2 & 0 \\
2 & 0 & 1 & 2 & 0 & 1 \\
1 & 2 & 0 & 1 & 2 & 0
\end{pmatrix}
\pmod{p}, 
\]
and hence 
\[
\begin{pmatrix}
b_{4,0} & b_{4,1} & b_{4,2} \\ 
b_{5,0} & b_{5,1} & b_{5,2} \\ 
b_{6,0} & b_{6,1} & b_{6,2}  
\end{pmatrix}
=B^{-1}
\equiv
\begin{pmatrix}
2 & 0 & 1 \\
2 & 0 & 0 \\
2 & 1 & 1
\end{pmatrix}
\pmod{p}. 
\]
Since $\mathfrak{p}_i^h=\mathfrak{a}_0^{l_{i,0}} \varpi_i\mathcal{O}_k$ for $i \in \{1,\ldots,6\}$ 
with 
$(l_{1,0},l_{2,0},l_{3,0},l_{4,0},l_{5,0},l_{6,0})=(8, 4, 2, 1, 5, 7)$ and 
\begin{align*}
\varpi_1 =\ & \textstyle{\big( -\frac{459549066139891}{6\delta}\alpha^5 - \frac{704093863348502}{3\delta}\alpha^4 - \frac{535153539380573}{\delta}\alpha^3 
}\\&\textstyle{
- \frac{6243827791610063}{6\delta}\alpha^2 - \frac{808834974927696857}{6\delta}\alpha - \frac{1354543489312812447}{2\delta} \big)} 41^{-8} , \\
\varpi_2 =\ & \textstyle{\big( -\frac{1138617949}{18\delta}\alpha^5 + \frac{9263007134}{9\delta}\alpha^4 + \frac{21513813602}{9\delta}\alpha^3 
}\\&\textstyle{
+ \frac{904069241621}{18\delta}\alpha^2 + \frac{561744316739}{6\delta}\alpha + \frac{4365781564589}{2\delta} \big)} 41^{-4} , \\
\varpi_3 =\ & \textstyle{\big( \frac{14439515}{162\delta}\alpha^5 + \frac{32950277}{54\delta}\alpha^4 + \frac{90525712}{27\delta}\alpha^3 
+ \frac{572674304}{81\delta}\alpha^2 + \frac{1219360249}{27\delta}\alpha + \frac{2077748773}{18\delta} \big)} 41^{-2} , \\
\varpi_4 =\ & \textstyle{\big( \frac{1649801}{729\delta}\alpha^5 - \frac{181297}{54\delta}\alpha^4 - \frac{893647}{243\delta}\alpha^3 
+ \frac{287666645}{1458\delta}\alpha^2 + \frac{1291588133}{486\delta}\alpha - \frac{8536877}{9\delta} \big)} 41^{-1} , \\
\varpi_5 =\ & \textstyle{\big( -\frac{322820922181}{1458\delta}\alpha^5 + \frac{76556486800}{81\delta}\alpha^4 + \frac{5754643430519}{486\delta}\alpha^3 
}\\&\textstyle{
- \frac{20003383223033}{1458\delta}\alpha^2 - \frac{222480668144683}{243\delta}\alpha + \frac{183586414913389}{27\delta} \big)} 41^{-5} , \\
\varpi_6 =\ & \textstyle{\big( -\frac{10129957962787}{6\delta}\alpha^5 + \frac{99318329682100}{9\delta}\alpha^4 - \frac{84269452005608}{3\delta}\alpha^3 
}\\&\textstyle{
+ \frac{488192503813457}{6\delta}\alpha^2 - \frac{49749237476440823}{18\delta}\alpha + \frac{119523731360647853}{6\delta} \big)} 41^{-7} , 
\end{align*}
we have 
\[
(l_{i.j})_{1 \le i \le 6 \atop 1 \le j \le 6}
\equiv
\begin{pmatrix}
0 & 0 & 2 & 2 & 1 & 0 \\
1 & 0 & 2 & 0 & 0 & 0 \\
0 & 0 & 0 & 2 & 1 & 0 \\ 
0 & 0 & 0 & 0 & 0 & 2 \\
0 & 0 & 0 & 1 & 0 & 1 \\
0 & 0 & 0 & 0 & 1 & 0 
\end{pmatrix}
\pmod{p}, 
\]
and 
\[
A = (a_{j,(i,i')})_{j,(i,i')} \equiv
\begin{pmatrix}
1& 0& 0& 2& 2& 0 \\
0& 2& 0& 2& 0& 2 \\
0& 0& 1& 0& 2& 1 \\
2& 2& 1& 1& 0& 1 \\
2& 0& 2& 1& 0& 2 \\
1& 2& 2& 1& 2& 2 
\end{pmatrix}
\pmod{p}. 
\]
Since 
$\det A \equiv 2 \not\equiv 0 \pmod{p}$, we have $X(\widetilde{k}) \simeq 0$. 
\end{example}

\bigskip
\begin{acknowledgements}
This work was partially supported by JSPS KAKENHI Grant Number JP25K06965. 
\end{acknowledgements}

\bigskip


\bigskip

{\footnotesize 
\noindent\textsc{Yasushi Mizusawa}\\%
Department of Mathematics, 
Rikkyo University, 
3-34-1, Nishi-Ikebukuro, 
Toshima-ku, 
Tokyo 171-8501, Japan. 
\ \ 
\texttt{mizusawa@rikkyo.ac.jp}
\\
\\
\noindent\textsc{Keiji Okano}\\%
Department of Teacher Education, 
Tsuru University, 
3-8-1 Tahara, 
Tsuru-shi, 
Yamanashi 402-0054, Japan.
\ \ 
\texttt{okano@tsuru.ac.jp}
}


\newpage
\pagestyle{empty}
\setcounter{section}{0}
\renewcommand{\thesection}{\bf \Alph{section}}


\section{Symbolic computations}

\subsection{Quartic field with $t=0$.}

{\scriptsize \begin{verbatim}
test(case)={ \\ C4 if 1, V4 if 2
if(case!=0, z12=-z11);
z13=z11; z14=z12; Z=[z11,z12,z13,z14]; 

l=[0,l12,l13,l14; l21,0,l23,l24; l13,l14,0,l12; l23,l24,l21,0];
if(case==1,
l=[0,l12,l13,l14; l14,0,l12,l13; l13,l14,0,l12; l12,l13,l14,0]);
if(case==2,
l=[0,l12,l13,l14; l12,0,l14,l13; l13,l14,0,l12; l14,l13,l12,0]);

c4=vector(3);for(i=1,3, c4[i]=-(z14^-1)*Z[i]); 
L=matrix(4,3);for(j=1,4, for(i=1,3,
  L[j,i]=(l[j,i]-(z14^-1)*Z[i]*l[j,4])/h;
));
AA=[
 L[1,2], L[1,3],      0;
-L[2,1],      0, L[2,3];
      0,-L[3,1],-L[3,2];
      0,      0,      0;
      1,     -1,      1
];
AA[4,1]=c4[1]*L[4,2]-c4[2]*L[4,1];
AA[4,2]=c4[1]*L[4,3]-c4[3]*L[4,1];
AA[4,3]=c4[2]*L[4,3]-c4[3]*L[4,2];
A=AA[1..4,];

Q=matid(5);Q[2,2]=z12;Q[3,1]=h;Q[4,2]=h;Q[3,3]=h;Q[4,4]=h;
P=matid(3);P[1,1]=z11;P[2,1]=-z12;P[2,2]=z12;P[3,2]=-z11;
QAAP=[
    -z12*L[3,1],     z12*L[1,3],             0;
-z11*z12*L[2,1],-z11*z12*L[2,3],    z12*L[2,3];
              0,              0, l[1,2]-l[1,4];
              0,              0,-l[2,1]+l[2,3];
        z11+z12,       -z11-z12,             1];
print1(QAAP==Q*AA*P);
D=(l[1,2]*z11-l[1,3]*z12)*(l[2,4]*z11-l[2,3]*z12)+(l[1,4]*z11-l[1,3]*z12)*(l[2,4]*z11-l[2,1]*z12);
print1(matdet(QAAP[1..2,1..2])==z11*D/h^2);

if(case==0,
R=matid(3);R[3,1]=-1;
S=matid(4);S[3,1]=1;S[4,2]=1;
print1(S*subst(h*A,z11,0)*R==[
 l[1,2]-l[1,4],l[1,3],             0;
-l[2,1]-l[2,3],     0,        l[2,3];
             0,     0, l[1,2]-l[1,4];
             0,     0,-l[2,1]+l[2,3]]);
print(subst(D,z11,0)==Z[2]^2*l[1,3]*(l[2,1]+l[2,3]));
);
if(case==1, print1(D==z11^2*((l12+l13)^2+(l13+l14)^2)));
if(case==2, print1(D==2*z11^2*(l12+l13)*(l13+l14)));
if(case!=0,
t=0;for(i=1,3, t=t+(L[2,i]-L[1,i]));print1(t==0);
M2=[
L[1,2],L[2,2]-(L[1,1]+L[1,2]+L[1,3]);
L[1,3],L[2,3]];
print1(matdet(M2)*h^2*z11^2==D);
);
if(case==1, 
N=[
-L[1,2]-L[1,3], L[2,1];
 L[1,2],        L[2,2]-(L[1,1]+L[1,2]+L[1,3]);
 L[1,3],        L[2,3]];
a=h*(l12 + l13)/(D/z11^2);
b=-h*(l13 + l14)/(D/z11^2);
print(N*[a;b]==[-1;1;0]);
);
if(case==2, 
N=[
-L[1,2]-L[1,3]-L[2,1],L[1,2]-L[2,2]+(L[1,1]+L[1,2]+L[1,3]),L[1,3]-L[2,3];
-L[1,2]-L[1,3]+L[2,1],L[1,2]+L[2,2]-(L[1,1]+L[1,2]+L[1,3]),L[1,3]+L[2,3]];
a = h/(2*(l12 + l13));
b = h/(2*(l13 + l14));
print([a,0;0,b]*N==[-1,1,0;0,-1,1]);
);
}
\end{verbatim}}
\hrulefill

{\scriptsize \begin{verbatim}
gp > test()
1111
\end{verbatim}}
\hrulefill

{\scriptsize \begin{verbatim}
gp > test(1)
111111
\end{verbatim}}
\hrulefill

{\scriptsize \begin{verbatim}
gp > test(2)
111111
\end{verbatim}}
\hrulefill

\subsection{Quartic field with $t=1$.}

{\scriptsize \begin{verbatim}
test()={
Z=[ z01,z02,z03,z04; z11,z12,z13,z14 ];
B=[ z03,z04; z13,z14 ];
D=matdet(B);
b=B^-1;
l=[0,l12,l13,l14; l21,0,l23,l24; l31,l32,0,l34; l41,l42,l43,0];
c=matrix(4,2);for(j=3,4,for(i=1,2,
  c[j,i]=-b[j-2,1]*Z[1,i]-b[j-2,2]*Z[2,i]
));
L=matrix(4,2);for(j=1,4,for(i=1,2,
  L[j,i]=(l[j,i]+c[3,i]*l[j,3]+c[4,i]*l[j,4])/h
));
A=[    -l10/h,             0,                      L[1,2];
            0,        -l20/h,                     -L[2,1];
-c[3,1]*l30/h, -c[3,2]*l30/h, c[3,1]*L[3,2]-c[3,2]*L[3,1];
-c[4,1]*l40/h, -c[4,2]*l40/h, c[4,1]*L[4,2]-c[4,2]*L[4,1]];
DhA=[
matdet([ 0,l10,0,0; 1,0,0,0; z01,z02,z03,z04; z11,z12,z13,z14 ]),
0,
-matdet([ 0,l12,l13,l14; 1,0,0,0; z01,z02,z03,z04; z11,z12,z13,z14 ]);
0,
-matdet([ l20,0,0,0; 0,1,0,0; z01,z02,z03,z04; z11,z12,z13,z14 ]),
-matdet([ l21,0,l23,l24; 0,1,0,0; z01,z02,z03,z04; z11,z12,z13,z14 ]);
matdet([ 0,l30,0,0; 0,0,1,0; z01,z02,z03,z04; z11,z12,z13,z14 ]),
-matdet([ l30,0,0,0; 0,0,1,0; z01,z02,z03,z04; z11,z12,z13,z14 ]),
-matdet([ l31,l32,0,l34; 0,0,1,0; z01,z02,z03,z04; z11,z12,z13,z14 ]);
matdet([ 0,l40,0,0; 0,0,0,1; z01,z02,z03,z04; z11,z12,z13,z14 ]),
-matdet([ l40,0,0,0; 0,0,0,1; z01,z02,z03,z04; z11,z12,z13,z14 ]),
-matdet([ l41,l42,l43,0; 0,0,0,1; z01,z02,z03,z04; z11,z12,z13,z14 ])
];
print(DhA==D*h*A);
}
\end{verbatim}}
\hrulefill

{\scriptsize \begin{verbatim}
gp > test()
1
\end{verbatim}}
\hrulefill

\subsection{Sextic field with $t=0$.}

{\scriptsize \begin{verbatim}
test(case)={ \\ C6 if 1
if(case==1, 
  z13=-z11-z12; z23=-z21-z22; 
  l21=l16;l23=l12;l24=l13;l25=l14;l26=l15;
  l31=l15;l32=l16;l34=l12;l35=l13;l36=l14;
);
Z=[z11,z12,z13,z11,z12,z13; z21,z22,z23,z21,z22,z23];
l=[
  0,l12,l13,l14,l15,l16;
l21,  0,l23,l24,l25,l26;
l31,l32,  0,l34,l35,l36;
l14,l15,l16,  0,l12,l13;
l24,l25,l26,l21,  0,l23;
l34,l35,l36,l31,l32,  0];
B=[Z[1,5],Z[1,6];Z[2,5],Z[2,6]];
d1=matdet(B);
b=B^-1;
c=matrix(6,4);for(j=5,6,for(i=1,4,
  c[j,i]=-b[j-4,1]*Z[1,i]-b[j-4,2]*Z[2,i]
));
L=matrix(6,4);for(j=1,6,for(i=1,4,
  L[j,i]=(l[j,i]+c[5,i]*l[j,5]+c[6,i]*l[j,6])/h
));
A=[
 L[1,2], L[1,3], L[1,4],      0,      0,      0;
-L[2,1],      0,      0, L[2,3], L[2,4],      0;
      0,-L[3,1],      0,-L[3,2],      0, L[3,4];
      0,      0,-L[4,1],      0,-L[4,2],-L[4,3];
      0,      0,      0,      0,      0,      0;
      0,      0,      0,      0,      0,      0];
for(j=5,6,
A[j,1]=c[j,1]*L[j,2]-c[j,2]*L[j,1];
A[j,2]=c[j,1]*L[j,3]-c[j,3]*L[j,1];
A[j,3]=c[j,1]*L[j,4]-c[j,4]*L[j,1];
A[j,4]=c[j,2]*L[j,3]-c[j,3]*L[j,2];
A[j,5]=c[j,2]*L[j,4]-c[j,4]*L[j,2];
A[j,6]=c[j,3]*L[j,4]-c[j,4]*L[j,3];
);
d2=matdet([Z[1,3],Z[1,1];Z[2,3],Z[2,1]]);
d3=matdet([Z[1,1],Z[1,2];Z[2,1],Z[2,2]]);
A11=d1*[
 (l12-l15),  (l13-l16),          0;
-(l21-l24),          0,  (l23-l26);
         0, -(l31-l34), -(l32-l35)];
A12=[
l14*d1+l15*d2+l16*d3, 0, 0;
0, l24*d1+l25*d2+l26*d3, 0;
0, 0, l34*d1+l35*d2+l36*d3];
A22=[
-(l14*d1+l12*d2+l13*d3 + l14*d1+l15*d2+l16*d3), (l12-l15)*d1, (l13-l16)*d1;
(l21-l24)*d2, -(l21*d1+l25*d2+l23*d3 + l24*d1+l25*d2+l26*d3), (l23-l26)*d2;  
(l31-l34)*d3, (l32-l35)*d3, -(l31*d1+l32*d2+l36*d3 + l34*d1+l35*d2+l36*d3)];
Q=matid(6);Q[4,1]=-1;Q[5,2]=-1;Q[6,3]=-1; 
P=matid(6);P[3,3]=0;P[4,4]=0;P[3,4]=1;P[4,3]=1;P[5,1]=1;P[6,2]=1; 
print1(Q*d1*h*A*P==matconcat([A11,A12;matrix(3,3),A22]));
M=[
 matdet([l14,l15,l13;z11,z12,z13;z21,z22,z23]),
-matdet([l14,l15,l16;z11,z12,z13;z21,z22,z23]),
-matdet([l14,l12,l13;z11,z12,z13;z21,z22,z23]);
-matdet([l24,l25,l23;z11,z12,z13;z21,z22,z23]),
 matdet([l24,l25,l26;z11,z12,z13;z21,z22,z23]),
-matdet([l21,l25,l26;z11,z12,z13;z21,z22,z23]);
-matdet([l31,l32,l36;z11,z12,z13;z21,z22,z23]),
-matdet([l34,l35,l36;z11,z12,z13;z21,z22,z23]),
 matdet([l31,l35,l36;z11,z12,z13;z21,z22,z23])];
print1(matdet(A22)==2*matdet(M));

if(case==1,
print1(matdet(A11)==d1^3*(-(l12-l15)^3-(l13-l16)^3));
m1 = l12 +l14 +l16;
m2 = (l12+2*l13+2*l14+l15)*(l13+2*l14+2*l15+l16) + (l12+l13-l15-l16)^2;
print1(matdet(M)==-d1^3*m1*m2);
AA=matconcat([A;[1,-1,0,1,0,0;1,0,-1,0,1,0;0,1,-1,0,0,1]]);
S=matconcat([
h*matid(3), matrix(3),matrix(3);
h*matid(3),h*matid(3),matrix(3);
 matrix(3), matrix(3), matid(3)]);
R=[
0, 0, 1, 0, 0, 0;
1, 0, 0, 0, 0, 0;
1, 1, 1, 0, 0,-1;
1, 0,-1, 1, 0, 0;
1, 1, 0, 0, 1,-1;
0, 1, 1, 0, 0, 0];
MM=S*AA*R;
print1(MM[1..3,1..3]==[
 l13+l14+l15,l14+l15+l16, l12+l14+l16;
 l12+l13+l14,l13+l14+l15,-l12-l14-l16;
-l14-l15-l16,l12+l13+l14, l12+l14+l16]);
print1(MM[4..9,1..3]==matrix(6,3));
print1(MM[7..9,4..6]==matid(3));
print1(matdet(MM[1..3,1..3])==m1*m2);

t=0;for(i=1,4, t=t+(L[2,i]-L[1,i]));print1(t==0);
s=0;for(i=1,4, s=s+(L[3,i]-L[1,i]));print1(s==0);
L1=L[1,1]+L[1,2]+L[1,3]+L[1,4];
M3=[
L[1,2],L[2,2]-L1,L[3,2];
L[1,3],L[2,3]   ,L[3,3]-L1;
L[1,4],L[2,4]   ,L[3,4]];
print1(matdet(M3)*h^3==m1*m2);

aa = l12 + l13 + l14;
bb = -(l13 + l14 + l15);
cc = l14 + l15 + l16;
a =  h*( aa^2 - bb*cc )/(m1*m2);
b = -h*( cc^2 - aa*bb )/(m1*m2);
c =  h*( bb^2 - aa*cc )/(m1*m2);
N=mattranspose(L[1..3,]);
L1=L[1,1]+L[1,2]+L[1,3]+L[1,4];
N[1,1]=N[1,1]-L1;
N[2,2]=N[2,2]-L1;
N[3,3]=N[3,3]-L1;
print1(N*[a;b;c]==[-1;1;0;0]);
print1(
m1*m2==aa^3+bb^3+cc^3-3*aa*bb*cc,
a+b == h*(aa-cc)/m2,
aa-cc == l12+l13-l15-l16,
(a+b)^2+3*c^2 == h^3*((l12 - l15)^3 +(l13 - l16)^3)/((a+b)*m1^2*m2^2) );
);
}
\end{verbatim}}
\hrulefill

{\scriptsize \begin{verbatim}
gp > test()
11
\end{verbatim}}
\hrulefill

{\scriptsize \begin{verbatim}
gp > test(1)
1111111111111111
\end{verbatim}}
\hrulefill

\subsection{Sextic field with $t=1$.}

{\scriptsize \begin{verbatim}
test()={
Z=[
z01,z02,z03,z04,z05,z06;
z11,z12,z13,z14,z15,z16;
z21,z22,z23,z24,z25,z26];
l=[
  0,l12,l13,l14,l15,l16;
l21,  0,l23,l24,l25,l26;
l31,l32,  0,l34,l35,l36;
l41,l42,l43,  0,l45,l46;
l51,l52,l53,l54,  0,l56;
l61,l62,l63,l64,l65,  0];
B=Z[,4..6];
D=matdet(B);
b=matconcat([matrix(3,3);B^-1]);
c=matrix(6,3);for(j=4,6, for(i=1,3, 
  c[j,i]=-b[j,1]*Z[1,i]-b[j,2]*Z[2,i]-b[j,3]*Z[3,i]
));
L=matrix(6,3);for(j=1,6, for(i=1,3,
  L[j,i]=(l[j,i]+c[4,i]*l[j,4]+c[5,i]*l[j,5]+c[6,i]*l[j,6])/h
));
l0=[l10,l20,l30,l40,l50,l60];
A=[
-l0[1]/h,       0,       0, L[1,2], L[1,3],      0;
       0,-l0[2]/h,       0,-L[2,1],      0, L[2,3];
       0,       0,-l0[3]/h,      0,-L[3,1],-L[3,2];
0,0,0,0,0,0;0,0,0,0,0,0;0,0,0,0,0,0];
for(i=4,6,
A[i,1]=-c[i,1]*l0[i]/h;
A[i,2]=-c[i,2]*l0[i]/h;
A[i,3]=-c[i,3]*l0[i]/h;
A[i,4]=c[i,1]*L[i,2]-c[i,2]*L[i,1];
A[i,5]=c[i,1]*L[i,3]-c[i,3]*L[i,1];
A[i,6]=c[i,2]*L[i,3]-c[i,3]*L[i,2];
);

d123=matdet(Mat([Z[,1],Z[,2],Z[,3]]));
d124=matdet(Mat([Z[,1],Z[,2],Z[,4]]));
d125=matdet(Mat([Z[,1],Z[,2],Z[,5]]));
d126=matdet(Mat([Z[,1],Z[,2],Z[,6]]));
d134=matdet(Mat([Z[,1],Z[,3],Z[,4]]));
d135=matdet(Mat([Z[,1],Z[,3],Z[,5]]));
d136=matdet(Mat([Z[,1],Z[,3],Z[,6]]));
d145=matdet(Mat([Z[,1],Z[,4],Z[,5]]));
d146=matdet(Mat([Z[,1],Z[,4],Z[,6]]));
d156=matdet(Mat([Z[,1],Z[,5],Z[,6]]));

d234=matdet(Mat([Z[,2],Z[,3],Z[,4]]));
d235=matdet(Mat([Z[,2],Z[,3],Z[,5]]));
d236=matdet(Mat([Z[,2],Z[,3],Z[,6]]));
d245=matdet(Mat([Z[,2],Z[,4],Z[,5]]));
d246=matdet(Mat([Z[,2],Z[,4],Z[,6]]));
d256=matdet(Mat([Z[,2],Z[,5],Z[,6]]));

d345=matdet(Mat([Z[,3],Z[,4],Z[,5]]));
d346=matdet(Mat([Z[,3],Z[,4],Z[,6]]));
d356=matdet(Mat([Z[,3],Z[,5],Z[,6]]));

d456=matdet(Mat([Z[,4],Z[,5],Z[,6]]));

A11=d456*[-l0[1],0,0;0,-l0[2],0;0,0,-l0[3]];
A12=[
l12*d456-l14*d256+l15*d246-l16*d245,
l13*d456-l14*d356+l15*d346-l16*d345,
0;
-l21*d456+l24*d156-l25*d146+l26*d145,
0,
l23*d456-l24*d356+l25*d346-l26*d345;
0,
-l31*d456+l34*d156-l35*d146+l36*d145,
-l32*d456+l34*d256-l35*d246+l36*d245
];
A21=[l0[4],0,0;0,l0[5],0;0,0,l0[6]]*[
 d156, d256, d356;
-d146,-d246,-d346;
 d145, d245, d345];
A22=[
 l41*d256-l42*d156+l45*d126-l46*d125,
 l41*d356-l43*d156+l45*d136-l46*d135,
 l42*d356-l43*d256+l45*d236-l46*d235;
-l51*d246+l52*d146-l54*d126+l56*d124,
-l51*d346+l53*d146-l54*d136+l56*d134,
-l52*d346+l53*d246-l54*d236+l56*d234;
 l61*d245-l62*d145+l64*d125-l65*d124,
 l61*d345-l63*d145+l64*d135-l65*d134,
 l62*d345-l63*d245+l64*d235-l65*d234];
print(D*h*A==matconcat([A11,A12;A21,A22]));
}
\end{verbatim}}
\hrulefill

{\scriptsize \begin{verbatim}
gp > test()
1
\end{verbatim}}
\hrulefill

\newpage

\section{Numerical computations}

\subsection{Examples \ref{ex:3V4t0} and \ref{ex:3V4not0}.}

{\scriptsize \begin{verbatim}
test(m)={\\ 11 or 17 or 65
p=3;n=7;
print([p,n,-m]);
f=x^4+2*(m-n)*x^2+(m+n)^2; \\ x = a + b, a = sqrt(n), b = sqrt(-m)
print(f," f");
k=bnfinit(f);
print(bnfcertify(k), " bnfcertify");
print(k.clgp[2], " k.clgp");
R=bnrinit(k,p^2);
print(matreduce(lift(R.clgp[2]*Mod(1,p)))[1,2], " d(G~)");

kill(a);kill(b);
Zk=listcreate();
for(i=1,4,
listput(Zk,lift(Mod(lift(Mod(subst(k.zk[i],x,a+b),a^2-n)),b^2+m)))
);
print(Vec(Zk), " k.zk");

a=Mod((x^2+n+m)/(2*x),f);print(lift(a^2));
b=Mod((x^2-n-m)/(2*x),f);print(lift(b^2));

if(m==11,ff=x+(1+y)/2,ff=(x+y)/2+1);
P=listcreate();
listput(P,idealismaximal(k,idealhnf(k, p, subst(subst(ff,x, a),y, b) )));
listput(P,idealismaximal(k,idealhnf(k, p, subst(subst(ff,x,-a),y, b) )));
listput(P,idealismaximal(k,idealhnf(k, p, subst(subst(ff,x, a),y,-b) )));
listput(P,idealismaximal(k,idealhnf(k, p, subst(subst(ff,x,-a),y,-b) )));
kill(a);kill(b);
for(i=1,4,print(P[i][1..4]," P",i));

pZk=idealmul(k, idealmul(k, idealmul(k, P[1], P[2]), P[3]), P[4]);
print(bnfisprincipal(k,pZk), " P1P2P3P4");

P1P3=bnfisprincipal(k,idealmul(k, P[1], P[3]));
print(P1P3, " P1P3");
print(lift(Mod(lift(Mod(subst(k.zk*P1P3[2],x,a+b),a^2-n)),b^2+m)), " gen P1P3");

P1P2=bnfisprincipal(k,idealmul(k, P[1], P[2]));
print(P1P2, " P1P2");
print(lift(Mod(lift(Mod(subst(k.zk*P1P2[2],x,a+b),a^2-n)),b^2+m)), " gen P1P2");

e=k.fu[1];
print(lift(Mod(lift(Mod(subst(lift(e),x,a+b),a^2-n)),b^2+m)), " e1");

z1=listcreate();
for(i=1,length(P),
for(z=0,p-1,
if(idealval(k,e^(p-1)-(1+p)^((p-1)*z),P[i]) > 1, listput(z1,z);break)
)
);
print(Vec(z1)," Z");

h=k.clgp[1];
w=listcreate();
for(i=1,4,listput(w,
Mod(k.zk*bnfisprincipal(k,idealpow(k,P[i],h))[2],f)
));
print(lift(w[1]), " w1");

W=listcreate();
for(i=1,4,
listput(W,lift(Mod(lift(Mod(subst(lift(w[i]),x,a+b),a^2-n)),b^2+m)))
);
for(i=1,4, print(W[i], " W",i));

l=matrix(4,4);
for(i=1,length(P),for(j=1,length(P), if(j==i,l[i,j]=0,
for(z=0,p-1,
if(idealval(k, w[i]^(p-1)-(1+p)^((p-1)*z),P[j]) > 1, l[i,j]=z;break)
)
)));
print(l, " l_{i,j} mod p");

print((l[1,2]-l[1,4])%p != 0 || (l[2,1]-l[2,3])%p != 0);
print(((l[1,2]+l[1,3])*(l[1,3]+l[1,4]))%p != 0);

D = (l[1,2]*z1[1]-l[1,3]*z1[2])*(l[2,4]*z1[1]-l[2,3]*z1[2])+(l[1,4]*z1[1]-l[1,3]*z1[2])*(l[2,4]*z1[1]-l[2,1]*z1[2]);
print(D%p != 0);

c4=vector(3); for(i=1,3, c4[i]=-(z1[4]^-1)*z1[i]); 
L=matrix(4,3); 
for(j=1,4, for(i=1,3,
  L[j,i]=(l[j,i]-(z1[4]^-1)*z1[i]*l[j,4])/h;
));
A=[
 L[1,2], L[1,3],      0;
-L[2,1],      0, L[2,3];
      0,-L[3,1],-L[3,2];
      0,      0,      0
];
A[4,1]=c4[1]*L[4,2]-c4[2]*L[4,1];
A[4,2]=c4[1]*L[4,3]-c4[3]*L[4,1];
A[4,3]=c4[2]*L[4,3]-c4[3]*L[4,2];

A=A*Mod(1,p);
print(matrank(A));
A=lift(A);
print(A, " A");

if(((l[1,2]+l[1,3])*(l[1,3]+l[1,4]))%p != 0, 
a=Mod(h/(2*(l[1,2]+l[1,3])),p); b=Mod(h/(2*(l[1,3]+l[1,4])),p); 
print([lift(a),lift(b)]," [a,b] mod p"); 
);
}
\end{verbatim}}
\hrulefill

{\scriptsize \begin{verbatim}
gp > test(11)
[3, 7, -11]
x^4 + 8*x^2 + 324 f
1 bnfcertify
[4] k.clgp
3 d(G~)
[1, (1/2*b - 1/2), -a, (-1/2*b - 1/2)*a] k.zk
7
-11
[3, [-1, 1, 0, 1]~, 1, 1] P1
[3, [-1, -1, -1, 0]~, 1, 1] P2
[3, [0, -1, -1, 0]~, 1, 1] P3
[3, [-2, -1, -1, 1]~, 1, 1] P4
[[0]~, [3, 0, 0, 0]~] P1P2P3P4
[[0]~, [-2, 0, -1, 0]~] P1P3
a - 2 gen P1P3
[[0]~, [0, 1, 0, 0]~] P1P2
(1/2*b - 1/2) gen P1P2
3*a + 8 e1
[1, 2, 1, 2] Z
1/72*x^3 - 23/36*x + 5/2 w1
-a + (-1/2*b + 5/2) W1
-a + (1/2*b - 5/2) W2
-a + (1/2*b + 5/2) W3
-a + (-1/2*b - 5/2) W4
[0, 0, 1, 1; 0, 0, 1, 1; 1, 1, 0, 0; 1, 1, 0, 0] l_{i,j} mod p
1
1
1
3
[2, 2, 0; 2, 0, 2; 0, 2, 2; 2, 2, 2] A
[2, 1] [a,b] mod p
\end{verbatim}}
\hrulefill

{\scriptsize \begin{verbatim}
gp > test(17)
[3, 7, -17]
x^4 + 20*x^2 + 576 f
1 bnfcertify
[20] k.clgp
3 d(G~)
[1, 1/2*a + 1/2*b, -1/2*a + 1/2*b, 1/2*b*a - 1/2] k.zk
7
-17
[3, [-1, -1, 0, 0]~, 1, 1] P1
[3, [1, 0, -1, -1]~, 1, 1] P2
[3, [1, 0, 1, -1]~, 1, 1] P3
[3, [1, -1, 0, 0]~, 1, 1] P4
[[0]~, [3, 0, 0, 0]~] P1P2P3P4
[[0]~, [-2, 1, -1, 0]~] P1P3
a - 2 gen P1P3
[[10]~, [901/157464, 2125/629856, 6859/944784, 3623/944784]~] P1P2
(3623/1889568*b - 7343/3779136)*a + (20093/3779136*b + 7189/1889568) gen P1P2
3*a + 8 e1
[1, 2, 1, 2] Z
1/8*x^3 + 39/4*x^2 - 19*x + 11 w1
(39/2*b - 49/2)*a + (-37/2*b - 173/2) W1
(-39/2*b + 49/2)*a + (-37/2*b - 173/2) W2
(-39/2*b - 49/2)*a + (37/2*b - 173/2) W3
(39/2*b + 49/2)*a + (37/2*b - 173/2) W4
[0, 0, 1, 0; 0, 0, 0, 1; 1, 0, 0, 0; 0, 1, 0, 0] l_{i,j} mod p
0
1
1
2
[0, 2, 0; 1, 0, 2; 0, 1, 0; 2, 0, 1] A
[1, 1] [a,b] mod p
\end{verbatim}}
\hrulefill

{\scriptsize \begin{verbatim}
gp > test(65)
[3, 7, -65]
x^4 + 116*x^2 + 5184 f
1 bnfcertify
[40, 2] k.clgp
3 d(G~)
[1, -a, 1/2*a - 1/2*b, -1/2*b*a + 1/2] k.zk
7
-65
[3, [-1, 1, 1, 0]~, 1, 1] P1
[3, [-2, 0, 1, 1]~, 1, 1] P2
[3, [-2, 0, -1, 1]~, 1, 1] P3
[3, [1, 1, 1, 0]~, 1, 1] P4
[[0, 0]~, [3, 0, 0, 0]~] P1P2P3P4
[[0, 0]~, [2, 1, 0, 0]~] P1P3
-a + 2 gen P1P3
[[30, 0]~, [16514123429/2248880205492486144, -3814055107/2248880205492486144, -6
669591541/2248880205492486144, 1732545589/2248880205492486144]~] P1P2
(-1732545589/4497760410984972288*b + 958518673/4497760410984972288)*a + (6669591
541/4497760410984972288*b + 34760792447/4497760410984972288) gen P1P2
3*a + 8 e1
[1, 2, 1, 2] Z
6664171/2*x^3 - 62398544*x^2 + 688638986*x - 4014380975 w1
(-124797088*b + 62206912)*a + (542027224*b - 395265423) W1
(124797088*b - 62206912)*a + (542027224*b - 395265423) W2
(-124797088*b - 62206912)*a + (542027224*b + 395265423) W3
(124797088*b + 62206912)*a + (542027224*b + 395265423) W4
[0, 1, 1, 2; 1, 0, 2, 1; 1, 2, 0, 1; 2, 1, 1, 0] l_{i,j} mod p
1
0
0
2
[1, 0, 0; 2, 0, 0; 0, 2, 1; 0, 1, 2] A
\end{verbatim}}
\hrulefill

\subsection{Example \ref{ex:3C4t0}.}

{\scriptsize \begin{verbatim}
pre()={
k=bnfinit(y^2-13);
print(bnfcertify(k), " bnfcertify");
bnr=bnrinit(k,[13,[1,1]]);
kummer=rnfkummer(bnr,2);print(kummer);
print(subst(X^2 + 130*X + 1300,X,-subst(kummer,x,0)*10));
}
\end{verbatim}}
\hrulefill

{\scriptsize \begin{verbatim}
gp > pre()
1 bnfcertify
x^2 + Mod(3/2*y + 13/2, y^2 - 13)
Mod(0, y^2 - 13)
\end{verbatim}}
\hrulefill

{\scriptsize \begin{verbatim}
test()={
p=3;
f=(x^2+5*13)^2-13*15^2; \\ x^4 + 130*x^2 + 1300
print(f," f");
k=bnfinit(f);
print(bnfcertify(k), " bnfcertify");
print(k.clgp[2], " k.clgp");
R=bnrinit(k,p^2);
print(matreduce(lift(R.clgp[2]*Mod(1,p)))[1,2], " d(G~)");

a=Mod((x^2+65)/15,f); print(lift(a^2)); \\ a = sqrt(13)
print(Mod(x^2,f) == -65+15*a);
e=k.fu[1]; 
print(lift(e), " e1");
print(e == (3+a)/2);
print(lift(e) == 1/30*x^2 + 11/3 );

g=nfgaloisconj(k)[4]; print(g, " gen Gal(k/Q)");
print(subst(lift(Mod(g*x,f)),x,sqrt(65-15*sqrt(13))*I), " < 0");
print(g == lift(x*e));
print((x*e)^2 == -65-15*a);

Pgen = 1/15*x^3 + 1/30*x^2 + 16/3*x + 11/3;
print(e+2*(e-1)*x == Mod(Pgen,f));
tmp=x;
P=listcreate();
for(i=1,4,
listput(P,idealismaximal(k,idealhnf(k, p, Pgen )));
Pgen=lift(Mod(subst(Pgen,x,lift(x*e)),f));
tmp=lift(Mod(subst(tmp,x,lift(x*e)),f));
);
print(tmp==x); \\ g^4 = id
for(i=1,4,print(P[i][1..4], " P",i));
pZk=idealmul(k, idealmul(k, idealmul(k, P[1], P[2]), P[3]), P[4]);
print(bnfisprincipal(k,pZk), " P1P2P3P4");

z1=listcreate();
for(i=1,length(P),
for(z=0,p-1,
if(idealval(k,e^(p-1)-(1+p)^((p-1)*z),P[i]) > 1, listput(z1,z);break)
)
);
print(Vec(z1), " Z");

h=k.clgp[1];
print(h, " h");
w=listcreate();
w0=-Mod(k.zk*bnfisprincipal(k,idealpow(k,P[1],h))[2],f)/e;
for(i=1,4,
listput(w,w0);
w0=Mod(subst(lift(w0),x,lift(x*e)),f);
);
print(w0==w[1]);
for(i=1,4,
print(lift(w[i]), " w",i)
);

l=matrix(4,4);
for(i=1,length(P),for(j=1,length(P), if(j==i,l[i,j]=0,
for(z=0,p-1,
if(idealval(k, w[i]^(p-1)-(1+p)^((p-1)*z),P[j]) > 1, l[i,j]=z;break)
)
)));

print(l, " l_{i,j} mod p");

print((l[1,2]-l[1,4])%p != 0 || (l[2,1]-l[2,3])%p != 0);
print(((l[1,2]+l[1,3])^2+(l[1,3]+l[1,4])^2)%p != 0);

D = (l[1,2]*z1[1]-l[1,3]*z1[2])*(l[2,4]*z1[1]-l[2,3]*z1[2])+(l[1,4]*z1[1]-l[1,3]*z1[2])*(l[2,4]*z1[1]-l[2,1]*z1[2]);
print(D%p != 0);

c4=vector(3); for(i=1,3, c4[i]=-(z1[4]^-1)*z1[i]); 
L=matrix(4,3); 
for(j=1,4, for(i=1,3,
  L[j,i]=(l[j,i]-(z1[4]^-1)*z1[i]*l[j,4])/h;
));
A=[
 L[1,2], L[1,3],      0;
-L[2,1],      0, L[2,3];
      0,-L[3,1],-L[3,2];
      0,      0,      0
];
A[4,1]=c4[1]*L[4,2]-c4[2]*L[4,1];
A[4,2]=c4[1]*L[4,3]-c4[3]*L[4,1];
A[4,3]=c4[2]*L[4,3]-c4[3]*L[4,2];

A=A*Mod(1,p);
print(matrank(A));
A=lift(A);
print(A, " A");

if(((l[1,2]+l[1,3])^2+(l[1,3]+l[1,4])^2)%p != 0, 
a=Mod( h*(l[1,2]+l[1,3])/((l[1,2]+l[1,3])^2+(l[1,3]+l[1,4])^2),p);
b=Mod(-h*(l[1,3]+l[1,4])/((l[1,2]+l[1,3])^2+(l[1,3]+l[1,4])^2),p);
print([lift(a),lift(b)]," [a,b] mod p"); 
);
}
\end{verbatim}}
\hrulefill

{\scriptsize \begin{verbatim}
gp > test()
x^4 + 130*x^2 + 1300 f
1 bnfcertify
[34, 2] k.clgp
3 d(G~)
13
1
1/30*x^2 + 11/3 e1
1
1
1/30*x^3 + 11/3*x gen Gal(k/Q)
-36.0555127546398929311922126747049594625129657384524621271045305622716694829301
04452046190820184907176735141820240635403760306782646978077051630171668927097577
426905642741526332338303949623469447962732 < 0
1
1
1
1
[3, [1, 1, -1, -1]~, 1, 1] P1
[3, [-1, -1, -1, 1]~, 1, 1] P2
[3, [1, 1, 1, 1]~, 1, 1] P3
[3, [-1, -1, 1, -1]~, 1, 1] P4
[[0, 0]~, [3, 0, 0, 0]~] P1P2P3P4
[1, 2, 1, 2] Z
68 h
1
1194376/3*x^3 - 4274169/2*x^2 + 30168356/3*x + 5545534 w1
-16868834/15*x^3 + 4274169/2*x^2 - 406945628/3*x + 283366519 w2
-1194376/3*x^3 - 4274169/2*x^2 - 30168356/3*x + 5545534 w3
16868834/15*x^3 + 4274169/2*x^2 + 406945628/3*x + 283366519 w4
[0, 2, 2, 1; 1, 0, 2, 2; 2, 1, 0, 2; 2, 2, 1, 0] l_{i,j} mod p
1
1
1
3
[2, 0, 0; 0, 0, 2; 0, 1, 2; 2, 1, 0] A
[2, 0] [a,b] mod p
\end{verbatim}}
\hrulefill

\subsection{Example \ref{ex:3V4t1}.}

{\scriptsize \begin{verbatim}
test()={
p=3;n=7;m=23;Na=29;
print1([p,n,-m,Na]);
print([subcycloiwasawa(x^2+m,p)[1],subcycloiwasawa(x^2+m*n,p)[1]]);
f=x^4+2*(m-n)*x^2+(m+n)^2; \\ x = a + b, a = sqrt(n), b = sqrt(-m)
print(f," f");
k=bnfinit(f);
print(bnfcertify(k), " bnfcertify");
print(k.clgp[2], " k.clgp");
R=bnrinit(k,p^2);
print(matreduce(lift(R.clgp[2]*Mod(1,p)))[1,2], " d(G~)");

a=Mod((x^2+m+n)/(2*x),f);print(lift(a^2));
b=Mod((x^2-m-n)/(2*x),f);print(lift(b^2));

P=listcreate();
listput(P,idealismaximal(k,idealhnf(k, p,  a + (1+b)/2 )));
listput(P,idealismaximal(k,idealhnf(k, p, -a + (1+b)/2 )));
listput(P,idealismaximal(k,idealhnf(k, p,  a + (1-b)/2 )));
listput(P,idealismaximal(k,idealhnf(k, p, -a + (1-b)/2 )));
kill(a);kill(b);

for(i=1,4,
print1(idealprimedec(k,p)[i]==P[1])
);print();

pZk=idealmul(k, idealmul(k, idealmul(k, P[1], P[2]), P[3]), P[4]);
print(bnfisprincipal(k,pZk), " P1P2P3P4");

h=k.clgp[1]/(p^valuation(k.clgp[1],p));print(h, " h");

A0=idealprimedec(k,Na)[1];
S=bnfisprincipal(k,A0)[1];print(S, " A0 in k.clgp");
print(k.zk*A0[2], " gen A0");

q0=1;
for(i=1,length(k.clgp[2]),
q0=lcm(k.clgp[2][i]/gcd(k.clgp[2][i],S[i]),q0)
);
e0=k.zk*bnfisprincipal(k,idealpow(k,A0,q0))[2];
print(q0, " q0");
print(lift(Mod(lift(Mod(subst(e0,x,a+b),a^2-n)),b^2+m)), " e0");

Z=matrix(2,4);

for(i=1,length(P),
for(z=0,p-1,
if(idealval(k, e0^(p-1)-(1+p)^((p-1)*z) , P[i]) > 1,
Z[1,i]=z; break
)));

e=k.fu[1];
for(i=1,length(P),
for(z=0,p-1,
if(idealval(k, e^(p-1)-(1+p)^((p-1)*z) , P[i]) > 1,
Z[2,i]=z; break
)));

print(lift(Mod(lift(Mod(subst(lift(e),x,a+b),a^2-n)),b^2+m)), " e1");

print(Z, " Z");

A0inv=idealinv(k,A0);

l0=listcreate();
w=listcreate();
for(i=1,length(P),Ph=idealpow(k,P[i],h);
for(z=1,k.clgp[1],
w0=bnfisprincipal(k,idealmul(k,Ph,idealpow(k,A0inv,z)));
if(w0[1][1]==0, 
listput(l0,z);listput(w,k.zk*w0[2]);break
)));
print(Vec(l0), " l_{i,0}");

W=listcreate();
for(i=1,4,
listput(W,lift(Mod(lift(Mod(subst(lift(w[i]),x,a+b),a^2-n)),b^2+m)))
);
for(i=1,4, print(W[i]*Na^l0[i], " W",i,"*",Na,"^",l0[i]));

l=matrix(4,4);
for(i=1,length(P),for(j=1,length(P), if(j==i,l[i,j]=0,
for(z=0,p-1,
if(idealval(k, w[i]^(p-1)-(1+p)^((p-1)*z),P[j]) > 1, l[i,j]=z;break)
)
)));
print(l, " l_{i.j} mod p");

l=l*Mod(1,p);
Z=Z*Mod(1,p);

B=Z[,3..4];
b=B^-1; print(lift(b), " B^-1");

c=matrix(4,2);
for(j=3,4,for(i=1,2,
  c[j,i]=-b[j-2,1]*Z[1,i]-b[j-2,2]*Z[2,i]
));

L=matrix(4,2);
for(j=1,4,for(i=1,2,
  L[j,i]=(l[j,i]+c[3,i]*l[j,3]+c[4,i]*l[j,4])/h
));

A=[
-l0[1]/h, 0,  L[1,2];
0, -l0[2]/h, -L[2,1];
-c[3,1]*l0[3]/h, -c[3,2]*l0[3]/h, c[3,1]*L[3,2]-c[3,2]*L[3,1];
-c[4,1]*l0[4]/h, -c[4,2]*l0[4]/h, c[4,1]*L[4,2]-c[4,2]*L[4,1]
];

kill(b);

A=A*Mod(1,p);
print(matrank(A));
A=lift(A);
print(A, " A");
}
\end{verbatim}}
\hrulefill

{\scriptsize \begin{verbatim}
gp > test()
[3, 7, -23, 29][1, 1]
x^4 + 32*x^2 + 900 f
1 bnfcertify
[24] k.clgp
3 d(G~)
7
-23
0010
[[0]~, [3, 0, 0, 0]~] P1P2P3P4
8 h
[2]~ A0 in k.clgp
x - 15 gen A0
12 q0
(891/2*b - 34839/2)*a + (287/2*b + 88779/2) e0
-3*a - 8 e1
[2, 2, 1, 1; 1, 2, 1, 2] Z
[4, 4, 8, 8] l_{i,0}
(-10091*b + 82468)*a + (-28948*b + 103404) W1*29^4
(5996*b - 77809)*a + (13556*b - 32004) W2*29^4
(1135888013/2*b - 3178361441/2)*a + (2975830185/2*b - 4780928131/2) W3*29^8
(-767151197/2*b + 2096910929/2)*a + (-2164820295/2*b + 99743021/2) W4*29^8
[0, 2, 0, 1; 1, 0, 2, 2; 2, 0, 0, 2; 1, 1, 1, 0] l_{i.j} mod p
[2, 2; 2, 1] B^-1
3
[1, 0, 1; 0, 1, 0; 0, 2, 1; 2, 0, 1] A
\end{verbatim}}
\hrulefill

\subsection{Example \ref{ex:5C4t1}.}

{\scriptsize \begin{verbatim}
pre()={
m=29;
k=bnfinit(y^2-m);
print(bnfcertify(k), " bnfcertify");
bnr=bnrinit(k,[m,[1,1]]);
kummer=rnfkummer(bnr,2);print(kummer);
print(subst(X^2 + 87*X + 261,X,-subst(kummer,x,0)*3));
}
\end{verbatim}}
\hrulefill

{\scriptsize \begin{verbatim}
gp > pre()
1 bnfcertify
x^2 + Mod(5/2*y + 29/2, y^2 - 29)
Mod(0, y^2 - 29)
\end{verbatim}}
\hrulefill

{\scriptsize \begin{verbatim}
test()={
p=5;m=29;n=3;Na=23;
f=(9/4)*( ((2/3)*x^2+29)^2-29*5^2 ); \\ x^4 + 87*x^2 + 261
print(f," f");
k=bnfinit(f);
print(bnfcertify(k), " bnfcertify");
print(k.clgp[2], " k.clgp");
R=bnrinit(k,p^2);
print(matreduce(lift(R.clgp[2]*Mod(1,p)))[1,2], " d(G~)");

a=Mod((2*x^2+3*29)/15,f); print(lift(a^2)); \\ a = sqrt(29)
print(Mod(x^2,f) == ((-29+5*a)/2)*3);
e=k.fu[1]; 
print(lift(e), " e1");
print(e == (5-a)/2);

g=nfgaloisconj(k)[3]; print(g, " gen Gal(k/Q)");
print(subst(lift(Mod(g*x,f)),x,sqrt((29-5*sqrt(29))*3/2)*I), " > 0");
print(Mod(g^2,f) == ((-29-5*a)/2)*3);

P=idealprimedec(k,p);
Pgen = k.zk*P[1][2]; print(Pgen, " gen P1");
tmp=x;
P=listcreate();
for(i=1,4,
listput(P,idealismaximal(k,idealhnf(k, p, Pgen )));
Pgen=lift(Mod(subst(Pgen,x,g),f));
tmp=lift(Mod(subst(tmp,x,g),f));
);
print(tmp==x); \\ g^4 = id
for(i=1,4,print(P[i][1..4], " P",i));
pZk=idealmul(k, idealmul(k, idealmul(k, P[1], P[2]), P[3]), P[4]);
print(bnfisprincipal(k,pZk), " P1P2P3P4");

h=k.clgp[1]/(p^valuation(k.clgp[1],p));print(h, " h");

A0=idealprimedec(k,Na)[1];
S=bnfisprincipal(k,A0)[1];print(S, " A0 in k.clgp");
print(k.zk*A0[2], " gen A0");

q0=1;
for(i=1,length(k.clgp[2]),
q0=lcm(k.clgp[2][i]/gcd(k.clgp[2][i],S[i]),q0)
);
e0=k.zk*bnfisprincipal(k,idealpow(k,A0,q0))[2];
print(q0, " q0");
print(e0, " e0");

Z=matrix(2,4);

for(i=1,length(P),
for(z=0,p-1,
if(idealval(k, e0^(p-1)-(1+p)^((p-1)*z) , P[i]) > 1,
Z[1,i]=z; break
)));

for(i=1,length(P),
for(z=0,p-1,
if(idealval(k, e^(p-1)-(1+p)^((p-1)*z) , P[i]) > 1,
Z[2,i]=z; break
)));

print(Z, " Z");

A0inv=idealinv(k,A0);

l0=listcreate();
w=listcreate();
for(i=1,length(P),Ph=idealpow(k,P[i],h);
for(z=1,k.clgp[1],
w0=bnfisprincipal(k,idealmul(k,Ph,idealpow(k,A0inv,z)));
if(w0[1][1]==0, 
listput(l0,z);listput(w,k.zk*w0[2]);break
)));
print(Vec(l0), " l_{i,0}");

for(i=1,4, print(w[i]*Na^l0[i], " w",i,"*",Na,"^",l0[i]));

l=matrix(4,4);
for(i=1,length(P),for(j=1,length(P), if(j==i,l[i,j]=0,
for(z=0,p-1,
if(idealval(k, w[i]^(p-1)-(1+p)^((p-1)*z),P[j]) > 1, l[i,j]=z;break)
)
)));

print(l, " l_{i,j} mod p");

l=l*Mod(1,p);
Z=Z*Mod(1,p);

B=Z[,3..4];
b=B^-1; print(lift(b), " B^-1");

c=matrix(4,2);
for(j=3,4,for(i=1,2,
  c[j,i]=-b[j-2,1]*Z[1,i]-b[j-2,2]*Z[2,i]
));

L=matrix(4,2);
for(j=1,4,for(i=1,2,
  L[j,i]=(l[j,i]+c[3,i]*l[j,3]+c[4,i]*l[j,4])/h
));

A=[
-l0[1]/h, 0,  L[1,2];
0, -l0[2]/h, -L[2,1];
-c[3,1]*l0[3]/h, -c[3,2]*l0[3]/h, c[3,1]*L[3,2]-c[3,2]*L[3,1];
-c[4,1]*l0[4]/h, -c[4,2]*l0[4]/h, c[4,1]*L[4,2]-c[4,2]*L[4,1]
];

A=A*Mod(1,p);
print(matrank(A));
A=lift(A);
print(A, " A");
}
\end{verbatim}}
\hrulefill

{\scriptsize \begin{verbatim}
gp > test()
x^4 + 87*x^2 + 261 f
1 bnfcertify
[10, 2] k.clgp
3 d(G~)
29
1
-1/15*x^2 - 2/5 e1
1
-1/15*x^3 - 27/5*x gen Gal(k/Q)
16.15549442140351209375213147462098866888536048493436651304116601004993788848297
60786299013450395064532561144653818372818230658642059446668808596501846267735860
15167578555644043707970039157344048871166 > 0
1
-2/15*x^3 - 1/15*x^2 - 44/5*x - 17/5 gen P1
1
[5, [-1, -1, -2, -2]~, 1, 1] P1
[5, [0, 1, 2, -2]~, 1, 1] P2
[5, [-1, -1, 2, 2]~, 1, 1] P3
[5, [0, 1, -2, 2]~, 1, 1] P4
[[0, 0]~, [5, 0, 0, 0]~] P1P2P3P4
4 h
[1, 1]~ A0 in k.clgp
x - 9 gen A0
10 q0
-41/15*x^3 + 322/15*x^2 - 3627/5*x + 2954/5 e0
[3, 3, 4, 0; 3, 2, 3, 2] Z
[4, 2, 6, 8] l_{i,0}
4892/15*x^3 + 16957/15*x^2 + 91424/5*x + 40034/5 w1*23^4
-19/15*x^3 - 194/15*x^2 - 673/5*x - 5758/5 w2*23^2
114234/5*x^3 - 1468532/15*x^2 + 7497169/5*x - 51582334/5 w3*23^6
141541772/15*x^3 + 221625902/15*x^2 + 3865644704/5*x + 4691851209/5 w4*23^8
[0, 3, 2, 1; 1, 0, 4, 4; 1, 4, 0, 4; 1, 0, 2, 0] l_{i,j} mod p
[4, 0; 4, 3] B^-1
3
[4, 0, 4; 0, 2, 4; 3, 3, 0; 2, 1, 0] A
\end{verbatim}}
\hrulefill

\subsection{Examples \ref{ex:3C6t0} and \ref{ex:7C6t0}.}

{\scriptsize \begin{verbatim}
test(case)={\\ (p,m,n)=(3,61,2) if 0, (p,m,n)=(7,19,3) if 1
if(case==0,p=3;m=61;n=2,p=7;m=19;n=3);
print([p,m,-n]);
cubepol=polsubcyclo(m,3);
f=polcompositum(cubepol,x^2+n)[1];
print(f, " f");
k=bnfinit(f);
print(bnfcertify(k), " bnfcertify");
print(k.clgp[2], " k.clgp");
h=k.clgp[1]; print(h, " h");
print([subcycloiwasawa(x^2+n,p)[1]], " quad lambda");
R=bnrinit(k,p^2);
print(matreduce(lift(R.clgp[2]*Mod(1,p)))[1,2], " d(G~)");

if(case==0,delta=11*241,delta=7*1699);

G=nfgaloisconj(k);
if(case==0,g=G[3],g=G[2]);
print(lift(nfgaloisapply(k,g,g)) != x);
print(lift(nfgaloisapply(k,nfgaloisapply(k,g,g),g)) != x);
print(g*delta," /delta gen Gal(k/Q)");

P1=idealprimedec(k,p)[1];
print(k.zk*P1[2], " gen P1");

w1=Mod(k.zk*bnfisprincipal(k,idealpow(k,P1,h))[2],f);
if(case==0,print(lift(w1*delta)," /delta w1"),print(lift(w1), " w1"));

P=listcreate();
w=listcreate();
for(i=1,6,
listput(P,P1);P1=nfgaloisapply(k,g,P1);
listput(w,w1);w1=nfgaloisapply(k,g,w1)
);
P=Vec(P);w=Vec(w);

e=k.fu;
print(lift(e[1])*delta," /delta e1");
print(lift(e[2])*delta," /delta e2");

Z=matrix(2,6);
for(i=1,2, for(j=1,6,
for(z=0,p-1,
if(idealval(k,e[i]^(p-1)-(1+p)^((p-1)*z),P[j]) > 1, Z[i,j]=z;break)
)
));

print(Z, " Z");

l=matrix(6,6);
for(i=1,length(P),for(j=1,length(P), if(j==i,l[i,j]=0,
for(z=0,p-1,
if(idealval(k, w[i]^(p-1)-(1+p)^((p-1)*z),P[j]) > 1, l[i,j]=z;break)
)
)));

print(l, " l_{i,j} mod p");
print( ((l[1,2]-l[1,5])^3+(l[1,3]-l[1,6])^3)%p," (l12-l15)^3+(l13-l16)^3 mod p");
m1=l[1,2]+l[1,4]+l[1,6];
print(m1%p," m1 mod p");
m2=(l[1,2]+2*l[1,3]+2*l[1,4]+l[1,5])*(l[1,3]+2*l[1,4]+2*l[1,5]+l[1,6]) + (l[1,2]+l[1,3]-l[1,5]-l[1,6])^2;
print(m2%p," m2 mod p");

B=Z[,5..6];
d1=matdet(B);
b=B^-1;
c=matrix(6,4);
for(j=5,6,for(i=1,4,
  c[j,i]=-b[j-4,1]*Z[1,i]-b[j-4,2]*Z[2,i]
));
L=matrix(6,4);
for(j=1,6,for(i=1,4,
  L[j,i]=(l[j,i]+c[5,i]*l[j,5]+c[6,i]*l[j,6])/h
));
A=[
 L[1,2], L[1,3], L[1,4],      0,      0,      0;
-L[2,1],      0,      0, L[2,3], L[2,4],      0;
      0,-L[3,1],      0,-L[3,2],      0, L[3,4];
      0,      0,-L[4,1],      0,-L[4,2],-L[4,3];
      0,      0,      0,      0,      0,      0;
      0,      0,      0,      0,      0,      0
];
for(j=5,6,
A[j,1]=c[j,1]*L[j,2]-c[j,2]*L[j,1];
A[j,2]=c[j,1]*L[j,3]-c[j,3]*L[j,1];
A[j,3]=c[j,1]*L[j,4]-c[j,4]*L[j,1];
A[j,4]=c[j,2]*L[j,3]-c[j,3]*L[j,2];
A[j,5]=c[j,2]*L[j,4]-c[j,4]*L[j,2];
A[j,6]=c[j,3]*L[j,4]-c[j,4]*L[j,3];
);

print(lift(matdet(A*Mod(1,p))), " det A");
print(lift(A*Mod(1,p)));
Aprime=[
1 , -1 ,  0 ,  1 ,  0 ,  0 ;
1 ,  0 , -1 ,  0 ,  1 ,  0 ;
0 ,  1 , -1 ,  0 ,  0 ,  1 
];
Aplus=matconcat([A;Aprime])*Mod(1,p);
print(matrank(Aplus));

l = l*Mod(1,p);
if((m1*m2)%p != 0, 
aa = l[1,2] + l[1,3] + l[1,4];
bb = -(l[1,3] + l[1,4] + l[1,5]);
cc = l[1,4] + l[1,5] + l[1,6];
a =  h*( aa^2 - bb*cc )/(m1*m2);
b = -h*( cc^2 - aa*bb )/(m1*m2);
c =  h*( bb^2 - aa*cc )/(m1*m2);
print([lift(a),lift(b),lift(c)]," [a,b,c] mod p"); 
);
}
\end{verbatim}}
\hrulefill

{\scriptsize \begin{verbatim}
gp > test()
[3, 61, -2]
x^6 - 2*x^5 - 33*x^4 + 50*x^3 + 398*x^2 - 396*x + 1089 f
1 bnfcertify
[43] k.clgp
43 h
[1] quad lambda
4 d(G~)
1
1
-4*x^5 - 25*x^4 + 440/3*x^3 + 1010*x^2 - 10265/3*x - 7425 /delta gen Gal(k/Q)
x - 2 gen P1
29881/3*x^5 - 93587/3*x^4 - 516835/3*x^3 + 958130/3*x^2 - 953187*x + 1625195 /de
lta w1
-232/9*x^5 + 1201/9*x^4 + 4972/9*x^3 - 36856/9*x^2 + 4852*x - 8085 /delta e1
-46/9*x^5 - 1613/9*x^4 + 3454/9*x^3 + 59333/9*x^2 - 5894*x - 92994 /delta e2
[2, 0, 1, 2, 0, 1; 1, 2, 0, 1, 2, 0] Z
[0, 2, 0, 0, 1, 2; 2, 0, 2, 0, 0, 1; 1, 2, 0, 2, 0, 0; 0, 1, 2, 0, 2, 0; 0, 0, 1
, 2, 0, 2; 2, 0, 0, 1, 2, 0] l_{i,j} mod p
2 (l12-l15)^3+(l13-l16)^3 mod p
1 m1 mod p
1 m2 mod p
2 det A
[1, 1, 0, 0, 0, 0; 0, 0, 0, 1, 1, 0; 0, 2, 0, 1, 0, 2; 0, 0, 1, 0, 1, 1; 2, 2, 2
, 1, 2, 1; 1, 1, 2, 1, 2, 0]
6
[1, 1, 1] [a,b,c] mod p
\end{verbatim}}
\hrulefill

{\scriptsize \begin{verbatim}
gp > test(1)
[7, 19, -3]
x^6 - 2*x^5 - 2*x^4 + 14*x^3 + 55*x^2 - 192*x + 343 f
1 bnfcertify
[3] k.clgp
3 h
[1] quad lambda
4 d(G~)
1
1
-272*x^5 - 305*x^4 + 860*x^3 + 1150*x^2 - 24925*x + 6300 /delta gen Gal(k/Q)
x - 7 gen P1
-x w1
-272*x^5 - 305*x^4 + 860*x^3 + 1150*x^2 - 13032*x + 6300 /delta e1
-162*x^5 - 444*x^4 + 862*x^3 + 2259*x^2 - 14058*x - 19684 /delta e2
[1, 1, 5, 1, 1, 5; 1, 5, 1, 1, 5, 1] Z
[0, 1, 0, 2, 6, 1; 1, 0, 1, 0, 2, 6; 6, 1, 0, 1, 0, 2; 2, 6, 1, 0, 1, 0; 0, 2, 6
, 1, 0, 1; 1, 0, 2, 6, 1, 0] l_{i,j} mod p
0 (l12-l15)^3+(l13-l16)^3 mod p
4 m1 mod p
6 m2 mod p
0 det A
[3, 2, 3, 0, 0, 0; 4, 0, 0, 3, 5, 0; 0, 2, 0, 2, 0, 1; 0, 0, 6, 0, 3, 2; 1, 4, 5
, 3, 1, 3; 2, 6, 4, 2, 5, 4]
6
[4, 0, 2] [a,b,c] mod p
\end{verbatim}}
\hrulefill

\subsection{Example \ref{ex:3C6t1}.}

{\scriptsize \begin{verbatim}
test()={
p=3;m=61;n=23;Na=41;
print([m,-n,Na]);
cubepol=polsubcyclo(m,3);
f=polcompositum(cubepol,x^2+n)[1];
print(f, " f");
k=bnfinit(f);
print(bnfcertify(k), " bnfcertify");
print(k.clgp[2], " k.clgp");
h=k.clgp[1]/(p^valuation(k.clgp[1],p)); print(h, " h");
print([subcycloiwasawa(x^2+n,p)[1]], " quad lambda");
R=bnrinit(k,p^2);
print(matreduce(lift(R.clgp[2]*Mod(1,p)))[1,2], " d(G~)");

delta=487*499;

G=nfgaloisconj(k);
g=G[3];
print(lift(nfgaloisapply(k,g,g)) != x);
print(lift(nfgaloisapply(k,nfgaloisapply(k,g,g),g)) != x);
print(g*delta," /delta gen Gal(k/Q)");

P1=idealprimedec(k,p)[1];
print(k.zk*P1[2]*delta," /delta gen P1");

P=listcreate();
for(i=1,6,
listput(P,P1);P1=nfgaloisapply(k,g,P1);
);
P=Vec(P);

A0=idealprimedec(k,Na)[1];
S=bnfisprincipal(k,A0)[1];print(S, " A0 in k.clgp");
print(k.zk*A0[2]," gen A0");

q0=1;
for(i=1,length(k.clgp[2]),
q0=lcm(k.clgp[2][i]/gcd(k.clgp[2][i],S[i]),q0)
);
e0=k.zk*bnfisprincipal(k,idealpow(k,A0,q0))[2];
print(q0, " q0");
print(e0*delta," /delta e0");

Z=matrix(3,6);

for(i=1,length(P),
for(z=0,p-1,
if(idealval(k, e0^(p-1)-(1+p)^((p-1)*z) , P[i]) > 1,
Z[1,i]=z; break
)));

e=k.fu;
print(lift(e[1])*delta," /delta e1");
print(lift(e[2])*delta," /delta e2");

for(j=1,2,for(i=1,length(P),
for(z=0,p-1,
if(idealval(k, e[j]^(p-1)-(1+p)^((p-1)*z) , P[i]) > 1,
Z[j+1,i]=z; break
))));

print(Z, " Z");

A0inv=idealinv(k,A0);

l0=listcreate();
w=listcreate();
for(i=1,length(P),Ph=idealpow(k,P[i],h);
for(z=1,k.clgp[1],
W=bnfisprincipal(k,idealmul(k,Ph,idealpow(k,A0inv,z)));
if(W[1][1]==0, 
listput(l0,z);listput(w,k.zk*W[2]);break
)));

print(Vec(l0), " l_{i,0}");

for(i=1,6,
print(lift(w[i]*delta*Na^l0[i])," /(delta*",Na,"^",l0[i],") w",i);
);

l=matrix(6,6);
for(i=1,length(P),for(j=1,length(P), if(j==i,l[i,j]=0,
for(z=0,p-1,
if(idealval(k, w[i]^(p-1)-(1+p)^((p-1)*z),P[j]) > 1, l[i,j]=z;break)
)
)));

print(l, " l_{i,j}");

l=l*Mod(1,p);
Z=Z*Mod(1,p);

B=Z[,4..6];
print(lift(B^-1), " B^-1");
b=matconcat([matrix(3,3);B^-1]);
c=matrix(6,3);
for(j=4,6, for(i=1,3, 
  c[j,i]=-b[j,1]*Z[1,i]-b[j,2]*Z[2,i]-b[j,3]*Z[3,i]
));
L=matrix(6,3);
for(j=1,6, for(i=1,3,
  L[j,i]=(l[j,i]+c[4,i]*l[j,4]+c[5,i]*l[j,5]+c[6,i]*l[j,6])/h
));

A=matconcat([
[-l0[1]/h,0,0;
 0,-l0[2]/h,0;
 0,0,-l0[3]/h],
[ L[1,2], L[1,3],      0;
 -L[2,1],      0, L[2,3];
       0,-L[3,1],-L[3,2]];
matrix(3,3),matrix(3,3)]);
for(i=4,6,
A[i,1]=-c[i,1]*l0[i]/h;
A[i,2]=-c[i,2]*l0[i]/h;
A[i,3]=-c[i,3]*l0[i]/h;
A[i,4]=c[i,1]*L[i,2]-c[i,2]*L[i,1];
A[i,5]=c[i,1]*L[i,3]-c[i,3]*L[i,1];
A[i,6]=c[i,2]*L[i,3]-c[i,3]*L[i,2];
);
A=A*Mod(1,p);
print(lift(matdet(A)), " det A");
print(lift(A), " A");
}
\end{verbatim}}
\hrulefill

{\scriptsize \begin{verbatim}
gp > test()
[61, -23, 41]
x^6 - 2*x^5 + 30*x^4 - 34*x^3 + 2015*x^2 - 1740*x + 43551 f
1 bnfcertify
[36, 4] k.clgp
16 h
[1] quad lambda
4 d(G~)
1
1
-388/9*x^5 - 5855/9*x^4 - 8360/9*x^3 - 24460/9*x^2 - 800795/3*x - 623540 /delta 
gen Gal(k/Q)
-469/2*x^5 + 971/3*x^4 + 9895/3*x^3 + 171851/6*x^2 - 981047/6*x + 221379/2 /delt
a gen P1
[20, 1]~ A0 in k.clgp
x - 5 gen A0
36 q0
1888524068369/18*x^5 + 6943574603762/9*x^4 - 28025383797112/9*x^3 - 897873595515
1/18*x^2 + 929227301965079/6*x + 4220669428821629/2 /delta e0
-354*x^5 - 1584*x^4 - 12638*x^3 + 231*x^2 - 176362*x + 1138112 /delta e1
-320*x^5 + 2687*x^4 - 16916*x^3 + 24922*x^2 - 137456*x - 860384 /delta e2
[2, 0, 2, 0, 2, 0; 2, 0, 1, 2, 0, 1; 1, 2, 0, 1, 2, 0] Z
[8, 4, 2, 1, 5, 7] l_{i,0}
-459549066139891/6*x^5 - 704093863348502/3*x^4 - 535153539380573*x^3 - 624382779
1610063/6*x^2 - 808834974927696857/6*x - 1354543489312812447/2 /(delta*41^8) w1
-1138617949/18*x^5 + 9263007134/9*x^4 + 21513813602/9*x^3 + 904069241621/18*x^2 
+ 561744316739/6*x + 4365781564589/2 /(delta*41^4) w2
14439515/162*x^5 + 32950277/54*x^4 + 90525712/27*x^3 + 572674304/81*x^2 + 121936
0249/27*x + 2077748773/18 /(delta*41^2) w3
1649801/729*x^5 - 181297/54*x^4 - 893647/243*x^3 + 287666645/1458*x^2 + 12915881
33/486*x - 8536877/9 /(delta*41^1) w4
-322820922181/1458*x^5 + 76556486800/81*x^4 + 5754643430519/486*x^3 - 2000338322
3033/1458*x^2 - 222480668144683/243*x + 183586414913389/27 /(delta*41^5) w5
-10129957962787/6*x^5 + 99318329682100/9*x^4 - 84269452005608/3*x^3 + 4881925038
13457/6*x^2 - 49749237476440823/18*x + 119523731360647853/6 /(delta*41^7) w6
[0, 0, 2, 2, 1, 0; 1, 0, 2, 0, 0, 0; 0, 0, 0, 2, 1, 0; 0, 0, 0, 0, 0, 2; 0, 0, 0
, 1, 0, 1; 0, 0, 0, 0, 1, 0] l_{i,j}
[2, 0, 1; 2, 0, 0; 2, 1, 1] B^-1
2 det A
[1, 0, 0, 2, 2, 0; 0, 2, 0, 2, 0, 2; 0, 0, 1, 0, 2, 1; 2, 2, 1, 1, 0, 1; 2, 0, 2
, 1, 0, 2; 1, 2, 2, 1, 2, 2] A
\end{verbatim}}
\hrulefill

\end{document}